\documentclass[twocolumn]{autart}

\usepackage{graphicx}
\usepackage{epsfig} 
\usepackage{mathptmx} 
\usepackage{times} 
\usepackage{amsmath} 
\usepackage{amssymb}  
\usepackage{amsfonts}
\usepackage{euscript}

\newcommand{\iterc}{{\kappa}}

\newcommand{\frakP}{\mathfrak{P}}
\newcommand{\frakD}{\mathfrak{D}}

\newcommand{\sfp}{\mathsf{p}}
\newcommand{\sfd}{\mathsf{d}}

\newcommand{\bx}{{\overline x}}

\newcommand{\bw}{{\overline w}}
\newcommand{\bv}{{\overline v}}

\usepackage{color}
\usepackage[dvipsnames]{xcolor}
\usepackage{multicol}
\usepackage{caption}
\usepackage{subcaption}
\usepackage{enumitem}

\definecolor{amethyst}{rgb}{1, 0, 1}
\definecolor{blue-violet}{rgb}{0.54, 0.17, 0.89}
\definecolor{brightturquoise}{rgb}{0.03, 0.91, 0.87}

\newtheorem{theorem}{Theorem}
\newtheorem{definition}[theorem]{Definition}

\newcommand{\mB}[1]{{\mathbb{#1}}}

\newcommand{\R}{{\mathbb R}}

\newcommand{\G}{{\mathcal{N}}}

\newcommand{\bB}{\mathbb{B}}

\newcommand{\Bp}{{\bf p}}

\newcommand{\prox}{\mbox{prox}}
\newcommand{\proj}{\mbox{proj}}

\newcommand{\Sc}[1]{{\mathcal{#1}}}

\newcommand{\calC}{\Sc{C}}
\newcommand{\calK}{\Sc{K}}

\newcommand{\calE}{\Sc{E}}

\newcommand{\calV}{\Sc{V}}

\newcommand{\del}{\delta}
\newcommand{\eps}{\epsilon}
\newcommand{\lam}{\lambda}

\newcommand{\sig}{\sigma}
\newcommand{\Gam}{\Gamma}

\newcommand{\map}[3]{#1:#2\rightarrow #3}
\newcommand{\abs}[1]{\vert #1\vert}

\newcommand{\norm}[1]{\left\Vert #1\right\Vert}
\newcommand{\tnorm}[1]{\left\Vert #1\right\Vert_2}
\newcommand{\onorm}[1]{\left\Vert #1\right\Vert_1}
\newcommand{\inorm}[1]{\left\Vert #1\right\Vert_\infty}
\newcommand{\dnorm}[1]{\left\Vert #1\right\Vert_\circ}
\newcommand{\polar}[1]{#1^\circ}

\newcommand{\argmin}{\mathop{\mathrm{argmin}}}

\newcommand{\support}[2]{\sig_{#2}\left(#1\right)}

\newcommand{\indicator}[2]{\delta_{#2}\left(#1 \right)}
\newcommand{\ip}[2]{\left\langle #1,\, #2\right\rangle}
\newcommand{\dom}[1]{\mathrm{dom}\left(#1\right)}
\newcommand{\lev}[2]{\mathrm{lev}_{#1}\left(#2\right)}
\newcommand{\epi}[1]{\mathrm{epi}\left(#1\right)}

\newcommand{\Skp}{\mathcal{S}^k_+}

\newcommand{\ncone}[2]{N\left(#1\left|\, #2\right.\right)}

\newcommand{\Rn}{\mB{R}^n}

\newcommand{\Nul}[1]{\mathrm{Nul}\left( #1\right)}

\newcommand{\set}[2]{\left\{#1\,\left|\, #2\right.\right\}}
\newcommand{\cl}[1]{\mathrm{cl}\left(#1\right)}
\newcommand{\intr}[1]{\mathrm{intr}\left(#1\right)}
\newcommand{\bdry}[1]{\mathrm{bdry}\left(#1\right)}
\newcommand{\rbdry}[1]{\mathrm{rbdry}\left(#1\right)}
\newcommand{\ri}[1]{\mathrm{ri}\left(#1\right)}
\newcommand{\conv}[1]{\mathrm{conv}\left(#1\right)}

\newcommand{\aff}[1]{\mathrm{aff}\left(#1\right)}

\usepackage{algorithm}
\usepackage{algorithmic}

\usepackage{pgfplots}

\begin{document}

\begin{frontmatter}

\thanks[footnoteinfo]{Corresponding author Gianluigi Pillonetto Ph. +390498277607. 
This research was supported  in part by the National 
Science Foundation grant no. DMS-1514559, by the Washington Research 
Foundation Data Science Professorship, 
by the MIUR FIRB project RBFR12M3AC-Learning
meets time: a new computational approach to learning in dynamic
systems, by the Progetto di Ateneo CPDA147754/14-New statistical learning approach for multi-agents adaptive estimation and coverage control  
as well as by the Linnaeus Center CADICS, funded by the
Swedish Research Council, and the ERC advanced grant LEARN, no 287381, funded by the European Research Council.}

\title{Generalized Kalman Smoothing: Modeling and Algorithms}

\author[First]{Aleksandr Aravkin}
\author[Second]{James V. Burke}
\author[Third]{Lennart Ljung}
\author[Fourth]{Aurelie Lozano}
\author[Fifth]{Gianluigi Pillonetto}
\address[First]{Department of Applied Mathematics, University of Washington, USA (e-mail: saravkin@uw.edu)}
\address[Second]{Department of Mathematics, University of Washington, Seattle, USA (e-mail: burke@math.washington.edu)}
\address[Third]{Division of Automatic Control, Link\"oping University, Link\"oping, Sweden (e-mail: ljung@isy.liu.se)}
\address[Fourth]{IBM T.J. Watson Research Center Yorktown Heights, NY, USA (e-mail: aclozano@us.ibm.com)}
\address[Fifth]{Department of Information  Engineering, University of Padova, Padova, Italy (e-mail: giapi@dei.unipd.it)}

\begin{abstract}
State-space smoothing has found many applications in science and
engineering.  Under linear and Gaussian assumptions, smoothed estimates can
be obtained using efficient recursions, for example Rauch-Tung-Striebel and
Mayne-Fraser algorithms. Such schemes are equivalent to linear
algebraic techniques that minimize a
convex quadratic objective function with structure induced by the dynamic
model.\\
\noindent These classical formulations fall short in many important circumstances.
For instance, smoothers obtained using quadratic penalties can fail when
outliers are present in the data, and cannot track impulsive inputs and
abrupt state changes. Motivated by these shortcomings, generalized Kalman
smoothing formulations have been proposed in the last few years, replacing quadratic models with
more suitable, often nonsmooth, convex functions.  In contrast to classical
models, these general estimators require use of iterated algorithms, and
these have received increased attention from control, signal
processing, machine learning, and optimization communities.\\
\noindent In this survey we show that the optimization viewpoint provides the control and signal processing
community great freedom in the development of novel modeling and inference
frameworks for dynamical systems. We discuss general statistical models for
dynamic systems, making full use of nonsmooth convex penalties and
constraints, and providing links to important models in signal
processing and machine learning. We also survey optimization techniques for
these formulations, paying close attention to dynamic problem
structure. Modeling
concepts and algorithms are illustrated with numerical examples.
\end{abstract}

\end{frontmatter}

\section{Introduction}

The linear state space model
\begin{subequations}  \label{eq:Lin}
\begin{align}
 x_{t+1} &=A_t x_t+B_t u_t +v_t\\
y_t &=C_t x_t + e_t
\end{align}
\end{subequations}
is the bread and butter for analysis and design in discrete time systems, control and signal processing \cite{kalman,KalBuc}. Applications areas are numerous, including navigation, tracking, healthcare and finance, to name a few.

For a system model, $y_t \in \R^{m}$ and $u_t \in \R^{p}$ are, respectively, the output and 
input evaluated at the time instant $t$. The dimensions $m$ and $p$ may depend on $t$, but we
treat them as fixed to simplify the exposition. 
In signal models, the input $u_t$ may be absent. 
The state vectors $x_t \in \R^n$ are the variables of interest; 
$A_t$ encodes the process transition, to the extent that it is known to the modeler,  
$C_t$ is the observation model, and $B_t$ describes the effect of the input on the transition.  
 The \emph{process disturbance} $v_t$ models stochastic deviations from the linear model $A_t$, while  $e_t$  model \emph{measurement errors}. We consider the {\it state estimation problem}, where the goal is to infer the values of $x_t$ from the input-output measurements. Given measurements
$$\mathcal{Z}^N_0:=\{u_0,y_1,u_1,y_2,\ldots,y_N,u_N\},$$ 
we are interested in obtaining an estimate $\hat{x}^N_t$ of $x_t$.  
If $N>t$ this is called a \emph{smoothing} problem, if $N=t$ it is a \emph{filtering} problem, and if $N<t$ it is a \emph{prediction} problem.

How well the state estimate fits the true state depends upon the choice of 
models for the stochastic term $v_t$, error term $e_t$, and possibly on the initial distribution of $x_0$.
While $u_t$ is usually a known deterministic sequence, the observations $y_t$ 
and states $x_t$ are stochastic processes. 
We can consider using several estimators $\hat{x}^N_t$ of the state sequence $\{x_t\}$ (all functions of $\mathcal{Z}^N_0$):
\begin{subequations}
\begin{align}
  \label{eq:cm}
E(x_t | \mathcal Z^N_0) \quad  &\mbox{conditional mean} \\  \label{eq:cm2}
 \max_{x_t} \Bp(x_t \big|  \mathcal Z^N_0)  \quad &\mbox{maximum {\it a posteriori} (MAP)} \\  \label{eq:cm3}
\nonumber  \min_{\hat{x_t}} E(\|x_t - \hat{x}_t\|^2) \quad  & \mbox{minimum expected} \\
  &  \mbox{mean square error (MSE)} \\  
 \min_{\hat{x_t} \in \mathrm{span}\left(\mathcal Z^N_0\right)} E(\|x_t - \hat{x}_t\|^2) \  
   & \mbox{minimum linear expected MSE}  \label{eq:cm4} 
\end{align}
\end{subequations}
When $e_t,v_t$ and the initial state $x_0$ are jointly Gaussian, all the four estimators coincide.
In the general setting, the estimators (\ref{eq:cm}) and (\ref{eq:cm3}) are the same.
Indeed, the conditional mean represents the minimum variance estimate.
In the general (non-Gaussian) case, computing~\eqref{eq:cm} may be difficult,  
while the MAP~\eqref{eq:cm2} estimator can be computed efficiently using optimization techniques 
for a range of disturbance and error distributions. 

Most models assume known means and variances 
for $v_t, e_t,$ and $x_0$. In the classic settings, these distributions are Gaussian:  
\begin{equation}
 \label{eq:wgn}
\begin{aligned}
e_t  & \sim \mathcal{N}(0,R_t) \\ 
v_t & \sim \mathcal{N}(0,Q_t)\\
x_0 & \sim \mathcal{N}(\mu,\Pi)
\end{aligned}, \qquad \text{all variables are mutually independent.}
\end{equation}
Under this assumption, all the
$y_t$ and $x_t$ become  jointly Gaussian stochastic processes, which implies that the conditional mean (\ref{eq:cm}) becomes a linear function of the data $\mathcal{Z}^N_0$. This is a general property of Gaussian variables. Many explicit expressions and recursions for this linear filter  have been derived in the literature, 
some of which are discussed in this article.  We also consider a far more general setting, where the 
distributions in~\eqref{eq:wgn} can be selected from a range of densities, and discuss applications and general inference techniques.

We now make explicit the connection between {\it conditional mean}~\eqref{eq:cm} and {\it maximum likelihood}~\eqref{eq:cm2} in the Gaussian case. 
By Bayes' theorem and the independence assumptions (\ref{eq:wgn}),
the posterior  of the state sequence $\{x_t\}_{t=0}^N$ given the measurement sequence $\{y_t\}_{t=1}^N$ is 
\begin{eqnarray}\nonumber 
&& \Bp\left(\{x_t\} \big| \{y_t\}\right) = \frac{\Bp\left(\{y_t\}\big|\{x_t\}\right)\Bp\left(\{x_t\}\right)}{\Bp\left(\{y_t\} \right)} \\ \label{Bayes} 
\qquad &&= \frac{\Bp\left(x_0 \right)  \prod_{t=1}^N \Bp\left( y_t \big| x_t  \right)  \prod_{t=0}^{N-1} \Bp\left( x_{t+1} \big| x_t  \right)}{\Bp\left(\{y_t\} \right)} \\ \nonumber 
\qquad && \propto  \Bp\left(x_0 \right) \prod_{t=1}^N \Bp_{e_t}  \left( y_t - C_t x_t  \right) \prod_{t=0}^{N-1} \Bp_{v_t} \left( x_{t+1} -A_t x_t -B_t u_t  \right),
\end{eqnarray}
where we use $\Bp_{e_t}$ and $\Bp_{v_t}$ to denote the densities corresponding to $e_t$ and $v_t$. 
Under Gaussian assumptions~\eqref{eq:wgn}, and ignoring the normalizing constant, the posterior is given by 

\begin{equation}
\label{MAP}
\begin{aligned}
e^{-\frac{1}{2}\left\|\Pi^{-1/2} (x_0 - \mu) \right\|^2} \prod_{t=0}^{N-1}  e^{-\frac{1}{2}\left\| Q_t^{-1/2}(x_{t+1} -A_t x_t -B_t u_t)\right\|^2} \\ 
\times  \prod_{t=1}^N e^{-\frac{1}{2}\left\|R_t^{-1/2}(y_t-C_t x_t)\right\|^2}.
\end{aligned}
\end{equation}
Note that state increments and measurement residuals appear explicitly in~\eqref{MAP}.
Maximizing~\eqref{MAP} is equivalent to minimizing its negative log: 
\begin{equation}
  \label{eq:opt}
\begin{aligned} 
\min_{x_0,\ldots,x_N}  \left\|\Pi^{-1/2} (x_0 - \mu) \right\|^2+\sum_{t=1}^N\left\|R_t^{-1/2}(y_t-C_t x_t)\right\|^2\\
+ \sum_{t=0}^{N-1}\left\| Q_t^{-1/2}(x_{t+1} -A_t x_t -B_t u_t)\right\|^2.
\end{aligned}
\end{equation}

 More general cases of correlated noise and singular covariance matrices are discussed in the Appendix. 
This result is also shown in e.g.~\cite{Bell1994} and \cite[Sec. 3.5, 10.6]{Kailath_Sayed_Hassibi:00} 
using a least squares argument.
The solution can be derived using various structure-exploiting linear recursions. 
For instance, the Rauch-Tung-Striebel (RTS) scheme derived in \cite{RTS} computes the state estimates
by forward-backward recursions, (see also  \cite{Ansley1982}
for a simple derivation through projections onto spaces spanned by suitable random variables.)
The Mayne-Fraser (MF) algorithm uses a two-filter formula to compute the smoothed estimate as a 
linear combination of forward and backward Kalman filtering estimates  \cite{Mayne1966,FraserPotter1969}.
A third scheme based on reverse recursion appears in \cite{Mayne1966}
under the name of Algorithm A. The relationships between these schemes, 
and their derivations from different perspectives are studied in \cite{Ljung1976,aravkin2013kalman}. 
Computational details for RTS and MF are presented in Section \ref{L2opt}. 

The {\it maximum a posteriori} (MAP) viewpoint~\eqref{eq:opt} easily generalizes to new settings. 
Assume, for example, that the noises $e_t$ and $v_t$ are non-Gaussian, but rather have continuous probability densities 
defined by functions $V_t(\cdot)$ and $J_t(\cdot)$ as follows 
\begin{equation}\label{Generalpdf}
\Bp_{e_t}(e) \propto  \exp \left(- V_t\left(R_t^{-1/2} e\right)\right), \   \Bp_{v_t}(v) \propto \exp \left(- J_t\left(Q_t^{-1/2} v\right)\right).
\end{equation}
From (\ref{Bayes}), we obtain that the analogous MAP estimation problem for~\eqref{eq:opt} replaces all least squares 
$\|R_t^{-1/2}(y_t-C_t x_t)\|^2$ and $\|Q_t^{-1/2}(x_{t+1} - A_t x_t -B_t u_t)\|^2$ with more general terms $V_t \left(R_t^{-1/2}(y_t-C_t x_t)\right)$ and $J_t\left(Q_t^{-1/2}(x_{t+1} - A_t x_t -B_t u_t)\right)$, leading to  
\begin{equation}
  \label{eq:optVJ}
\begin{aligned}
\min_{x_0,\ldots,x_N}   &   -\log \Bp(x_0) +  \sum_{t=1}^NV_t \left(R_t^{-1/2}(y_t-C_t x_t)\right)\\
& + \sum_{t=0}^{N-1} J_t\left(Q_t^{-1/2}(x_{t+1} - A_t x_t -B_t u_t)\right).
\end{aligned}
\end{equation}
The initial distribution for $x_0$ can be non-Gaussian, and is specified by $\Bp(x_0)$.
An algorithm to solve~\eqref{eq:optVJ} is then required. 
In this paper, we will discuss general modeling of error distributions 
$\Bp_{e_t}$ and $\Bp_{v_t}$ in~\eqref{Generalpdf}, as well as tractable algorithms 
for the solutions of these formulations.

Classic {\it{Kalman}} filters, predictors and smoothers have been enormously successful, and the literature detailing their properties and applications is rich and pervasive.
Even if Gaussian assumptions~\eqref{eq:wgn} are violated, but the $v_t$, $e_t$
are still white with covariances $Q_t$ and $R_t$, problem (\ref{eq:opt}) gives the best {\it linear} estimate, 
i.e. among all linear functions of the data $\mathcal{Z}^N_0$, the Kalman smoother residual has the smallest variance. However, this does not ensure successful performance, giving 
strong motivation to consider extensions to the Gaussian framework!  
%
For instance, impulsive disturbances often occur in process models, including
target tracking, where one has to deal with force disturbances describing maneuvers for the tracked object,
fault detection/isolation, where impulses model additive faults, and load disturbances. 
Unfortunately, smoothers that use the quadratic penalty on the state increments are not able to follow 
fast jumps in the state dynamics \cite{Ohlsson2011}. 
This problem is also relevant in the context of 
identification of switched linear regression models where the 
system states can be seen as time varying parameters which can be subject to abrupt changes 
\cite{Ohlsson2013,Nied2013}. In addition, constraints on the states arise naturally in many settings,  
and estimation can be improved by taking these constraints into account. 
Finally, estimates corresponding to quadratic losses applied to data misfit residuals are vulnerable to outliers, i.e. to unexpected deviations of the noise errors from Gaussian assumptions.
In these cases, a Gaussian model for $e$ gives poor estimates. 
Two examples are described below, the first focusing on impulsive disturbances, and second on measurement outliers. 

\subsection{DC motor example}\label{MotEx}

A DC motor can be modeled as a dynamic system, where the input is applied torque
while the output is the angle of the motor shaft, see also pp. 95-97 in \cite{Ljung:99}.
The state comprises angular velocity and angle of the motor shaft, 
and with system parameters and discretization 
as in Section 8 of \cite{Ohlsson2011}, we have the following discrete-time model: 
\begin{equation}
  \label{DCmotor}
\begin{aligned}
x_{t+1} &= \left(\begin{array}{cc} 0.7 & 0 \\0.084 & 1\end{array}\right) x_t
+\left(\begin{array}{c}11.81 \\0.62 \end{array}\right) (u_t + d_t)\\
y_t &= \left(\begin{array}{cc}0 & 1\end{array}\right) x_t + e_t
\end{aligned}
\end{equation}
where $d_t$ denotes a disturbance process while the measurements $y_t$ are noisy samples
of the angle of the motor shaft.
\begin{figure*}
  \begin{center}
   \begin{tabular}{cccc}
\hspace{.1in}
 { \includegraphics[scale=0.46]{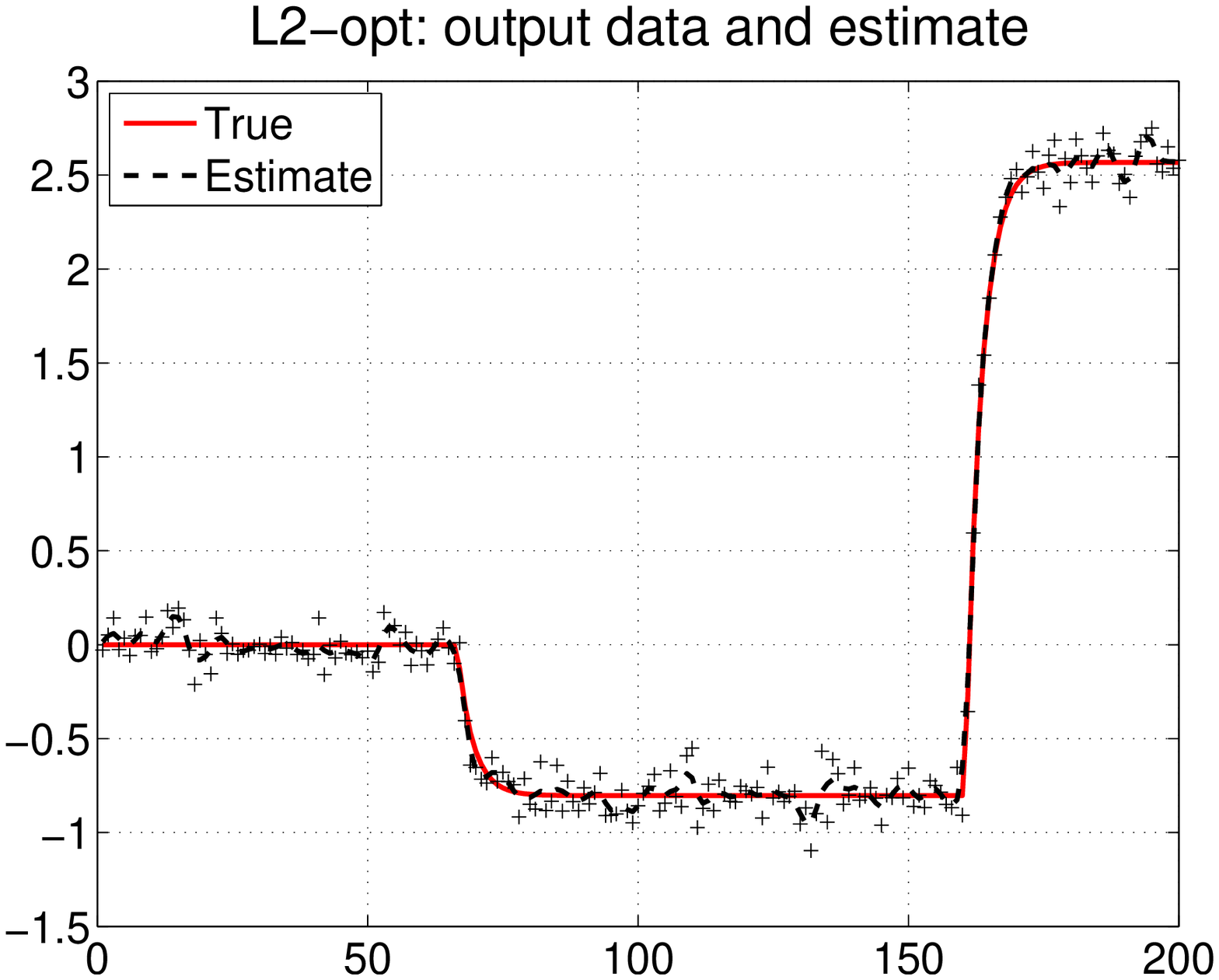}} 
\hspace{.1in}
 { \includegraphics[scale=0.46]{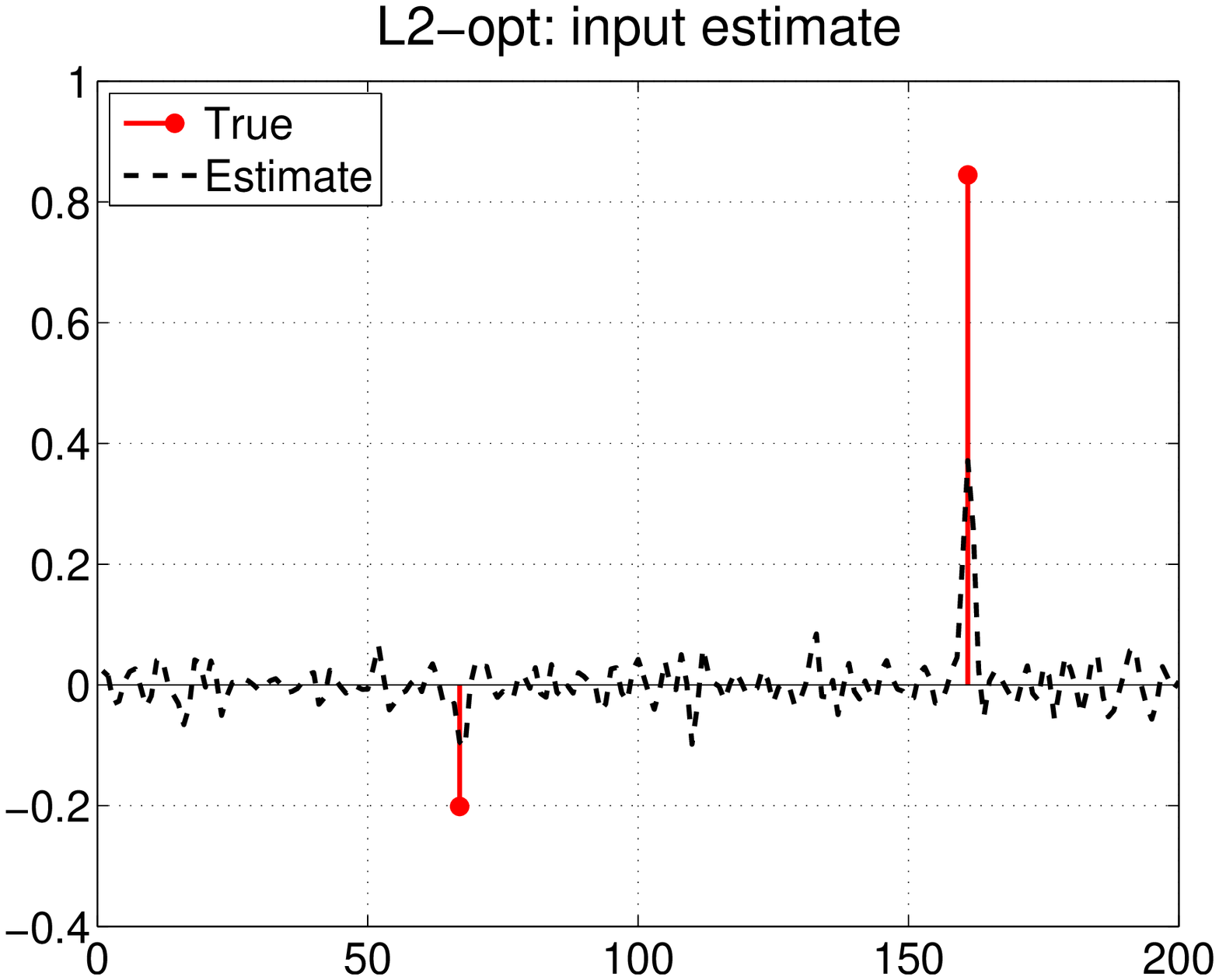}} 
    \end{tabular}
 \caption{{\bf{DC motor and impulsive disturbances}}. {\it{Left:}} noiseless output (solid line), measurements  ($+$) and output reconstruction by the optimal linear smoother L$_2$-opt (dashed line). 
 {\it{Right:}} impulsive disturbance and reconstruction by L$_2$-opt (dashed line). } 
    \label{Fig12Imp}
     \end{center}
\end{figure*}

{\bf{Impulsive inputs:}} In the DC system design, the
disturbance torque acting on the motor shaft plays an important role and
an accurate reconstruction of $d_t$ 
can greatly improve model robustness with respect to  load variations.
Since the non observable input is often impulsive, 
we model the $d_t$ as independent random variables such that
$$
d_t  = 
\left\{ \begin{array}{cl}
    0 & \mbox{with probability} ~ 1-\alpha \\
    \G(0,1) & \mbox{with probability} ~ \alpha 
\end{array} \right.
$$
According to (\ref{eq:Lin}), this corresponds to a zero-mean (non-Gaussian) noise $v_t$,
with covariance
\(
Q_t=\alpha 
\binom{11.81}{0.62} (11.81,\, 0.62).
\)
We consider the problem of reconstructing $d_t$ from noisy output samples
generated under the assumptions
$$
x_0 =  \left(\begin{array}{c}0 \\0\end{array}\right),   \ u_t=0, \quad \alpha=0.01, \quad e_t \sim \G(0,0.1^2).
$$
An instance of the problem is shown in Fig. \ref{Fig12Imp}.
The left panel displays the noiseless output (solid line) and the 
measurements ($+$). The right panel displays the $d_t$ (solid line)
and their estimates (dashed line) obtained by 
the Kalman smoother\footnote{
Note that the covariance matrices $Q_t$ are singular. 
In this case, the smoothed estimates have been computed using the RTS 
scheme \cite{RTS}, as e.g. described in Section 2.C of \cite{Kitagawa1985}, where invertibility 
of the transition covariance matrices are not required. This scheme provides the solution of the 
generalized Kalman smoothing objective (\ref{eq:optPseudo2}), and is explained in the Appendix.
} 
and given by
\(
\hat{d}^N_t = \left(1/11.81 \ \ 0 \right) \left( \hat{x}^N_{t+1}  -A_t \hat{x}^N_{t+1} \right).
\)

This estimator, denoted L$_2$-opt, uses only information on the means and covariances of the noises.
It solves problem (\ref{eq:cm4})  and, hence, corresponds to the best linear estimator. 
However, it is apparent that the 
disturbance reconstruction is not satisfactory. The smoother estimates of the impulses
are poor, and the largest peak, centered at $t=161$, is highly underestimated.  
\begin{figure*}
  \begin{center}
   \begin{tabular}{cccc}
\hspace{.1in}
 { \includegraphics[scale=0.46]{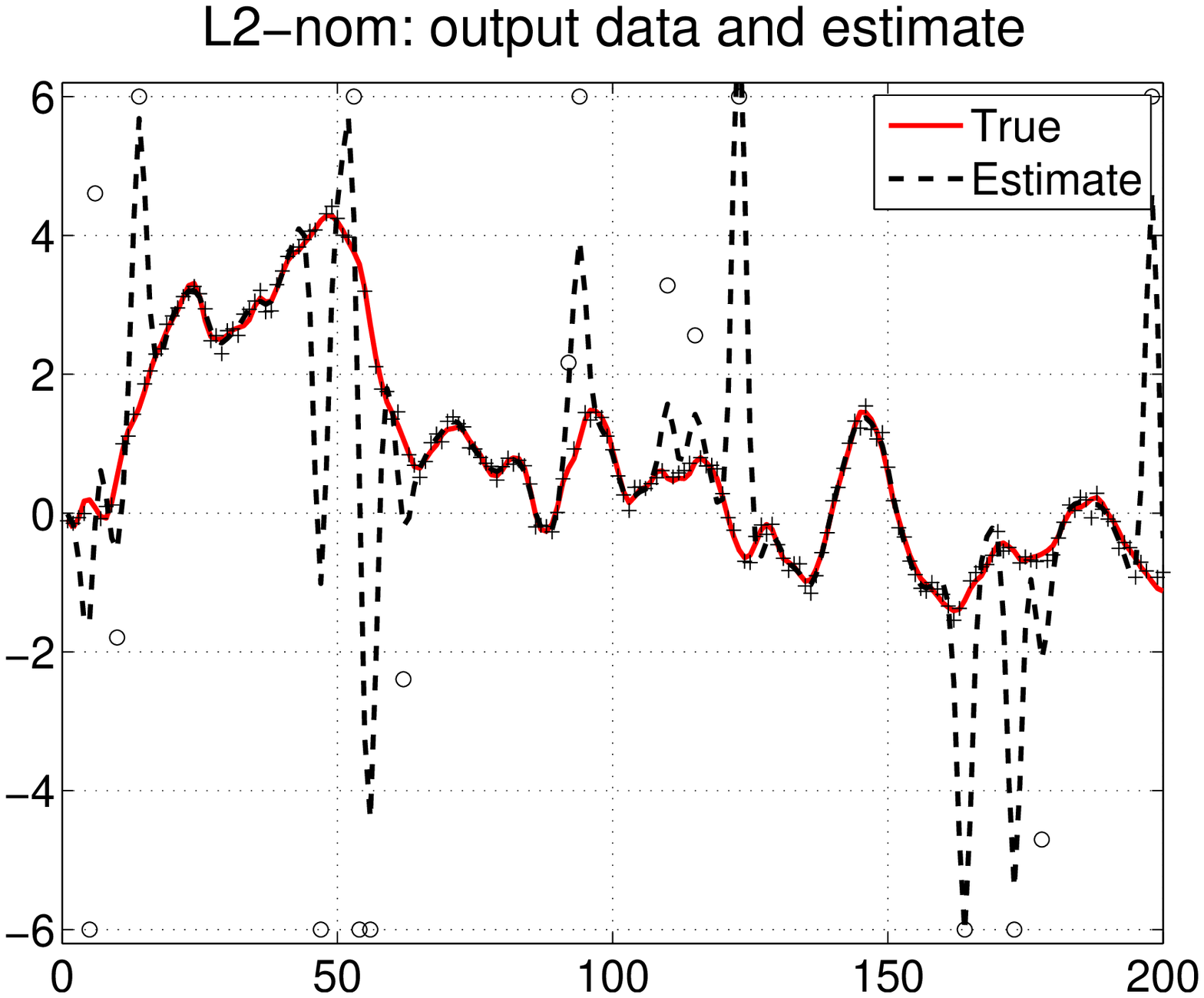}} 
\hspace{.1in}
 { \includegraphics[scale=0.46]{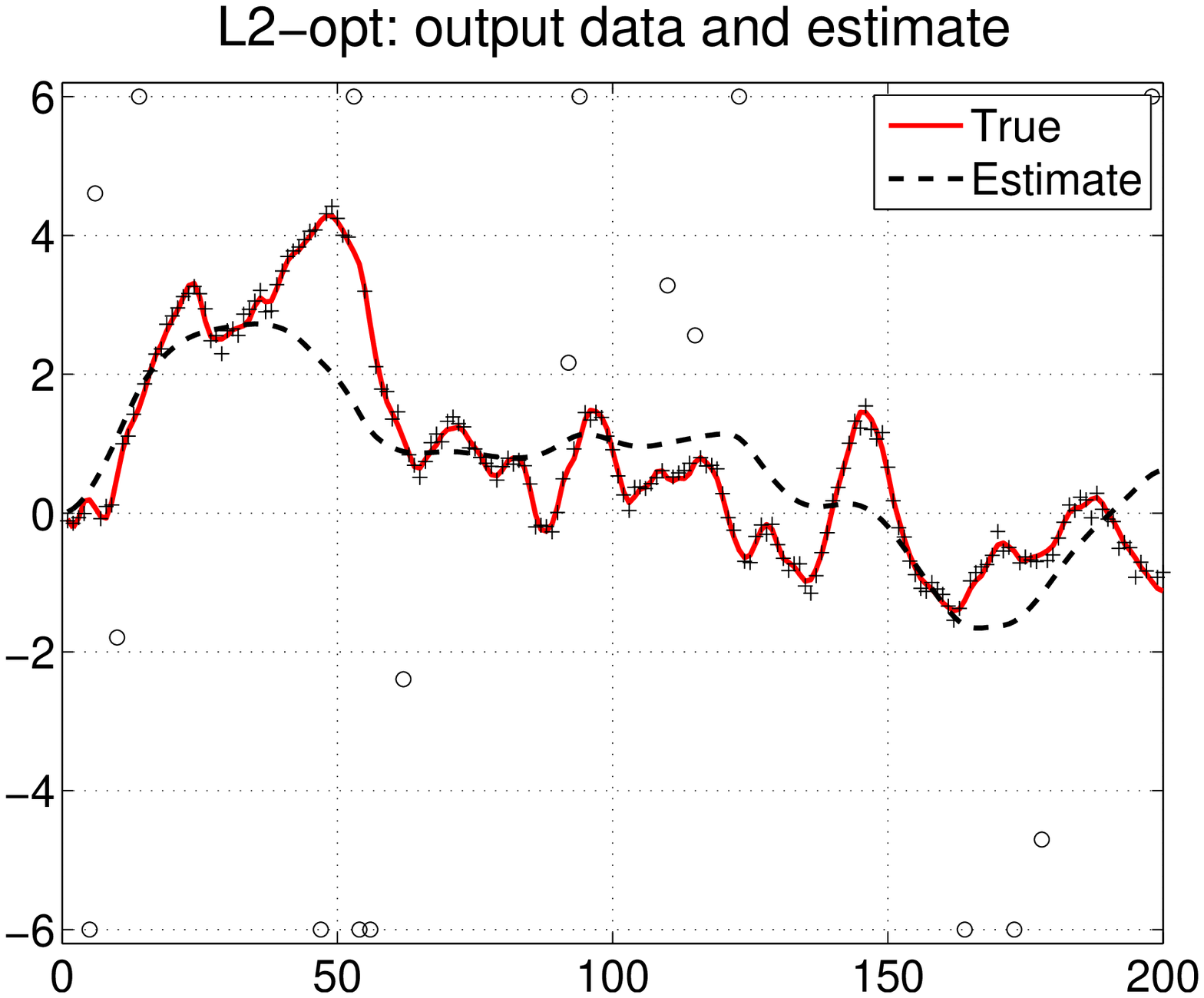}} 
    \end{tabular}
  \caption{{\bf{DC motor with Gaussian disturbances and outliers in output measurements.}} 
   Noiseless output (solid line), measurements  ($+$) and outliers ($\circ$). 
  {\it{Left:}}  Kalman estimates (dashed line) with assumed nominal measurement error variance (0.01).  {\it{Right:}} Kalman estimates 
  (dashed line) from the optimal linear smoother which uses the correct measurement error variance (10.009). } 
    \label{FigOut12}
     \end{center}
\end{figure*}

{\bf{Outliers corrupting output data:}} Consider now a situation where 
the disturbance $d_t$ can be well modeled as a Gaussian process. 
So, there is no impulsive noise entering the system. 
In particular, we set
$d_t \sim \G(0,0.1^2)$, so that $v_t$ is now Gaussian with covariance
$$
Q_t=0.1^2 \left(\begin{array}{c}11.81 \\0.62 \end{array}\right) \left(\begin{array}{cc}11.81 & 0.62 \end{array}\right).
$$
The outputs $y_t$ are instead contaminated by outliers, 
i.e. unexpected measurements noise model deviations.
In particular, output data are corrupted by a mixture of 
two normals with a fraction of outliers contamination equal to $\alpha=0.1$; i.e.,
$$
e_t \sim (1-\alpha) \G(0,\sigma^2) + \alpha \G (0,(100\sigma)^2). 
$$
Thus, outliers occur with probability
$0.1$, and are generated from a distribution with standard deviation 100 times 
greater than that of the nominal. 
We consider the problem of reconstructing the angle of the motor shaft 
(the second state component which corresponds to the noiseless output) setting 
$$
x_0 =  \left(\begin{array}{c}0 \\0\end{array}\right),   \ u_t=0, \ \sigma^2=0.1^2.
$$
An instance of the problem is shown in Fig. \ref{FigOut12}.
The two panels display the noiseless output (solid line), the  accurate
measurements affected by the noise with nominal variance (denoted by $+$) 
and the outliers (denoted by $\circ$ with values outside the range 
$\pm 6$ displayed on the boundaries of the panel). 
The left panel displays the estimate (dashed line) obtained by 
the classical Kalman smoother, called L$_2$-nom, with the variance noise set to $\sigma^2$.

Note that this estimator does not match any of the criteria (\ref{eq:cm}-\ref{eq:cm4}).
In fact, this example represents a situation where the contamination is totally unexpected
and the smoother is expected to work under nominal conditions. 
One can see that the reconstructed profile is very sensitive to outliers. The right panel 
shows the estimate (dashed line) returned by the optimal linear estimator L$_2$-opt
(\ref{eq:cm4}), obtained by setting
the noise variance to $(1-\alpha) \sigma^2 + \alpha (100\sigma)^2$. 

In this case, the smoother is aware of the true variance of the signal; nonetheless, 
the reconstruction is still not satisfactory, since it cannot track the true output profile 
given the high measurement variance; the best linear estimate essentially averages the signal. 
Manipulating noise statistics 
is clearly not enough; to improve the estimator performance, 
we must change our model for the underlying distribution of the errors $e_t$.

\subsection{Scope of the survey}

In light of this discussion and examples, it is natural to turn to the optimization (MAP) interpretation (\ref{eq:opt}) to design formulations 
and estimators that perform well in alternative and more general situations. 
The connection between numerical analysis and optimization and various kinds of smoothers 
has been growing stronger over the years~\cite{Ljung1976,Paige1977,Bell1993,aravkin2013kalman}. 
It is now clear that many popular algorithms in the engineering literature, including Rauch-Tung-Striebel (RTS) smoother
and the Mayne-Fraser (MF) smoother, can be viewed as specific linear algebraic techniques to solve an optimization objective
whose structure is closely tied to dynamic inference. 
Indeed, recently, Kalman smoothing has seen a remarkable renewal in terms of modern techniques
and extended formulations based on emerging practical needs. 
This resurgence has been coupled with the development of
new computational techniques and the intense progress in convex optimization in the last two decades has led to a vast literature 
 on finding good state estimates in these more general cases.  
Many novel contributions to theory and algorithms related to Kalman smoothing, 
and to dynamic system inference in general, have come from statistics, 
engineering, and numerical analysis/optimization communities. 
However, while the statistical and engineering viewpoints are pervasive in the literature, 
the optimization viewpoint and its accompanying modeling and computational 
power is less familiar to the control community. Nonetheless, the optimization perspective
has been the source of a wide range of astonishing recent advances across the board in signal 
processing, control, machine learning, and large-scale data analysis. In this survey,  we will
 show how the optimization viewpoint allows the control and signal processing community
great freedom in the development of novel modeling and inference frameworks
for dynamical systems.

Recent approaches in dynamic systems inference replace  
quadratic terms, as in (\ref{eq:opt}), with suitable convex functions, as in (\ref{eq:optVJ}).
In particular, new smoothing schemes deal with 
sparse dynamic models~\cite{Angelosante2009}, 
methods for  tracking abrupt changes~\cite{Ohlsson2011}, robust formulations~\cite{Farahmand2011,Aravkin2011tac}, 
inequality constraints on the state~\cite{Bell2008}, 
and sum of norms~\cite{Ohlsson2011}, 
many of which can be modeled using the general class called piecewise linear quadratic (PLQ) penalties~\cite{JMLR:v14:aravkin13a,RTRW}.
All of these approaches are based on an underlying body of theory and methodological tools developed 
in statistics, machine learning, kernel methods \cite{Bottou07,Hofmann08,Scholkopf01b,Chan2011},
 and convex optimization \cite{Boyd04}. Advances in sparse tracking~\cite{Boydl1trend,Angelosante2009,Ohlsson2011} are 
based on LASSO (group LASSO) or elastic net techniques~\cite{Lasso96,LARS2004,Yuan2006,EN_2005},
which in turn use coordinate descent, see e.g.~\cite{Bert,FriedCD10,Dinuzzo11}.  
Robust methods~\cite{Aravkin2011tac,Farahmand2011,JMLR:v14:aravkin13a,Chang2013,Aga2011}
rely on Huber \cite{Huber09} or Vapnik losses, leading to support vector regression \cite{Drucker97,Gunter07,Chia2012} for state space models,
and take advantage of interior point optimization methods~\cite{KMNY91,NN94,Wright:1997}.
Domain constraints are important for most applications, 
including camera tracking, fault diagnosis, chemical processes, vision-based systems, 
target tracking, biomedical systems, robotics, and navigation~\cite{haykin2001kalman,simon2010kalman}.
Modeling these constraints allows {\it a priori} information to be encoded into dynamic inference formulations, 
and the resulting optimization problems can also be solved using interior point methods~\cite{bell2009inequality}.

Taking these developments into consideration, the aims of this survey are as follows. 
First, our goal is to firmly establish the connection between classical algorithms, including the RTS and MF smoothers, to the optimization perspective in the least squares case. 
This allows the community to view existing efficient algorithms as modular subroutines that can be exploited in new formulations. 
Second, we will survey modern regression approaches from statistics and machine learning, based on
new convex losses and penalties, highlighting their usefulness in the context of dynamic inference.
These techniques are effective both in designing models for process disturbances $v_t$
as well as robust statistical models for measurement errors $e_t$. 
Our final goal is two-fold: we want to survey algorithms for generalized smoothing formulations, 
but also to understand the theoretical underpinnings for the design and analysis of 
such algorithms. To this end, we include a self-contained tutorial of convex analysis, 
developing concepts of duality and optimality conditions from fundamental principles, 
and focused on the general Kalman smoothing context.  With this foundation, 
we review optimization techniques to solve 
all general formulations of Kalman smoothers, 
including both first-order splitting methods, and second order (interior point) methods. 

In many applications, process and measurement models may be nonlinear. These cases fall outside the scope of the 
current survey, since they require solving a nonconvex problem. In these cases, 
particle filters~\cite{arulampalam2002tutorial} and unscented methods~\cite{wan2000unscented} are very popular. 
An alternative is to exploit the composite structure of these problems, and apply a generalized Gauss-Newton 
method~\cite{burke1995gauss}. For detailed examples, see~\cite{aravkin2011laplace} and~\cite{aravkin2014robust}.

{\bf Roadmap of the paper}:
In section~\ref{L2opt}, we show the explicit connection between RTS and MF smoothers and the least squares formulation. This builds the foundation for efficient general methods that exploit underlying state space structure
of dynamic inference. In section~\ref{sec:stats}, we present a general modeling framework where 
error distributions~\eqref{eq:wgn} can come from a large class of log-concave densities, and discuss 
important applications to impulsive disturbances and robust smoothing. We also show how to incorporate 
state-space constraints.  In section~\ref{sec:numerics}, 
we present empirical results for the examples in the paper, showing the practical effect of the proposed methods. 
All examples are implemented using an open source software package IPsolve\footnote{https://github.com/saravkin/IPsolve}. 
A few concluding remarks end the paper. Two appendices are provided. 
The first discusses smoothing under correlated noise and singular covariance matrices, and the second
a brief tutorial on the tools from convex analysis that are useful
to understand the algorithms presented in section~\ref{sec:algo} and applied in section~\ref{sec:numerics}.

\section{Kalman smoothing, block tridiagonal systems and classical schemes}\label{L2opt}

To build an explicit correspondence between least squares problems and classical smoothing schemes, 
we first introduce data structures that explicitly embed the entire state sequence, measurement sequence, 
covariance matrices, and initial conditions into a simple form. Given a sequence of column vectors $\{ v_k \}$
and matrices $ \{ T_k \}$ let
\[
\mbox{vec} ( \{ v_k \} )
:=
\begin{bmatrix}
v_1 \\ v_2  \\ \vdots \\ v_N
\end{bmatrix}
\; , \;
\mbox{diag} ( \{ T_k \} )
:=
\begin{bmatrix}
T_1    & 0      & \cdots & 0 \\
0      & T_2    & \ddots & \vdots \\
\vdots & \ddots & \ddots & 0 \\
0      & \cdots & 0      & T_N
\end{bmatrix} .
\]
We make the following definitions:
\begin{equation}
  \label{defs}
\begin{aligned}
R       & :=  \mbox{diag} ( \{R_1,  R_2, \dots, R_N\}  )\in \mathbb{R}^{mN \times mN}
\\
Q       & :=  \mbox{diag} ( \{\Pi,Q_0,  Q_1, \dots, Q_{N-1}\}  ) \in\mathbb{R}^{n(N+1)\times n(N+1)}
\\
x       & := \mbox{vec} ( \{x_0,x_1,  x_2, \dots, x_N\} ) \in \mathbb{R}^{n(N+1)\times 1}
\\
y      & := \mbox{vec} (\{y_1,  y_2, \dots, y_N\}) \in \mathbb{R}^{mN\times 1}
\\
z      &  := \mbox{vec} (\{\mu, B_0u_0, \dots, B_{N-1}u_{N-1}\})  \in\mathbb{R}^{n(N+1)\times 1}
\end{aligned}
\end{equation}

and
\begin{equation}
\label{processG}
\begin{aligned}
A  & := \begin{bmatrix}
    \mbox{I}   & 0      &          &
    \\
    -A_0   & \mbox{I}  & \ddots   &
    \\
        & \ddots &  \ddots  & 0
    \\
        &        &   -A_{N-1}  & \mbox{I}
\end{bmatrix} ,
\quad   C & := 
\begin{bmatrix}
0 & C_1    & 0      & \cdots & 0 \\
0 & 0      & C_2    & \ddots & \vdots \\
\vdots & \vdots & \ddots & \ddots & 0 \\
0 & 0      & \cdots & 0      & C_N
\end{bmatrix} ,
\end{aligned}
\end{equation}
where $A  \in\mathbb{R}^{n(N+1) \times n(N+1)}$ and $C  \in\mathbb{R}^{mN \times n(N+1)}$.
Using definitions~\eqref{defs} and~\eqref{processG}, 
problem~\eqref{eq:opt} can be efficiently stated as  
\begin{equation}\label{fullLS}
\min_{x}   \left\|R^{-1/2}(y- Cx) \right\|^2 + \left\|Q^{-1/2}(z - Ax)\right\|^2\;.
\end{equation}
The solution to~\eqref{fullLS} can be obtained by solving the linear system
\begin{equation}\label{smoothingSol}
(C^\top R^{-1} C + A^\top Q^{-1} A) x =  r
\end{equation}
where 
$$
r:= C^\top R^{-1} y + A^\top Q^{-1}z\;.
$$

The linear operator in~\eqref{smoothingSol} is a positive definite symmetric block tridiagonal (SBT) system. 
Direct computation gives   
$$
C^\top R^{-1} C + A^\top Q^{-1} A
=
\begin{bmatrix}
F_0 & G_0^\top & 0 & \cdots&0\\
G_0 & F_1 & G_1^\top & \cdots &\vdots\\
\vdots&G_1 &\ddots & \ddots &  \vdots  \\
0&\cdots& \ddots&& G_{N-1}^T\\
0& \cdots&0& G_{N-1} & F_N
\end{bmatrix} ,
$$
%
a symmetric positive definite block tridiagonal system in $\mathbb{R}^{n(N+1) \times n(N+1)}$,
with $F_t \in \R^{n\times n}$ and $G_t \in \R^{n\times n}$ defined as follows:
$$
\begin{aligned}
F_0&:=\Pi^{-1}+ A_0^\top Q_0^{-1} A_0\\
F_t &:= Q_{t-1}^{-1} + A_{t}^\top Q^{-1}_{t} A_{t} +C_t^\top R_t^{-1}C_t, \quad t=1,\ldots,N\\
G_t &:= -Q_t^{-1}A_t, \quad t=0,\ldots,N-1
\end{aligned}
$$
using the convention $A_{N}^\top Q^{-1}_{N} A_{N}=0$.

We now present two popular smoothing schemes, the RTS and MF.
In our algebraic framework, both of them return the solution of the Kalman smoothing problem (\ref{fullLS})
by efficiently solving the block tridiagonal system (\ref{smoothingSol}), which can be rewritten as
\begin{equation}
\label{BlockTridiagonalEquation}
\left( \begin{matrix}
F_0 & G_0^\top & 0 & \cdots&0\\
G_0 & F_1 & G_1^\top & \cdots &\vdots\\
\vdots&G_1 &\ddots & \ddots &  \vdots  \\
0&\cdots& \ddots&& G_{N-1}^T\\
0& \cdots&0& G_{N-1} & F_N
\end{matrix} \right)
\left( \begin{array}{c} 
	x_0 \\ \rule{0em}{1.5em} x_1 \\ \vdots \\ x_{N-1} \\ x_N
\end{array} \right)
=
\left( \begin{matrix} 
	r_0 \\ \rule{0em}{1.5em} r_1 \\ \vdots \\  r_{N-1} \\ r_N
\end{matrix} \right).
\end{equation}
In particular, the RTS scheme coincides with the 
forward backward algorithm as described in \cite[algorithm 4]{Bell2000}
while the MF scheme can be seen as a block tridiagonal solver
exploiting two filters running in parallel. Full analysis of these algorithms, as well as others, are presented in~\cite{aravkin2013kalman}.

\begin{algorithm}
\caption{Rauch Tung Striebel (Forward Block Tridiagonal scheme)}\label{ForwardAlgorithm}
The inputs to this algorithm are 
\( \{ G_t \}_{t=0}^{N-1} \),
\( \{ F_t \}_{t=0}^N \),
and
\( \{ r_t \}_{t=0}^N \) where, for each $t$,
\( G_t \in \R^{n \times n} \),
\( F_t \in \R^{n \times n} \), and
\( r_t \in \R^{m } \).
The output is the sequence   $\{\hat{x}^N_t\}_{t=0}^N$ that solves equation~(\ref{BlockTridiagonalEquation}),
with each $\hat{x}^N_t \in \R^{n}$. 
\begin{enumerate}

\item
Set \( d_0^f = F_0  \) and \( s_0^f = r_0 \).

\noindent
For  \( t = 1 \)  to  \( N \) :

\begin{itemize}
\item
Set \( d_t^f = F_{t} - G_{t-1} (d_{t-1}^f)^{-1} G_{t-1}^\top  \).
\item
Set \( s_t^f = r_t - G_{t-1} (d_{t-1}^f)^{-1} s_{t-1} \).
\end{itemize}

\item
Set \( \hat{x}^N_t = (d_N^f)^{-1} s_N \).

\noindent
For \( t = N-1 \) to  \( 0 \) :
\begin{itemize}
\item
Set \( \hat{x}^N_t = (d_t^f)^{-1} ( s_t^f - G_{t}^\top \hat{x}^N_{t+1} ) \).
\end{itemize}

\end{enumerate}
\end{algorithm}

\begin{algorithm}
\caption{Mayne Fraser (Two Filter Block Tridiagonal scheme)}\label{algoDiag}
The inputs to this algorithm are 
\( \{ G_t \}_{t=0}^{N-1} \),
\( \{ F_t \}_{t=0}^N \),
and
\( \{ r_t \}_{t=0}^N \) where, for each $t$,
\( G_t \in \R^{n \times n} \),
\( F_t \in \R^{n \times n} \), and
\( r_t \in \R^{m } \).
The output is the sequence   $\{\hat{x}^N_t\}_{t=0}^N$ that solves equation~(\ref{BlockTridiagonalEquation}),
with each $\hat{x}^N_t \in \R^{n}$. 
\begin{enumerate}

\item
Set \( d_0^f = F_0  \) and \( s_0^f = r_0 \).

\noindent
For  \( t = 1 \)  to  \( N \) :

\begin{itemize}
\item
Set \( d_t^f = F_{t} - G_{t-1} (d_{t-1}^f)^{-1} G_{t-1}^\top  \).
\item
Set \( s_t^f = r_t - G_{t-1} (d_{t-1}^f)^{-1} s_{t-1} \).
\end{itemize}

\item
Set \( d_N^b = F_N  \) and \( s_N^b = r_N \).

\noindent
For \( t = N-1, \ldots , 0 \),
\begin{itemize}
\item
Set 
\( d_t^b = F_t - G_{t}^\top (d_{t+1}^b)^{-1} G_{t} \).
\item
Set
\( s_t^b = r_t - G_{t}^\top (d_{t+1}^b)^{-1} s_{t+1} \).
\end{itemize}

\item
For \( t = 1 , \ldots , N \)
\begin{itemize}
\item
Set 
\( \hat{x}^N_t  = (d_k^f + d_k^b - b_k)^{-1} ( s_k^f + s_k^b - r_k) \).
\end{itemize}
\end{enumerate}
\end{algorithm}


\section{General formulations: convex losses and penalties, and  statistical properties of the resulting estimators}
\label{sec:stats}

In the previous section, we showed that Gaussian assumptions on process disturbances $v_t$ and measurement errors $e_t$
lead to least squares formulations  \eqref{eq:opt} or \eqref{fullLS}. One can then view classic smoothing algorithms as numerical subroutines for solving these least squares problems. 
In this section, we generalize the Kalman smoothing model to allow log-concave  distributions for $v_t$ and $e_t$ in model~\eqref{eq:wgn}. 
This allows more general {\it convex} disturbance and error measurement models, 
and the log-likelihood (MAP) problem~\eqref{fullLS} becomes a more general {\it convex} inference problem.

In particular, we consider the following general convex formulation:
\begin{equation}\label{fullGen}
\min_{x \in \mathcal X }   V\left(R^{-1/2}(y- Cx)\right)  +  \gamma J\left(Q^{-1/2}(z - Ax)\right). 
\end{equation}
where $x\in \mathcal{X}$ specifies a feasible domain for the state,  
$V: \mathbb{R}^{mN} \rightarrow\mathbb{R}$ measures the discrepancy between observed and 
predicted data (due to noise and outliers),
while $J: \mathbb{R}^{n(N+1)} \rightarrow \mathbb{R}$
measures the discrepancies between predicted and observed state transitions, 
due to the net effect of factors outside the process model;  we can think of these discrepancies 
as `process noise'. The structure of this problem is related to Tikhonov regularization and 
inverse problems~\cite{Tikhonov,Bertero1,DeVito2005}.  
In this context, $\gamma$ is called the {\it regularization parameter}
and has a link to the (typically unknown) scaling of the pdfs of
$e_t$ and $v_t$ in (\ref{Generalpdf}).
The choice of $\gamma$ controls
the tradeoff between bias and variance, and it has to be tuned from data. 
Popular tuning methods include cross-validation or generalized cross-validation~\cite{Rice1986,Golub79,Hastie01}.

Problem \eqref{fullGen} is overly general. In practice we restrict $V$ and $J$ to be functions following the block structure of their arguments,
i.e. sums of terms $V_t \left(R_t^{-1/2}(y_t-C_t x_t)\right)$ and $J_t\left(Q_t^{-1/2}(x_{t+1} - A_t x_t -B_t u_t)\right)$, leading to the objective
already reported in (\ref{eq:optVJ}). The terms $V_t: \R^m \rightarrow \R$ and $J_t: \R^n \rightarrow \R$
can then be linked to the MAP interpretation of the state estimate 
 (\ref{Generalpdf})-(\ref{eq:optVJ}), so that $V_t$ is a version of $- \log \Bp_{e_t}$ and $J_t$ is a version of $-\log \Bp_{v_t}$.
 Possible choices for such terms are depicted in Fig. \ref{fig:quadratic}-\ref{fig:enet} and Fig. \ref{GLT-KF}. \\

Domain constraints $x \in \mathcal{X}$ provide a disciplined framework for 
incorporating prior information into the inference problem, which 
improves performance for a wide range of applications. 
Analogously, general $J$ and $V$ allow the modeler to incorporate
information about {\it uncertainty}, both in the process and measurements. 
This freedom in designing~\eqref{fullGen} has numerous benefits. The modeler can choose $J$ to reflect prior knowledge 
on the structure  of the process noise; important examples include sparsity (see Fig.~\ref{Fig12Imp}) and smoothness. 
In addition, she can robustify the formulation in the presence of outliers or non-gaussian errors (see Fig.~\ref{FigOut12}), by 
selecting penalties $V$ that perform well in spite of data contamination. To illustrate, 
we present specific choices for the functions $V$ and $J$ and explain how they can be used 
in a range of modeling scenarios; we also highlight the potential for constrained formulations. 

\subsection{General functions $J$ for modeling process noise}

As mentioned in the introduction, a widely used assumption for the process noise is that it is Gaussian. 
This yields the quadratic loss $\| Q^{-1/2}(z - Ax)\|^2.$ 
However, in many applications prior knowledge on the process disturbance dictates alternative loss functions. 
A simple example is the DC motor in Section 1.1.  We assumed that the process disturbance $v_t$ is impulsive.  One therefore expects that the disturbance $v_t$ should be zero most of the time, 
while taking non-zero values at a few unknown time points.  
If each $v_t$ is scalar, a natural way to regulate the number of non-zero components in $\textrm{vec}\left(\{v_t\}\right)$ 
is to use the $\ell_0$ norm for $J$ in~\eqref{fullGen}: 
\[
J(z - Ax; Q)= \| Q^{-1/2}(z - Ax) \|_0,
\]
where $\|z\|_0$ counts the number of nonzero elements of $z$. 

\begin{figure}
        \centering
        \resizebox {\columnwidth} {!} {
\begin{tikzpicture}
    \begin{axis}[grid = major,
        xmin=-1.7,xmax=1.7,ymin=-1.7,ymax=1.7,
    ]
    \addplot[domain = -1:1.3, color = blue, very thick]{1-0.5*x};
    \addplot+[fill, mark = none] coordinates
        {(0,1) (1,0) (0,-1) (-1,0)} --cycle;
    \end{axis}
\end{tikzpicture}
\begin{tikzpicture}
   \begin{axis}[grid = major,
        xmin=-1.7,xmax=1.7,ymin=-1.7,ymax=1.7,
    ]
    \addplot[domain = -1:1.3, color = blue, very thick]{1.13-0.5*x};
        \addplot +[fill,domain=0:2*pi,samples=50,mark=none]({cos(deg(x))},{sin(deg(x))});
    \end{axis}
\end{tikzpicture}}
 \caption{\label{fig:balls} When minimizing $\|Ax - b\|$ subject to a $\ell_1$-norm constraint (left panel), the  solution tends to land on a corner, 
 where many coordinates are $0$; in 2D the cartoon, the $x$-coordinate is zero. 
 An $\ell_2$-norm constraint (right panel) does not have this effect.  }
\end{figure}
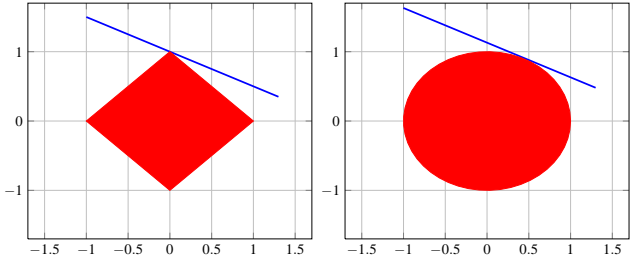

\begin{figure*}[t!]
    \begin{subfigure}[t]{0.33\textwidth}
       \centering
\begin{tikzpicture}
  \begin{axis}[
    thick,
    height=2cm,
    xmin=-2,xmax=2,ymin=0,ymax=1,
    no markers,
    samples=50,
    axis lines*=left, 
    axis lines*=middle, 
    scale only axis,
    xtick={-1,1},
    xticklabels={},
    ytick={0},
    ] 
\addplot[blue, domain=-2:+2]{.5*x^2};
  \end{axis}
  \end{tikzpicture}
  \caption{\label{fig:quadratic}quadratic}
    \end{subfigure}%
    \begin{subfigure}[t]{0.33\textwidth}
        \centering
\begin{tikzpicture}
  \begin{axis}[
    thick,
    height=2cm,
    xmin=-2,xmax=2,ymin=0,ymax=1,
    no markers,
    samples=100,
    axis lines*=left, 
    axis lines*=middle, 
    scale only axis,
    xtick={-1,1},
    xticklabels={},
    ytick={0},
    ] 
  \addplot[red, dashed, domain=-2:+2]{abs(x)};
  \end{axis}
\end{tikzpicture}
    \caption{\label{fig:1norm}$\ell_1$ norm}
    \end{subfigure}
   \begin{subfigure}[t]{0.33\textwidth}
   \centering
   \begin{tikzpicture}
  \begin{axis}[
    thick,
    height=2cm,
    xmin=-2,xmax=2,ymin=0,ymax=1,
    no markers,
    samples=50,
    axis lines*=left, 
    axis lines*=middle, 
    scale only axis,
    xtick={-1,1},
    xticklabels={},
    ytick={0},
    ] 
\addplot[red,domain=-2:-1,densely dashed]{-x-.5};
\addplot[blue, domain=-1:+1]{.5*x^2};
\addplot[red,domain=+1:+2,densely dashed]{x-.5};
\addplot[blue,mark=*,only marks] coordinates {(-1,.5) (1,.5)};
  \end{axis}
\end{tikzpicture}
\caption{\label{fig:Huber}Huber, $\kappa = 1$}
\end{subfigure}    
   \begin{subfigure}[t]{0.32\textwidth}
   \centering
\begin{tikzpicture}
  \begin{axis}[
    thick,
    height=2cm,
    xmin=-2,xmax=2,ymin=0,ymax=1,
    no markers,
    samples=50,
    axis lines*=left, 
    axis lines*=middle, 
    scale only axis,
    xtick={-0.5,0.5},
    xticklabels={},
    ytick={0},
    ] 
    \addplot[red,domain=-2:-0.5,densely dashed] {-x-0.5};
    \addplot[domain=-0.5:+0.5] {0};
    \addplot[red,domain=+0.5:+2,densely dashed] {x-0.5};
    \addplot[blue,mark=*,only marks] coordinates {(-0.5,0) (0.5,0)};
  \end{axis}
\end{tikzpicture}
\caption{\label{fig:Vapnik} Vapnik, $\epsilon = 0.5$}
\end{subfigure}  
   \begin{subfigure}[t]{0.32\textwidth}
   \centering
\begin{tikzpicture}
  \begin{axis}[
    thick,
    height=2cm,
    xmin=-2,xmax=2,ymin=0,ymax=1,
    no markers,
    samples=50,
    axis lines*=left, 
    axis lines*=middle, 
    scale only axis,
    xtick={-1,1, -.5, .5},
    xticklabels={},
    ytick={0},
    ] 
\addplot[domain=-0.25:+0.25] {0};
\addplot[red,domain=-2:-1,densely dashed]{-x-.5-.5*.25};
\addplot[blue, domain=-1:-.25]{.5*x^2-.5*.25};
\addplot[blue, domain=.25:1]{.5*x^2-.5*.25};
\addplot[red,domain=+1:+2,densely dashed]{x-.5-.5*.25};
\addplot[blue,mark=*,only marks] coordinates {(-1,.5-.5*.25) (1,.5-.5*.25)(-.45, 0) (.45, 0)};
  \end{axis}
\end{tikzpicture}
\caption{\label{fig:sel} Huber ins. loss}
\end{subfigure}  
   \begin{subfigure}[t]{0.328\textwidth}
   \centering
\begin{tikzpicture}
  \begin{axis}[
    thick,
    height=2cm,
    xmin=-2,xmax=2,ymin=0,ymax=1,
    no markers,
    samples=100,
    axis lines*=left, 
    axis lines*=middle, 
    scale only axis,
    xtick={-1,1},
    xticklabels={},
    ytick={0},
    ] 
\addplot[amethyst, domain=-2:+2]{.5*x^2 + 0.5*abs(x)};
  \end{axis}
\end{tikzpicture}
\caption{\label{fig:enet} elastic net}
\end{subfigure}  
    \caption{Important penalties for errors and process models. }
\end{figure*}
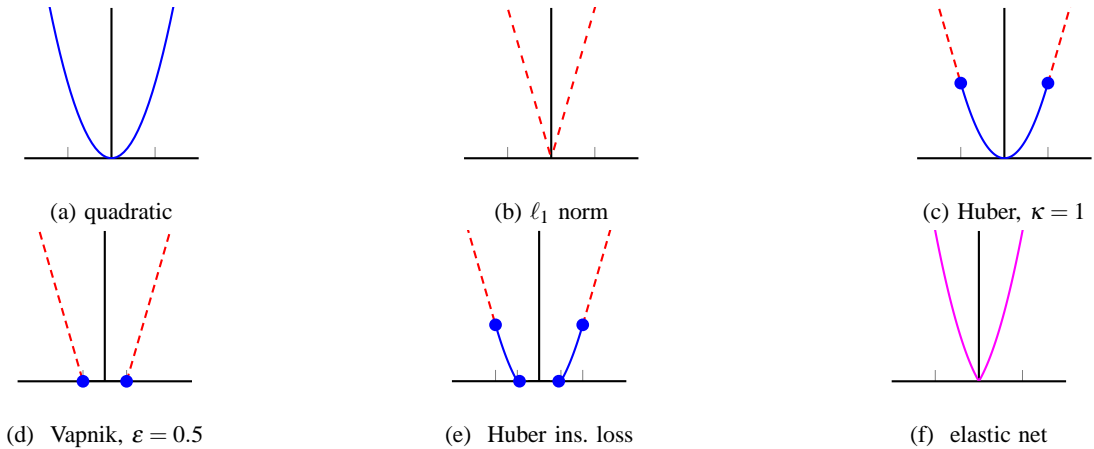

\textit{Sparsity promotion via $\ell_1$ norm.}
The $\ell_0$ norm, however, is non-convex, and solving optimization problems involving the $\ell_0$ norm is NP-hard (combinatorial). 
Tractable approaches can be designed by replacing the $\ell_0$ norm with a convex relaxation, 
the $\ell_1$ norm, $\|x\|_1 = \sum |x_i|$. 
The $\ell_1$ norm is nonsmooth and encourages sparsity, see Fig.~\ref{fig:1norm}.
The use of the $\ell_1$ norm in lieu of the $\ell_0$ norm is now common practice,
especially in compressed sensing~\cite{candes2006nos,Donoho2006} and statistical learning, see e.g.~\cite{Hastie01}. 
The reader can gain some intuition by considering the intersection of a general 
hyperplane with the $\ell_1$ ball and $\ell_2$ ball in Fig.~\ref{fig:balls}. 
The intersection is likely to land on a corner, 
which means that adding a $\ell_1$ norm constraint (or penalty) tends to select solutions with many zero elements.\\ 
For the case of scalar-valued process disturbance $v_t$, we can set $J$ to be the $\ell_1$ norm and obtain the  problem
\begin{equation}\label{fullGenLasso}
\min_{x }   \frac{1}{2}\|R^{-1/2}(y- Cx)\|^2  +  \gamma \|  Q^{-1/2}(z - Ax) \|_1,
\end{equation}
where $\gamma$ is a penalty parameter controlling the tradeoff between measurement fit and number of non-zero components in process disturbance --- larger $\gamma$ implies a larger number of zero process disturbance elements, at the cost of increasing the bias of the estimator.\\ 
Note that the vector norms in (\ref{fullGenLasso}) translate to term-wise norms of the time components as in (\ref{eq:optVJ}). Problem~\eqref{fullGenLasso} is analogous to the LASSO problem~\cite{Tibshirani96}, 
originally proposed in the context of linear regression. 
Indeed, the LASSO problem minimizes the sum of squared residuals regularized by the $\ell_1$ penalty on the regression coefficients. In the context of regression, the LASSO has been shown to have strong statistical guarantees, including prediction error consistency~\cite{GeerBuhl09}, consistency of the parameter estimates in $\ell_2$ or some other norm~\cite{GeerBuhl09,MeinshausenYu09}, as well as variable selection consistency~\cite{Meinshausen06,Wainwright06,Zhao06}. However, this connection is limited in the dynamic context: if we think of Kalman smoothing as linear regression, note from~\eqref{fullGenLasso} that the measurement vector $y$ is a {\it single} observation of the parameter (state sequence) $x$, so asymptotic consistency results are not relevant. 
More important is the general idea of using the $\ell_1$ norm to promote sparsity of the 
right object, in this case, the residual $ Q^{-1/2}(z - Ax)$, which corresponds 
to our model of impulsive disturbances.

\textit{Elastic net penalty.}
Suppose we need a penalty that is nonsmooth at the origin, but has quadratic growth 
in the tails. For example, taking $J$ with these properties is useful in the context of our model for 
impulsive disturbances, if we believed them to be sparse, and also considered large disturbances unlikely. 
The {\it elastic net} shown in Fig.~\ref{fig:enet} has these properties --- it is a weighted sum 
$\alpha \|\cdot\|_1 + (1-\alpha)\|\cdot\|_2^2$. The elastic net penalty has been widely used for sparse regularization with correlated predictors~\cite{ZouHuiHastie:2005,EN_2005,li2010Bayesian,de2009elastic}. 
Using an elastic net constraint has a {\it grouping effect}~\cite{ZouHuiHastie:2005}.
Specifically, when minimizing $\frac{1}{2}\|Ax-b\|^2$ with an elastic net constraint,  
the distance between estimates $\hat x_i$ and $\hat x_j$ is proportional to $\sqrt{1-\kappa_{ij}}$, 
where $\kappa_{ij}$ is the correlation between the corresponding columns of $A$. 
In our context, in case of nearly perfectly correlated impulsive disturbances (either all present or all absent), the
elastic net can discover the entire group, while the $\ell_1$ norm alone usually picks a single member of the group. 

\textit{Group sparsity.}
If the process disturbance is known to be grouped (e.g. a disturbance vector is always present or absent 
for each time point), $J(\cdot)$ can be set to the mixed $\ell_{2,1}$ norm, where the $\ell_2$ norm is applied to each block of $Q_t^{-1/2}(z_t - A_tx_t)$, yielding the following Kalman smoothing formulation:
\begin{equation}\label{fullGenGrLasso}
\min_{x }  \|R^{-1/2}(y- Cx)\|^2  +  \gamma \sum_{t=1}^N  \left\|Q_t^{-1/2}(z_t - A_tx_t)\right\|_2,
\end{equation}
where $\gamma$ is again a penalty parameter controlling the tradeoff between measurement fit and number of non-zero components in process disturbance.  
Note that the objective is still of the type (\ref{eq:optVJ}) with a  
penalty term that now corresponds to the
sparsity inducing $\ell_1$ norm applied to groups of process disturbances $v_t$, where the 
$\ell_2$ norm used as the intra-group penalty. 
This group penalty has been widely used in statistical learning where it is referred to as the ``group-LASSO'' penalty. 
Its purpose is to select important factors, each represented by a group of derived variables, for joint model selection and estimation in regression. 
In the state estimation context, the estimator (\ref{fullGenGrLasso}) was proposed in \cite{Ohlsson2011} and will be used later on in Section \ref{sec:examples}
to solve the  impulsive inputs problem described in section \ref{MotEx}.
The group $\ell_{2,1}$ penalty was originally proposed in the context of linear regression in~\cite{YuaLi06}. 
General $\ell_{q,1}$ regularized least squares formulations (with $q \ge 2$) were subsequently studied in~\cite{Tro06,ZhaRoc06,YuaLi06,Jac09} and shown to have strong statistical guarantees, including convergence rates in $\ell_2$-norm~\cite{Lou09,Bar08}) as well as model
selection consistency~\cite{OboWaiJor08,NegWai08}.

\subsection{General functions $V$ to model measurement errors}

Gaussian assumption on model deviation is not valid in many cases.
Indeed, heavy tailed errors are frequently observed in applications such 
glint noise~\cite{HewerMartin1987}, air turbulence~\cite{FSI:2013}, 
and asset returns~\cite{Rachev:2003} among others.
The resulting state estimation problems can be addressed by adopting the penalties $J$ introduced above.
But, in addition, corrupted measurements might occur due to equipment malfunction, secondary sources of noise or other anomalies. The quadratic loss is not robust with respect to the presence of outliers in the data~\cite{Hub,Aravkin2011tac,Farahmand2011,Gao2008}, as seen in Fig.~\ref{FigOut12},
leading to undesirable behavior of resulting estimators.
This calls for the design of new losses $V$.

One way to derive a robust approach is to assume that the noise comes from a probability density 
with tail probabilities larger (heavier) than those of the Gaussian, and consider the maximum a posteriori (MAP)
problem derived from the corresponding negative log likelihood function. 
For instance the Laplace distribution $c\exp(-\|x\|_1)$ corresponds to the $\ell_1$ loss function by this approach, see Fig~\ref{fig:1norm}.  The {\it tail probabilities}  $P(|x| > t)$ of the standard Laplace distribution are greater than that of the Gaussian; so larger observations are more likely under this error model.   
Note, however, that the $\ell_1$ loss also has a nonsmooth feature at the origin (which is exactly why we considered it as a choice for $J$ in the previous section). 
In the current context, when applied to the measurement residual $Hx- z$, 
the approach will sparsify the residual, i.e. fit a portion of the data exactly. 
Exact fitting of some of the data may be reasonable in some contexts, 
but undesirable in many others, where we mainly care about guarding against outliers, 
and so only the {\it tail behavior} is of interest. 
In such settings, the Huber Loss~\cite{huber81} (see Fig~\ref{fig:Huber}) is a more suitable model,
as it combines the $\ell_2$ loss for small errors with the absolute loss for larger errors.
Huber~\cite{huber81} showed that this loss is optimal over a particular class of errors
\[
(1-\epsilon)\mathcal{N} + \epsilon \mathcal{M},
\]
where $\mathcal{N}$ is Gaussian, and $\mathcal{M}$ is unknown; the level $\epsilon$ is then related to the Huber parameter $\kappa$.

Another important loss function is the Vapnik $\epsilon$-insensitive loss~\cite{Drucker97},
sometimes known as the `deadzone' penalty, see Fig.~\ref{fig:Vapnik},
defined as
\[
V_\epsilon(r) : = \max \{ 0, |r|-\epsilon\},
\]
where $r$ is the (scalar) residual. The $\epsilon$-insensitive loss was originally considered in support vector regression~\cite{Drucker97}, where the `deadzone' helps identify active support vectors, i.e. data elements that determine the solution. 
This penalty has a Bayesian interpretation, as a mixture of Gaussians that may have nonzero means~\cite{pontil2000noise}. In particular, its use yields smoothers that are robust to minor fluctuations 
below a noise floor (as well as to large outliers). Note that the radius of the deadzone $\epsilon$ defines
a noise floor beyond which one cannot resolve the signal. This penalty can also be `huberized', yielding a penalty called `smooth insensitive loss'~\cite{chu2001unified,lee2005epsi,dekel2005smooth}, see Fig.~\ref{fig:sel}.

The process of choosing penalties based on behavior in the tail, near the origin, or at other specific regions
of their subdomains makes it possible to customize the formulation of~\eqref{fullGen} to address a range of situations. 
We can then associate statistical densities to all the penalties in Figs.~\ref{fig:quadratic}-\ref{fig:enet}, and use
this perspective to incorporate prior knowledge about mean and variance of the residuals and process disturbances~\cite[Section 3]{JMLR:v14:aravkin13a}. This allows one to incorporate variance information 
on process components; as e.g. available in the example of Fig.~\ref{FigOut12}.
 
 \textit{Asymmetric extensions.} All of the PLQ losses in Figs.~\ref{fig:quadratic}-\ref{fig:enet} have asymmetric analogues.   
 For example, the asymmetric 1-norm~\cite{koenker1978regression} and asymmetric Huber~\cite{aravkin2014orthogonal} have been used for analysis of heterogeneous datasets, especially in high dimensional inference.

\begin{figure} 
\centering
\begin{tikzpicture}
  \begin{axis}[
    thick,
    width=.3\textwidth, height=2cm,
    xmin=-4,xmax=4,ymin=0,ymax=1,
    no markers,
    samples=50,
    axis lines*=left, 
    axis lines*=middle, 
    scale only axis,
    xtick={-1,1},
    xticklabels={},
    ytick={0},
    ] 
\addplot[domain=-4:+4,densely dashed]{exp(-.5*x^2)/sqrt(2*pi)};
 \addplot[red, domain=-4:+4]{0.5*exp(-.5*abs(x))};
  \addplot[blue, domain=-4:+4]{0.3*exp(-.5*ln(1 + x^2))};
  \end{axis}
\end{tikzpicture}
\begin{tikzpicture}
  \begin{axis}[
    thick,
    width=.3\textwidth, height=2cm,
    xmin=-3,xmax=3,ymin=0,ymax=2,
    no markers,
    samples=50,
    axis lines*=left, 
    axis lines*=middle, 
    scale only axis,
    xtick={-1,1},
    xticklabels={},
    ytick={0},
    ] 
\addplot[domain=-3:+3,densely dashed]{.5*x^2};
 \addplot[red, domain=-3:+3]{.5*abs(x)};
  \addplot[blue, domain=-3:+3]{.5*ln(1 + x^2)};
  \end{axis}
\end{tikzpicture}
\begin{tikzpicture}
  \begin{axis}[
    thick,
    width=.3\textwidth, height=2cm,
    xmin=-3,xmax=3,ymin=-2,ymax=2,
    no markers,
    samples=50,
    axis lines*=left, 
    axis lines*=middle, 
    scale only axis,
    xtick={-1,1},
    xticklabels={},
    ytick={0},
    ] 
\addplot[domain=-3:3,densely dashed]{x};
 \addplot[red, domain=-3:0]{-.5};
  \addplot[red, domain=0:3]{.5};
\addplot[blue, domain=-3:3,]{x/(1 + x^2)};
  \end{axis}
\end{tikzpicture}
    \caption{\label{GLT-KF}
Gaussian (black dashed), Laplace (red solid), and Student's t (blue solid) Densities, Corresponding Negative Log Likelihoods, and Influence Functions.}
\end{figure}
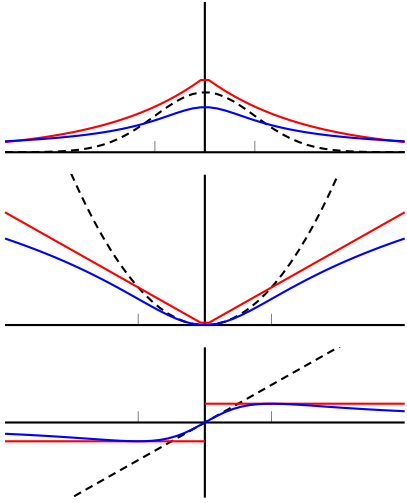

\textit{Beyond convex approaches.} All of the penalty options for $J$ and $V$ presented so far are 
{\it convex}. Convex losses make it possible to provide strong guarantees --- for example, if both 
$J$ and $V$ are convex in~\eqref{fullGen}, then any stationary point is a global minimum. In addition, 
if $J$ has compact level sets (i.e. there are no directions where it stays bounded), then at least 
one global minimizer exists. From a modeling perspective, however, it may be beneficial to choose a non-convex penalty in order to strengthen a particular feature.
In the context of residuals, the need for non-convex loss is motivated by considering the {\it influence function}. This function measures the derivative of the loss with respect to the residual, quantifying the effect of the size of a residual on the loss. 
For nonconstant convex losses, linear growth is the limiting case, and this gives each residual constant influence. 
Ideally the influence function should redescend towards zero for large residuals, so that these are basically ignored.  But redescending influence corresponds to sublinear growth, which excludes convex loss functions. We refer the reader to~\cite{Hampel} for a review of influence-function approaches to robust statistics, including redescending influence functions. An illustration is presented in Figure~\ref{GLT-KF}, contrasting the density, negative log-likelihood, and influence 
function of the heavy-tailed student's t penalty with those of gaussian (least squares) and laplace ($\ell_1$) densities and penalties. 
More formally, consider any scalar density $p$ arising from a symmetric convex
coercive and differentiable penalty $\rho$ via $p(x) = \exp(-\rho(x))$, 
and take any point $x_0$ with $\rho'(x_0) = \alpha_0 > 0$. 

Then, for 
all $x_2 > x_1\geq x_0$ it is shown in~\cite{AravkinFHV:2012} that  the conditional tail distribution induced by $\Bp(x)$
satisfies 
\begin{equation}
\label{memoryFreeIneq}
\Pr(|y| > x_2 \mid |y| > x_1) \leq \exp(-\alpha_0[x_2 - x_1])\;.
\end{equation}
When $x_1$ is large, the condition $|y| > x_1$ indicates that 
we are looking at an outlier. However, as shown by (\ref{memoryFreeIneq}), {\it any}
log-concave statistical model treats the outlier conservatively, 
dismissing the chance that $|y|$ could be significantly bigger than $x_1$. 
Contrast this behavior with that of the Student's t-distribution.
With one degree of freedom, the Student's t-distribution is simply the Cauchy
distribution, with a density proportional to $1/(1 + y^2)$.  Then we
have that
\[
 \lim_{x\to\infty} \Pr(|y|>2x \mid |y|>x) =
 \lim_{x\to\infty} \frac{\frac{\pi}{2}-\arctan(2x)}{\frac{\pi}{2} - \arctan(x)}
 = \frac{1}{2}.
\] 
See~\cite{aravkin2014robust} for a more detailed discussion of non-convex robust approaches to Kalman smoothing using the Student's t distribution. 

Non-convex functions $J$ have also been frequently applied to modeling process noise. 
In particular, see \cite{Wipf_IEEE_TSP_2007,Wipf_ARD_NIPS_2007,Wipf_IEEE_TIT_2011} for a link
between penalized regression problems like LASSO and Bayesian methods.
One classical approach is ARD~\cite{McKayARD}, which exploits hierarchical hyperpriors with
`hyperparameters'  estimated via maximizing the marginal likelihood, 
following the Empirical Bayes paradigm \cite{Maritz:1989}. 
In addition, see~\cite{loh2013regularized,ABCP14} for statistical results in the nonconvex case. 
Although the nonconvex setting is essential in this context, 
it is important to point out that solution methodologies in the above examples 
are based on iterative convex approximations, which is our main focus. 

\subsection{Incorporating Constraints}

Constraints can be important for improving estimation. In state estimation problems, 
constraints arise naturally in a variety of ways. When estimating biological quantities
such as concentration, or physical quantities such height above ground level, we know 
these to be {\it non-negative}. Prior information can induce other constraints; for example, 
if maximum velocity or acceleration is known, this gives {\it bound constraints}. 
Some problems also offer up other interesting constraints: in the absence of maintenance, 
physical systems degrade (rather than improve), giving {\it monotonicity constraints}~\cite{simon2002kalman}.  Both unimodality and monotonicity can be formulated using linear inequality constraints~\cite{aravkin2013new}. 

All of these examples motivate the constraint $x\in\mathcal{X}$ in~\eqref{fullGen}. Since we focus
only on the convex case, we require that $\mathcal{X}$ should be convex. In this paper, we focus on two types of convex sets: 
\begin{enumerate}
\item $\mathcal{X}$ is polyhedral, i.e. given by $\mathcal{X} = \{x: D^Tx \leq d\}$. 
\item $\mathcal{X}$ has a simple projection operator $\proj_{\mathcal{X}}$, where 
\[
\proj_{\mathcal{X}}(y) := \arg\min_{x\in \mathcal{X}} \frac{1}{2}\|x-y\|_2^2.
\]
\end{enumerate}

The cases are not mutually exclusive, for example box constraints are polyhedral and easy to project onto. 
The set $\mathbb{B}_2 := \{x: \|x\|_2 \leq 1\}$ is not polyhedral, but has an easy projection operator:
\[
\proj_{\mathbb{B}_2} (y) = 
\begin{cases}
y/\|y\|_2 &\mbox{if} \quad \|y\|_2 > 1 \\
y & \mbox{else}.
\end{cases}
\]
In general, we let $\mathbb{B}$ denote a closed unit ball for a given norm, and for the $\ell_p$ norms,
this unit ball is denoted by $\mathbb{B}_p$.
These approaches extend to the nonconvex setting. A class of nonconvex Kalman smoothing problems, 
where $\mathcal{X}$ is given by functional inequalities, is studied in~\cite{Bell2008}. 
We restrict ourselves to the convex case, however. 

\section{Efficient algorithms for Kalman smoothing}
\label{sec:algo}

In this section, we present an overview of smooth and nonsmooth methods for convex problems, 
and tailor them specifically to the Kalman smoothing case.  
The section is organized as follows. We begin with a few basic facts about convex sets 
and functions, and review gradient descent and Newton methods  
for smooth convex problems. 
Next, extensions to nonsmooth convex functions are discussed beginning 
with a brief exposition of 
sub-gradient descent and its associated (slow) convergence rate.
We conclude by showing how first- and second-order methods 
can be extended to develop efficient algorithms for the nonsmooth case using the proximity operator,  
splitting techniques, and interior point methods.

\subsection{Convex sets and functions}
\label{sec:convexBasics}

A subset $\calC$ of $\R^n$ is said to be convex if it contains every line segment whose
endpoints are in $\calC$, i.e.,
$$
(1-\lam)x+\lam y\in \calC\quad \forall\, \lam\in[0,1]\quad \mbox{ whenever }\ x,y\in \calC.
$$
For example, the unit ball $\mathbb{B}$ for any norm is a convex set.  

A function $\map{f}{\R^n}{\R\cup\{\infty\}}$ is said to be convex if the secant line between
any two points on the graph of $f$ always lies above the graph of the function, i.e. $\forall \lambda \in [0,1]$: 
$$
f((1-\lam)x+\lam y)\le (1-\lam)f(x)+\lam f(y), \ \  \forall\ x,y\in\R^n.
$$

These ideas are related by the epigraph of $f$:
$$
\epi{f}:=\set{(x,\mu)}{f(x)\le\mu}\subset\R^n\times \R.
$$
A function $\map{g}{\R^n}{\R\cup\{\infty\}}$ is convex if and only if $\epi{g}$ is a convex set.
A function $f$ is called {\it closed} if $\epi{f}$ is a closed set, or equivalently, 
if $f$  is {\it lower semicontinuous} (lsc).

Facts about convex sets can be translated into facts about convex functions. 
The reverse is also true with the aid of the convex indicator functions:
\begin{equation}
\label{eq:indicator}
\indicator{x}{\mathcal{C}} := 
\begin{cases} 0 & \mbox{if} \quad x \in \mathcal{C} \\
\infty & \mbox{else}. 
\end{cases}
\end{equation}
Examples of convex sets include subspaces and their translates (affine sets) 
as well as the lower level sets
of convex functions:
$$
\lev{f}{\tau}:=\set{x}{f(x)\le \tau}.
$$
Just as with closed sets, the intersection of an arbitrary collection of convex sets is also convex. 
For this reason
we define the convex hull of a set $\calE$ to be the the intersection of all convex sets that contain it, 
denoted by $\conv{\calE}$.

The convex sets of greatest interest to us are the convex polyhedra,
$$
\mathcal{W}:=\set{x}{H^Tx\le h}\quad\mbox{ for some $H\in\R^{n\times m}$ and $h\in\R^m$,}
$$
while the convex functions of greatest interest are the 
piecewise linear-quadratic (PLQ) penalties, shown in Figs.~\ref{fig:quadratic}-\ref{fig:enet}. 
As discussed in Section~\ref{sec:stats}, these penalties allow 
us to model impulsive disturbances in the process (see Figs.~\ref{fig:1norm} and~\ref{fig:enet}),
to develop robust distributions for measurements (see Fig.~\ref{fig:Huber}) 
and implement support vector regression (SVR) in the context of dynamic systems (see fig.~\ref{fig:Vapnik}).

\subsection{Smooth case: first- and second-order methods}

Consider the problem 
\[
\min_x f(x),
\]
together with an iterative procedure indexed by $\iterc$ that is initialized at $x^1$. 
When $f$ is a $C^1$-smooth function with $\beta$-Lipschitz continuous gradient, 
i.e. {\it $\beta$-smooth}: 
\begin{equation}
\label{eq:beta-smoth}
\|\nabla f(x) - \nabla f(y)\| \leq \beta\|x-y\|,
\end{equation}
 $f$ admits the upper bounding quadratic model
\begin{equation}\label{eq:quad upper bd}
f(x) \leq m_\iterc(x):=f(x^\iterc)+\langle\nabla f(x^\iterc),x-x^\iterc\rangle+\frac{\beta}{2}\|x-x^\iterc\|^2.
\end{equation}
If we minimize $m_\iterc(x)$ to obtain $x^{\iterc +1}$, this gives the iteration 
\[
x^{\iterc +1} := x^\iterc - \frac{1}{\beta}\nabla f(x^\iterc), 
\]
or steepest descent. The upper bound \eqref{eq:quad upper bd} shows we have strict descent:
$$
\begin{aligned}
f(x^{\iterc +1})&\leq f(x^\iterc)-\langle \nabla f(x^\iterc),\beta^{-1}\nabla f(x^\iterc)\rangle +\frac{\beta}{2}\|\beta^{-1}\nabla f(x^\iterc)\|^2 \\
&=f(x^\iterc) -\frac{\|\nabla f(x^\iterc)\|^2}{2\beta}.
\end{aligned}
$$
If, in addition, $f$ is convex, and a minimizer $x^*$ exists, we obtain 
\[
\begin{aligned}
f(x^\iterc)  -\frac{\|\nabla f(x^\iterc)\|^2}{2\beta} & 
\leq f^* + \langle \nabla f(x^\iterc), x^\iterc - x^*\rangle -\frac{\|\nabla f(x^\iterc)\|^2}{2\beta}\\
& = f^* + \frac{\beta}{2} \left(\|x^\iterc - x^*\|^2 - \|x^{\iterc + 1} - x^*\|^2\right),
\end{aligned}
\]
where $f^*  = f(x^*)$ is the same at any minimizer by convexity. 

Adding up, we get an $O\left(\frac{1}{\iterc}\right)$ convergence rate on function values: 
\[
f(x^{\iterc})-f^* \le \frac{\beta \|x^1-x^*\|^2}{2\iterc}.
\]

For the least squares Kalman smoothing problem~\eqref{fullLS}, we also know that $f$ is {\it $\alpha$-strongly convex}, i.e. 
$f(x) - \frac{\alpha}{2} \|x\|^2$ is convex with $\alpha \geq 0$.   
Strong convexity  can be used to obtain a much better rate for steepest descent: 
\[
f(x^\iterc) - f^{*}  \leq \frac{\beta}{2}(1- \alpha/\beta)^{\iterc}\|x^1 - x^{*}\|^2. 
 \] 
Note that  $0 \leq \frac{\alpha}{\beta}\leq 1$.   

When minimizing a strongly convex function, the minimizer $x^*$ is unique, and we can also obtain 
a rate on the squared distance between $x^\iterc$ and $x^*$: 
\[
\|x^\iterc - x^{*}\|^2 \leq (1- \alpha/\beta)^{\iterc}\|x^1 - x^{*}\|^2.
\]
These rates can be further improved by considering {\it accelerated-gradient} methods (see e.g.~\cite{nest_lect_intro}) which achieve the much faster rate $(1-\sqrt{\alpha/\beta})^\iterc$. 

Each iteration of steepest descent in the classic least squares formulation~\eqref{fullLS} of the 
Kalman smoothing problem gives a 
fractional reduction in both function value and distance to optimal solution. In this case, computing the gradient 
requires only matrix-vector products, which require $O(Nn^2)$ arithmetic operations.   
Thus, either gradient descent or conjugate gradient (which has the same rate as accelerated gradient 
methods in the least squares case) is a reasonable option if $n$ is very large.

The solution to (12) can also be obtained by solving a linear system using $O(Nn^3)$ arithmetic operations,
since~\eqref{smoothingSol} is block-tridiagonal positive definite. This complexity is tractable for moderate state-space dimension $n$.
The approach is equivalent to a single iteration on the full quadratic model of the Newton's method, discussed below. 

Consider the problem of minimizing a $C^2$-smooth function $f$. 
Finding a critical point $x$ of $f$ can be recast as the problem of solving the nonlinear equation $\nabla f(x)=0$. 
For a smooth function $\map{G}{\Rn}{\Rn}$, Newton's method is designed to locate solutions to
the equation $G(x)=0$.
Given a current iterate $x^\iterc$, Newton's method linearizes $G$ at $x^\iterc$ and solves the equation 
$G(x^\iterc)+\nabla G(x^\iterc)(y-x^\iterc)=0$ for $y$. 
Provided that $\nabla G(x^\iterc)$ is invertible, the Newton iterate is given by 
\begin{equation}\label{eq:Newton iter}
x^{\iterc+1}:=x^\iterc-[\nabla G(x^\iterc)]^{-1}G(x^\iterc).
\end{equation}
When $G:=\nabla f$, the Newton iterate \eqref{eq:Newton iter} is the unique critical point of the best quadratic  
approximation of $f$ at $x^\iterc$, namely
\begin{eqnarray*}
\qquad \qquad Q(x^\iterc;y)&:=& f(x^\iterc)+\langle \nabla f(x^\iterc),y-x^\iterc\rangle\\
&+&\frac{1}{2}\langle \nabla^2 f(x^\iterc)(y-x^\iterc),y-x^\iterc\rangle,
\end{eqnarray*}
provided that the Hessian $\nabla^2 f(x^\iterc)$ is invertible. 

If $G$ is a $C^1$-smooth function with $\beta$-Lipschitz Jacobian $\nabla G$ that is locally invertible for all $x$
near a point $x^*$ with $G(x^*) = 0$, 
then near $x^*$ the Newton iterates \eqref{eq:Newton iter}
satisfy
$$\|x^{\iterc +1}-x^*\|\leq \frac{\beta}{2}\|\nabla G(x^\iterc)^{-1}\|\cdot\|x^\iterc-x^*\|^2.$$

Once we are {\it close enough} to a solution, Newton's method gives a {\it quadratic} rate of convergence.  
Consequently, locally the number of correct digits double for each iteration. Although the solution 
may not be obtained in one step (as in the quadratic case), only a few iterations are required to converge 
to machine precision.

In the remainder of the section, we generalize steepest descent and Newton's methods 
to nonsmooth problems of type~\eqref{fullGen}. 
In Section~\ref{sec:sub-gradient}, we describe the {\it sub-gradient} descent method, 
and show that it converges very slowly.  
In Section~\ref{sec:prox-gradient}, we describe the 
proximity operator and proximal-gradient methods, which are applicable when working with separable nonsmooth 
terms in~\eqref{fullGen}. In Section~\ref{sec:splitting}, we show how to solve 
more general nonsmooth problems~\eqref{fullGen} using splitting techniques, including ADMM and Chambolle-Pock iterations. 
Finally, in Section~\ref{sec:IP}, we show how second-order interior point
methods can be brought to bear on all problems of interest of type~\eqref{fullGen}.

\subsection{Nonsmooth case: subgradient descent}
\label{sec:sub-gradient}

Given a convex function $f$, a vector $v$ is a {\it subgradient} of $f$ at a point $x$
if 
\begin{equation}
\label{eq:sgineq}
f(y) \geq f(x) + \langle v, y-x\rangle \quad \forall y.
\end{equation}
The set of all subgradients at $x$ is called the {\it subdifferential}, and is denoted by $\partial f(x)$. 
Subgradients generalize the notion of gradient; in particular, $\partial f(x) = \{v\} \iff v = \nabla f(x)$~\cite{RTR}. A more comprehensive discussion of the subdifferential is presented in Appendix A2.

Consider the absolute value function shown in Figure~\ref{fig:1norm}. It is differentiable at all points except for $x = 0$, and so the subdifferential is precisely the gradient for all $x \neq 0$. The subgradients at $x=0$ are the slopes of lines passing through the origin and lying below the graph of the absolute value function. Therefore, $\partial |\cdot|(0) = [-1,1]$.

Consider the following simple algorithm for minimizing a Lipschitz continuous (but nonsmooth) convex $f$. 
Given an oracle that delivers some $v^\iterc\in \partial f(x^\iterc)$, set
\begin{equation}
\label{eq:sub}
x^{\iterc +1} := x^\iterc - \alpha_\iterc  v^\iterc,
\end{equation}
for a judiciously chosen stepsize $\alpha_\iterc $. Suppose we are minimizing $|x|$ and start at $x=0$, the global minimum. The oracle could return any value $v\in[-1,1]$, and so we will move away from $0$ when using~\eqref{eq:sub}! 
In general, the function value need not decrease at each iteration, and we see that 
$\alpha_\iterc $ must decrease to $0$ for any hope of convergence. 
On the other hand, if $\sum_\iterc  \alpha_\iterc   = R < \infty$, we can never reach $x^*$
if $\|x^1 - x^*\| > R$, where $x^1$ is the initial point and $x^*$ the minimizer. 
Therefore, we also must have $\sum_\iterc \alpha_\iterc = \infty$.\\
Setting $l_\iterc := f(x^\iterc) + \langle v^\iterc, x^* - x^\iterc\rangle$, by~\eqref{eq:sgineq} we have $l_k \leq f(x^*) \leq f(x^\iterc)$ for 
$v\in \partial f(x^\iterc)$. The subgradient method closes the gap between $l_\iterc$ and $f(x^\iterc)$.  
The Liptschitz continuity of $f$ implies that $\|v^\iterc\| \leq L$, and so, by \eqref{eq:sgineq}, 
\[
\begin{aligned}
0\!\le\! \|x^{\iterc +1} \!-\! x^*\|^2 \!&=\! \|x^\iterc \!-\! x^*\|^2 \!+\! 2\alpha_\iterc  \langle v^\iterc, x^* \!-\! x^\iterc  \rangle \!+\! \alpha_\iterc ^2 \|v^\iterc\|^2 \\
\!& \leq\! \|x^1 \!-\! x^*\|^2 \!+\! \sum_{i=1}^\iterc 2\alpha_i \langle v^i, x^* \!-\! x^i \rangle  \!+\! L^2 \sum_{i=1}^\iterc \alpha_i^2\\
\!& =\! \|x^1 \!-\! x^*\|^2 \!+\! \sum_{i=1}^\iterc 2\alpha_i(l_i\!-\!f(x^i))\!+\! L^2 \sum_{i=1}^\iterc \alpha_i^2.
\end{aligned}
\]
Rewriting this inequality gives 
\begin{equation}
  \label{eq:bound}
\begin{aligned}
0 \leq \min_{i = 1, \dots, \iterc} (f(x^i) - l_i) & \leq \sum_{i=1}^\iterc \frac{\alpha_i}{\sum_{i=1}^\iterc \alpha_i}(f(x^i) - l_i ) \\
& \leq \frac{\|x^1 - x^*\|^2  + L^2 \sum_{i=1}^\iterc \alpha_i^2}{2\sum_{i=1}^\iterc \alpha_i}.
\end{aligned}
\end{equation}
In particular, if $\{\alpha_\iterc\}$ are square summable but not summable, convergence 
of $\displaystyle \min_{i = 1, \dots, \iterc} \{f(x^i) - l_i\} $ to $0$ is guaranteed.  
But there is a fundamental limitation of the subgradient method. In fact, suppose 
that we know $\|x^1 - x^*\|$, and want to choose steps $\alpha_i$ to minimize the gap in $t$ iterations. 
By minimizing the right hand side of~\eqref{eq:bound}, we find that the optimal step sizes (with respect 
to the error bound) are 
\[
\alpha_i = \frac{\|x^1 - x^*\|}{L\sqrt{\iterc}}.
\]
Plugging these back in, and defining $f_{\mbox{best}}^\iterc = \min_{i=1,\dots, \iterc} f(x^i)$, we have  
\[
f_{\mbox{best}}^\iterc  - f^* \leq \frac{\|x^1 - x^*\|L}{\sqrt{\iterc}}.
\]
Consequently, the best provable subgradient descent method is extremely slow. This rate 
can be significantly improved by exploiting the structure of the nonsmoothness in $f$.

\subsection{Proximal gradient methods and accelerations}
\label{sec:prox-gradient}
For many convex functions, and in particular for a range of general smoothing formulations~\eqref{fullGen}, 
we can design algorithms that are much faster than $O(1/\sqrt{\iterc})$. 
Suppose we want to minimize the sum 
\[
f(x) + g(x),
\]
where $f$ is convex and $\beta$-smooth~\eqref{eq:beta-smoth}, while $g$ is any convex function. Using the bounding 
model~\eqref{eq:quad upper bd} for $f$, we can get a global upper bound for the sum:
\begin{eqnarray*}
&&  f(x) + g(x) \leq m_\iterc(x)\\
&& m_\iterc(x) := f(x^\iterc)+\langle\nabla f(x^\iterc),x-x^\iterc\rangle+\frac{\beta}{2}\|x-x^\iterc\|^2 + g(x).
\end{eqnarray*}
We immediately see that setting 
\begin{equation}
\label{eq:pstep}
x^{\iterc +1} := \arg\min_{x}\; m_\iterc(x) 
\end{equation}
 ensures descent for $f+g$, since 
\[
f(x^{\iterc +1}) + g(x^{\iterc +1}) \leq m_\iterc(x^{\iterc +1}) \leq m_\iterc(x^{\iterc}) = f(x^\iterc) + g(x^\iterc).
\] 
One can check that $m_\iterc(x^{\iterc +1})= m_\iterc(x^\iterc)$ if and only if $x^\iterc$ is a global minimum of $f+g$. 
Rewriting~\eqref{eq:pstep} as 
\[
x^{\iterc +1} := \arg\min_{x} \beta^{-1} g(x) + \frac{1}{2}\|x - (x^\iterc- \frac{1}{\beta}\nabla f(x^\iterc))\|^2,
\]
and define the {\it proximity} operator for $\eta g$ \cite{BauCom} by
\begin{equation}
\label{eq:proxOp}
\prox_{\eta g}(y) := \arg\min_{x} \eta g(x) + \frac{1}{2}\|x - y\|^2.
\end{equation}
We see that~\eqref{eq:pstep} is precisely the proximal gradient method: 
\begin{equation}
\label{eq:pgiter}
x^{\iterc +1} := \prox_{\beta^{-1}g}\left(x^\iterc - \frac{1}{\beta}\nabla f(x^\iterc)\right).
\end{equation}

The proximal gradient iteration~\eqref{eq:pgiter} 
converges with the same rate as gradient descent, in particular with rate 
$O\left(1/\iterc\right)$ for convex functions and $O\left((1-\alpha/\beta)^\iterc\right)$ 
for $\alpha$-strongly convex functions. These rates are in a completely different 
class than the $O\left(1/\sqrt{\iterc}\right)$ rate obtained by  
the subgradient method, since they exploit the 
additive structure of $f+g$. 
Proximal gradient algorithms can also be accelerated, 
achieving rates of $O\left(1/\iterc^{2}\right)$  and $O\left((1-\sqrt{\alpha/\beta})^\iterc\right)$ 
respectively, using techniques from~\cite{nest_lect_intro}. 

In order to implement~\eqref{eq:pgiter}, we must be able to efficiently compute the proximity 
operator for $\eta g$. 
For many nonsmooth functions $g$, this operator can be computed in $O(n)$ or $O(n\log n)$ time. 
An important example is the {\it convex indicator} function~\eqref{eq:indicator}. 
In this case, the prox-operator is the projection operator:
\begin{equation}
\label{eq:proxproj}
\begin{aligned}
\prox_{\eta \indicator{x}{\mathcal{C}}}(y) &=&\indicator{x}{\mathcal{C}} +  \min_x \frac{1}{2}\|x - y\|^2 \\
&=& \min_{x\in \mathcal{C}} \frac{1}{2}\|x - y\|^2 = \proj_{\mathcal{C}} (y).
\end{aligned}
\end{equation}
In particular, when minimizing $f$ over a convex set $\mathcal{C}$, iteration~\eqref{eq:pgiter} 
recovers the {\it projected gradient} method if we choose $g(x) = \indicator{x}{\mathcal{C}}$.

Many examples and identities useful for computing proximal operators are collected in~\cite{combettes2011proximal}. One important example is the Moreau identity (see e.g.~\cite{RTRW}):
\begin{equation}
\label{eq:moreau}
\prox_{f}(y) + \prox_{f^*}(y) = y\ . 
\end{equation}
Here, $f^*$ denotes the {\it convex conjugate} of $f$:
\begin{equation}
\label{eq:conjIntro}
f^*(\omega) := \sup_{y}( \langle y,\omega \rangle - f(y)),
\end{equation}
whose properties are explained in Appendix A2, 
in the context of convex duality. 
Identity~\eqref{eq:moreau} shows that the prox of $f$ can be used to compute the prox of $f^*$, and vice versa. 

{\bf Example: proximity operator for the $\ell_1$-norm}. 
Consider the example $g(x) =  \|x\|_1$, often used in applications to induce sparsity of $x$. 
The proximity operator of this function is can be computed  by reducing to 
the 1-dimensional setting and considering cases. Here, we show how to compute it using~\eqref{eq:moreau}:
\[
\prox_{ \eta \|\cdot\|_1}(y) = y - \prox_{(\eta \|\cdot\|_1)^*} (y).
\]
The convex conjugate of the scaled 1-norm is given by 
\begin{equation}
\label{eq:ell1conj}
(\eta \|\cdot\|_1)^*(\omega) = \sup_{x} \langle x, \omega\rangle - \eta \|x\|_1 
= 
\begin{cases} 0 & \mbox{if} \quad \|\omega\|_\infty \leq \eta \\
\infty & \mbox{otherwise} 
\end{cases},
\end{equation} 
which is precisely the indicator function of $\eta\mathbb{B}_{\infty}$, the scaled $\infty$-norm unit ball. 
As previously observed, the proximity operator for and indicator function is the projection. 
Consequently, the identity~\eqref{eq:moreau} simplifies to 
\[
\prox_{ \eta \|\cdot\|_1}(y) = y - \proj_{\eta \mathbb{B}_\infty}(y)
\]
whose $i$th element is given by 
\begin{equation}
\label{eq:proxell1}
\prox_{ \eta \|\cdot\|_1}(y)_i = 
\begin{cases}
y_i -y_i = 0 & \mbox{if} \quad |y_i| \leq \eta \\
y_i - \eta\mbox{sign}(y_i) & \mbox{if} \quad |y_i| > \eta\,
\end{cases}
\end{equation}
which corresponds to {\it soft-thresholding}. 
Computing the proximal operator for the 1-norm and projection onto the $\infty$-norm ball both 
require $O(n)$ operations. Projection onto the 1-norm ball $\mathbb{B}_1$ can be implemented using a sort,
and so takes $O(n\log(n))$ operations, see e.g.~\cite{BergFriedlander:2008}.  
\begin{flushright}
$\blacksquare$
\end{flushright}

To illustrate the method in the context of Kalman smoothing, consider taking the general formulation~\eqref{fullGen} with $V$ and $J$ both smooth,
$\gamma = 1$, and $x \in \tau\mathbb{B}$ any norm-ball for which we have a fast projection 
(common cases are $2$-norm, $1$-norm, or $\infty$-norm):
$$
\min_{x\in\tau\mathbb{B}} V(R^{-1/2}(y - Cx))+ J(Q^{-1/2}(z - Ax )). 
$$

\begin{algorithm}
\caption{Proximal Gradient for Kalman Smoothing, $J$ and $V$ Huber or quadratic}\label{alg:ProxGrad}
\begin{enumerate}

\item
Initialize $x^1 = 0$, compute $ d^1 = \nabla f(x^1)$.
Let $\beta = \|C^TR^{-1}C + A^TQ^{-1}A\|_2$.

\noindent
\item While $\|\prox_{g}(x^\iterc  - d^\iterc )\| > \epsilon$

\begin{itemize}

\item Set $\iterc = \iterc+1$.

\item
update $x^{\iterc} = \prox_{\beta^{-1} g}(x^{\iterc-1} - \beta^{-1} d^{\iterc-1})$.

\item
Compute $d^\iterc  = \nabla f(x^\iterc )$.

\end{itemize}
\item Output $x^\iterc $.

\end{enumerate}
\end{algorithm}
\begin{algorithm}
\caption{FISTA for Kalman Smoothing, $J$ and $V$ Huber or quadratic}\label{alg:FISTA}
\begin{enumerate}

\item
Initialize $x^1 = 0$, $s_1 = 1$, compute $ d^1 = \nabla f(x^1)$.
Let $\beta = \|C^TR^{-1}C + A^TQ^{-1}A\|_2$.

\noindent
\item While $\|\prox_{g}(\omega^\iterc  - d^\iterc )\| > \epsilon$

\begin{itemize}

\item Set $\iterc = \iterc+1$.

\item
update $x^{\iterc} = \prox_{\beta^{-1} g}(\omega^{\iterc-1} - \alpha d^{\iterc-1})$.

\item
set $s_\iterc = \frac{1 + \sqrt{1 + 4 s_{\iterc-1}^2}}{2}$

\item 
set $\omega^\iterc  = x^\iterc  + \frac{s_{\iterc-1} - 1}{s_\iterc}(x_\iterc - x_{\iterc-1})$.

\item
Compute $g^\iterc  = \nabla f(x^\iterc )$.

\end{itemize}
\item Output $\omega^\iterc $.

\end{enumerate}
\end{algorithm}

The gradient for the system is given by 
\begin{eqnarray*}
\nabla f(x) &=& C^TR^{-1/2}\nabla V(R^{-1/2}(Cx - y))\\ 
&& +A^TQ^{-1/2}\nabla J(Q^{-1/2}Ax - z)).
\end{eqnarray*}
When $V$ and $J$ are quadratic or Huber penalties, the Lipschitz constant
$\beta$ of $\nabla f$ is bounded by the largest singular value of $C^TR^{-1}C + A^TQ^{-1}A$, which we can obtain using power iterations. This system is block tridiagonal, so matrix-vector multiplications are far more efficient than for 
general systems. Specifically, for Kalman smoothing, the systems $C, Q, R$ are block diagonal, while $A$ is block bidiagonal. As a result, products with $A, A^T, C, Q^{-1/2}, R^{-1/2}$ can all be computed using $O(Nn^2)$ arithmetic operations, rather than $O(N^2n^2)$ operations
 as for a general system of the same size. A simple proximal gradient method is given by Algorithm~\ref{alg:ProxGrad}.
Note that soft thresholding for Kalman smoothing has complexity $O(nN)$, while e.g. projecting
onto the 1-norm ball has complexity $O(nN\log(nN))$. 
Therefore the $O(n^2N)$ cost of computing the gradient $\nabla f(x^\iterc )$ is dominant.

Algorithm~\ref{alg:ProxGrad} has at worst 
$O\left(\iterc^{-1}\right)$ rate of convergence. If $J$ is taken to be a quadratic, 
$f$ is strongly convex, in which case
we achieve the much faster rate $O\left((1-\alpha/\beta)^\iterc \right)$.\\
Algorithm~\ref{alg:FISTA} illustrates
the FISTA scheme~\cite{beck2009fast} applied to Kalman smoothing. 
This acceleration uses two previous iterates rather than just one, and
achieves a worst case rate of $O\left(\iterc^{-2}\right)$.
This can be further improved to $O\left((1-\sqrt{\alpha/\beta})^\iterc \right)$
when $J$ is a convex quadratic using techniques in~\cite{nest_lect_intro}, 
or periodic restarts of the step-size sequence $s_\kappa$.

\subsection{Splitting methods}
\label{sec:splitting}

Not all smoothing formulations~\eqref{fullGen} are the sum of a smooth function and
a separable nonsmooth function. In many cases, the composition of a nonsmooth penalty with a general
linear operator can preclude the approach of the previous section; for example, 
the {\it robust} Kalman smoothing problem in~\cite{aravkin2011laplace}:
\begin{equation}
\label{eq:KalmanL1}
\min_{x} \| R^{-1/2}(y - Cx)\|_1 + \frac{1}{2} \|Q^{-1/2}(z - Ax )\|^2. 
\end{equation}
Replacing the quadratic penalty with the 1-norm allows the development of a robust smoother
when a portion of (isolated) measurements
are contaminated by outliers. The composition of the nonsmooth 1-norm with a general linear form
makes it impractical to use the proximal gradient method 
since the evaluation of the prox operator
\[
\prox_{\eta \|y - C(\cdot)\|_1} (y) = \arg\min_{x} \frac{1}{2}\|y-x\|^2 + \eta \|y - Cx\|_1
\]
requires an iterative solution scheme for general $C$. However, it is possible to design a primal-dual method using a range of strategies known as splitting methods.
Convex duality theory and related concepts are explained in Appendix A2.\\

A well-known splitting method, popularized by~\cite{boyd2011distributed}, is  
the Alternating Direction Method of Multipliers (ADMM), which is equivalent to 
Douglas-Rachford splitting on an appropriate dual problem~\cite{lions1979splitting}. 
The ADMM scheme is applicable to general problems of type 
\begin{equation}
\label{eq:ADMMtemplate}
\min_{x,\omega} f(x) + g(\omega) \quad \mbox{s.t.} \quad K_1x + K_2\omega = c.
\end{equation}
A fast way to derive the approach is to consider the Augmented Lagrangian \cite{RTR74}
dualizing the equality constraint in~\eqref{eq:ADMMtemplate}:
\[\begin{aligned}
\mathcal{L}(x,\omega, u,\tau) \!:=\! f(x) \!+\! g(\omega) \!&+\! u^T(K_1x \!+\! K_2\omega \!-\! c) \\ 
\!&+\! \frac{\tau}{2}\|K_1x \!+\! K_2\omega \!-\! c\|^2,
\end{aligned}
\]
where $\tau > 0$.
The ADMM method proceeds by using alternating minimization of $\mathcal{L}$ in $x$
and $\omega$ with appropriate dual updates (which is equivalent to the Douglas-Rachford method on the dual of~\eqref{eq:ADMMtemplate}. The iterations are explained fully in Algorithm~\ref{alg:admm_gen}.

\begin{algorithm}
\caption{ADMM algorithm for~\eqref{eq:ADMMtemplate}}\label{alg:admm_gen}
\begin{enumerate}
\item
Input $x^1, \omega^0 \neq \omega^1$. Input $\tau > 0$,  $\epsilon$. 
\item While $\|K_1x^\iterc  + K_2\omega^\iterc  - c\| > \epsilon$ and \\
$\|\tau K_1^TK_2(\omega^{\iterc +1} -\omega^\iterc )\|> \epsilon$
\begin{itemize}
\item Set $\iterc := \iterc+1$.
\item
update  
\small
$$x^{\iterc +1} := \arg\min_x \left\{
\begin{aligned}f(x) &+ (u^\iterc)^TK_1x \\ &+ \frac{\tau}{2}\|K_1x + K_2\omega^\iterc - c\|^2
\end{aligned}
\right\}$$
 \normalsize
\item 
update 
\small 
$$\omega^{\iterc +1} := \arg\min_\omega \left\{ 
\begin{aligned}g(\omega) &+ (u^\iterc)^TK_2\omega \\
&+  \frac{\tau}{2}\|K_1x^{\iterc +1} + K_2\omega - c\|^2 
\end{aligned}
\right\}$$
\normalsize 
\item 
update $u^{\iterc +1} := u^\iterc + \tau(K_1x^{\iterc +1} + K_2\omega^{\iterc +1} - c)$
\end{itemize}
\item Output $(x^\iterc,\omega^\iterc) $.
\end{enumerate}
\end{algorithm}

ADMM has convergence rate $O(1/\iterc)$, but can be accelerated 
under sufficient regularity conditions (see e.g.~\cite{davis2015three}). 
For the Laplace $\ell_1$ smoother \eqref{eq:KalmanL1}, the transformation to template~\eqref{eq:ADMMtemplate}
is given by 
\begin{equation}
\label{eq:KalmanL1Ref}
\min_{x,\omega}\! \left\{ \|\omega\|_1\! +\! \frac{1}{2} \|Q^{-1/2}(z - Ax)\|^2 \, \left| \, \omega \!+\!  R^{-1/2} Cx = R^{-1/2} y  \right.\right\}.
\end{equation}
ADMM specialized to~\eqref{eq:KalmanL1Ref} is given by Algorithm~\ref{alg:admm}.

\begin{algorithm}
\caption{ADMM algorithm for~\eqref{eq:KalmanL1Ref}}\label{alg:admm}
\begin{enumerate}
\item
Input $x^1, \omega^0 \neq \omega^1$. Input $\tau > 0$,  $\epsilon$. 
\item While $\|\omega^\iterc  +  R^{-1/2} Cx^\iterc  - R^{-1/2} y\| > \epsilon$ and\\
 $\|\tau C^TR^{-T/2}(\omega^{\iterc +1} -\omega^\iterc )\|> \epsilon$
\begin{itemize}
\item Set $\iterc := \iterc+1$.
\item
update
\small
\begin{eqnarray*}
x^{\iterc +1} := \arg\min_x \frac{1}{2} \|Q^{-1/2}(z - Ax) \|^2 + x^Tu^\iterc \\ + 
\frac{\tau}{2}\| R^{-1/2}(C x-y)+ \omega^\iterc\|^2
\end{eqnarray*}
 \normalsize
\item 
update
\small 
\begin{eqnarray*}
\omega^{\iterc +1} \!:=\! \arg\min_\omega  \|\omega\|_1 \!+\!  \frac{\tau}{2}\left\|\omega \!+\! u^\iterc/\tau 
+ R^{-1/2} (C x^{\iterc +1}  - y)\right\|^2
 \end{eqnarray*}
 \normalsize
\item 
update $$u^{\iterc +1} := u^\iterc + \tau(R^{-1/2} C x^{\iterc +1} +  \omega^{\iterc +1} - R^{-1/2} y)$$
\end{itemize}
\item Output $x^\iterc $.
\end{enumerate}
\end{algorithm}

We make two observations. First, note that the $x$-update requires solving a least squares problem, 
in particular inverting $A^TQ^{-1}A + C^TR^{-1}C$. 
Fortunately, in problem~\eqref{eq:KalmanL1Ref} this system does not change between iterations, 
and can be factorized once in $O(n^3N)$ arithmetic operations and stored. 
Each iteration of the $x$-update can be obtained in $O(n^2N)$ arithmetic operations which has the same 
complexity as a matrix-vector product. Splitting schemes that
avoid factorizations are described below. 
However, avoiding factorizations is not always the best strategy since
the choice of splitting scheme can have a dramatic effect on the performance.
Performance differences between various splitting are explored in the numerical
section.  Second, the $\omega$-update has a convenient closed form representation
in terms of the proximity operator~\eqref{eq:proxOp}:
\[
\omega^{\iterc +1} := \prox_{\tau^{-1} \|\cdot\|_1} (u^\iterc/\tau + R^{-1/2} (C x^{\iterc +1}  - y)).
\]
The overall complexity of each iteration of the ADMM $\ell_1$-Kalman smoother is $O(n^2N)$, after 
the initial $O(n^3N)$ investment to factorize $A^TQ^{-1}A + C^TR^{-1}C$. 

There are several types of splitting schemes, including Forward-Backward~\cite{passty1979ergodic}, Peaceman-Rachford~\cite{lions1979splitting}, and others.  A survey of these algorithms is beyond the scope of this paper. 
See~\cite{BauCom,davis2015three} for a discussion of splitting methods and the relationships 
between them. See also~\cite{davis2014faster}, for a detailed analysis of convergence rates of several splitting schemes under regularity assumptions.\\
We are not aware of a detailed study or comparison of these techniques for general Kalman smoothing problems, 
and future work in this direction can have a significant impact in the community.   
To give an illustration of the numerical behavior and variety of splitting algorithms, we 
present the algorithm of Chambolle-Pock~(CP)~\cite{chambolle2011first}, for convex problems of type  
\begin{equation}
\label{eq:cpForm}
\min_{x} f(Kx) + g(x),
\end{equation}
where $f$ and $g$ are convex functions with computable proximity operators, 
while $L$ is the largest singular value of $K$. 
The CP iteration is specified in Algorithm~\ref{alg:CP}. 
\begin{algorithm}
\caption{Chambolle-Pock algorithm for~\eqref{eq:cpForm}}\label{alg:CP}
\begin{enumerate}
\item
Input $x^0 \neq x^1, \omega^0 \neq \omega^1$. Input $\tau, \sigma$ s.t. $\tau\sigma L^2 <1$. Input $\epsilon$.
\item While $(\|\omega^{\iterc +1} - \omega^\iterc \| + \|x^{\iterc +1} - x^\iterc \|> \epsilon)$ 
\begin{itemize}
\item Set $\iterc = \iterc+1$.
\item
update $\omega^{\iterc +1} = \prox_{\sigma f^*} (\omega^\iterc  + \sigma K (2 x^{\iterc} - x^{\iterc-1}))$
\item 
update $x^{\iterc +1} = \prox_{\tau g}(x^\iterc  - \tau K^T \omega^{\iterc +1})$

\end{itemize}
\item Output $x^\iterc $.
\end{enumerate}
\end{algorithm}

Algorithm~\ref{alg:CP} requires only the proximal operators for $f^*$ and $g$ to be implementable. 
Like ADMM, it has a convergence 
rate of $O(1/\iterc)$, and can be accelerated to $O(1/\iterc^2)$ under specific regularity assumptions. 
When $g$ is strongly convex, one such acceleration is presented in \cite{chambolle2011first}.\\

There are multiple ways to apply the CP scheme to a given Kalman smoothing formulation.  
Some schemes allow CP to solve large-scale  
smoothing problems~\eqref{fullGen} using 
only matrix-vector products,  avoiding large-scale matrix solves entirely. 
However, this may not be the best approach, as we show in our numerical study 
in the following section. 
General splitting schemes such as Chambolle-Pock  
can achieve at best $O(1/\iterc^2)$ convergence rate 
for general nonsmooth Kalman formulations. Faster rates require much stronger assumptions, 
e.g. smoothness of the primal or dual problems~\cite{chambolle2011first}. 
When these conditions are present, the methods can be remarkably efficient.

\subsection{Formulations Using Piecewise Linear Quadratic (PLQ) Penalties \cite{RTRW}}
\label{sec:IP}

When the state size $n$ is moderate, so that $O(n^3N)$ is an acceptable cost to pay, 
we can obtain very general and fast methods for Kalman smoothing systems.
We recover {\it second-order behavior} and fast local convergence rates by 
developing interior point methods for the entire class~\eqref{fullGen}.
These methods can be developed for any piecewise linear quadratic $V$ and $J$, and
allow the inclusion of polyhedral constraints that link adjacent time points. This can be accomplished 
using $O(n^3N)$ arithmetic operations, the same complexity as solving the least squares Kalman smoother. 

To see how to develop second-order interior point methods for these PLQ smoothers, we first define 
the general PLQ family and consider its conjugate representation and optimality conditions.

\begin{definition}[PLQ functions and penalties]
\label{generalPLQ}A piecewise linear quadratic (PLQ) function is any function 
$\map{\rho(c, C, b, B, M; \cdot)}{\R^n}{\R\cup\{\infty\}}$ 
admitting representation
\begin{equation}
 \label{PLQpenalty}
\begin{aligned}
\rho(c, C, b, B, M; x) 
&:=
\sup_{v \in\calV}
\left\{ \langle v,b + Bx \rangle - \half\langle v, Mv\rangle \right\} \\
& = \left(\half\|\cdot\|_M^2 + \indicator{\cdot}{\calV}\right)^*(b + Bx)\;,
\end{aligned}
\end{equation}
where $\calV$ is the polyhedral set specified by $H\in\mathbb{R}^{k \times \ell}$ and $h\in\mathbb{R}^\ell$ as follows
$$
\calV = \{v: H^Tv \leq h\}\; ,
$$
 $M\in \Skp$ the set of real symmetric positive semidefinite matrices,
$b + Bx$ is an injective affine transformation in $x$, with $B\in\mB{R}^{k\times n}$, 
so, in particular, $n \leq k$ and $null(B) = \{0\}$. If $0\in \calV$, 
then the PLQ is necessarily non-negative and hence represents a {\it penalty}. 
\end{definition}
The last equation in~\eqref{PLQpenalty} 
is seen immediately using~\eqref{eq:conjIntro}. 
In what follows we reserve the symbol $\rho$ for a PLQ penalty often writing $\rho(x)$
and suppressing the litany
of parameters that precisely define the function. When detailed knowledge of these parameters is
required, they will be specified.\\
Below we show how the six loss functions illustrated in 
Figure \ref{fig:quadratic}-\ref{fig:enet} can be represented as members of
the PLQ class. In each case, the verification of the representation is straightforward.
These dual (conjugate) representations facilitate the general optimization approach. 

{\bf Examples of scalar PLQ} 

\begin{enumerate}

\item {\bf quadratic} ($\ell_2$) penalty, Fig.~\ref{fig:quadratic}: 
\[
 \sup_{v \in \mathbb{R} }\left\{vx - \frac{1}{2} v^2 \right\}    
\]
\item {\bf  absolute value} ($\ell_1$) penalty, Fig.~\ref{fig:1norm}:
\[
 \sup_{v\in [-1,1] }\left\{vx \right\}   
\]
\item {\bf Huber} penalty, Fig.~\ref{fig:Huber}: 
\[
 \sup_{v\in [-\kappa, \kappa]}\left\{vx - \frac{1}{2} v^2 \right\}  
\]
\item {\bf Vapnik} penalty, Fig~\ref{fig:Vapnik}:
\[
 \sup_{v \in \left[0,1\right]^2 }
\left\{
\left\langle \begin{bmatrix}x- \epsilon \\ 
-x - \epsilon \end{bmatrix} , v \right\rangle
\right\} 
\] 
\item {\bf Huber insensitive} loss, Fig.~\ref{fig:sel}:
\[
 \sup_{v \in [0,1]^2}
\left\{
\left\langle \begin{bmatrix}x- \epsilon \\ 
-x - \epsilon \end{bmatrix} , v \right\rangle
 - \frac{1}{2}v^Tv 
\right\}  
\]
\item {\bf Elastic net}, Fig~\ref{fig:enet}:
\[
\sup_{v \in [0,1] \times \mathbb{R}}
\left\{
\left\langle 
\begin{bmatrix}1  \\ 
1 \end{bmatrix}x , v \right\rangle
 - \frac{1}{2}v^T \begin{bmatrix} 0 & 0 \\ 0 & 1\end{bmatrix}v\right\}  
\]
\end{enumerate}

Note that the set $\calV$ is shown explicitly, and in each case can be easily represented as $\calV:=\{v\;:\;D^Tv \leq d\}$.
In addition, $H$ and $M$ are very sparse in all examples. 
\begin{flushright}
$\blacksquare$
\end{flushright}

Consider now optimizing a PLQ penalty subject to inequality constraints: 
\begin{equation}
\label{genELQP}
\begin{aligned}
\min_{x}\quad \rho( x) \\
\text{s.t.} \quad D^Tx \leq d
\end{aligned}.
\end{equation}
Using the techniques of convex duality theory developed in Appendix A2, 
the Lagrangian for \eqref{genELQP} is given by
\begin{eqnarray*}
\mathcal{L}(x, v, \omega) &=& \ip{\omega}{D^Tx-d}-\indicator{\omega}{\R_{+}^{n_1 }}+\ip{v}{b+Bx} \\
&-& \half v^TMv-\indicator{C^Tv-c}{\R_{-}^{n_2}},
\end{eqnarray*}
where $n_1$ and $n_2$ are dimensions of $d$ and $c$. 
The dual problem associated to this Lagrangian is 
\begin{equation}
\label{genELQP dual}
\begin{aligned}
\min_{(v,\omega)} &\quad \ip{d}{\omega}+\half v^TMv-\ip{b}{v}
\\
\mbox{s.t.} &\quad B^Tv+D\omega=0,\quad C^Tv \leq c,\quad 0\le \omega\ .
\end{aligned}
\end{equation}
The optimality conditions for this primal-dual pair are 
\begin{equation}
\label{eq:optimPLQ}
\begin{aligned}
&\omega,  w \geq 0\\
& D  \omega + B^T  v = 0 \\
&Mv+Cw=Bx+b\\
&C^Tv\le c\\
&D^Tx\le d\\
&\omega_j(D^Tx-d)_j=0,\ j=1,\dots,n_1\\
&w_j(C^Tv-c)_j=0,\ j=1,\dots,n_2.
\end{aligned}
\end{equation}
The final two conditions in \eqref{eq:optimPLQ} are called
{\it complementary slackness} conditions.
If $(\bx, \bv, \overline{\omega},\bw))$ satisfy all of the conditions in \eqref{eq:optimPLQ}, then $\bx$ solves the primal
problem \eqref{genELQP} and $(\bv, \overline{\omega})$ solves the dual problem \eqref{genELQP dual}.
The optimality criteria \eqref{eq:optimPLQ} are known as the Karush-Kuhn-Tucker (KKT) conditions
for \eqref{genELQP} and are used in the interior point method described in the next section. 

\subsection{Interior point (second-order) methods for PLQ functions}

Interior point methods directly target the KKT system~\eqref{eq:optimPLQ}.
In essence, they apply a damped Newton's method to a relaxed KKT system~\cite{KMNY91,NN94,Wright:1997}, recovering second-order behavior (i.e. 
superlinear convergence rates) for nonsmooth problems. 

To develop an interior point method for the previous section, we first introduce slack variables 
\[
s:= d-D^T  x\ge 0 \quad \mbox{and}  \quad r := c-C^T  v\ge 0\ .
\]
Complementarity slackness conditions~\eqref{eq:optimPLQ} can now be stated as 
$$
  \Omega S = 0\quad\mbox{and} \quad  WR = 0,
$$
where $ \Omega,S,  W,R$ are diagonal matrices with diagonals $ \omega,s,  w,r$, respectively. 
Let $\mathbf{1}$ denote the vector of all ones of the appropriate dimension.
Given $\mu>0$, we apply damped Newton iterations to the {\it relaxed} KKT system
$$
F_\mu(x, v, s, r, \omega, w) := 
\begin{bmatrix}
& D  \omega + B^T  v\\
& M  v +Cw- B x - b \\
& D^T x - d + s\\
& C^T v - c +r \\ 
& \Omega s - \mu \mathbf{1} \\
& Wr - \mu \mathbf{1}
\end{bmatrix}=0,
$$
where  $ \omega,s,  w,r \geq 0$ is enforced by the line search. 

Interior point methods apply damped Newton iterations to find a solution to $F_\mu = 0$
(with $ \omega,s,  w,r$ nonnegative) as $\mu$ is driven to $0$, so that cluster points are necessarily
KKT points of the original problem. Damped Newton iterations take the following form. 
Let $\xi := [x^T, v^T, s^T, r^T, \omega^T, w^T]^T$. Then the iterations are given by 
$$
\xi^{\iterc + 1} := \xi^\iterc - \gamma (F_{\mu_\iterc}^{(1)})^{-1} F_{\mu_\iterc},
$$
with $\gamma$ chosen so that $\omega^{\iterc+1}, w^{\iterc+1}, s^{\iterc+1}, r^{\iterc+1} \geq 0$ is satisfied, 
and some merit function (often $\|F_{\mu_\iterc}(\xi^{\iterc+1})\|$) is decreased. The homotopy parameter
${\mu_\iterc}$ is decreased at each iteration in a manner that preserves a measure of centrality within the 
feasible region.

While interior point methods have 
a long history (see e.g.~\cite{NN94,Wright:1997}), using them in this manner to solve any PLQ problem 
in a uniform way was proposed in~\cite{JMLR:v14:aravkin13a} to which we refer the
reader for further implementation details.
In particular, the Kalman smoothing case is fully developed in~\cite[Section 6]{JMLR:v14:aravkin13a}.
Each iteration of the resulting conjugate-PLQ interior point method can be implemented
with a complexity of $O(N(n^3 + m^3))$, which scales linearly in with $N$, just as for the classic smoother. 
The local convergence rate for IP methods is superlinear or quadratic in many circumstances~\cite{ye2011interior}, 
which in practice means that few iterations are required.


\section{Numerical experiments and illustrations}
\label{sec:numerics}

We now present a few numerical results to illustrate the formulations
and algorithms discussed above. 
In Section~\ref{sec:rates}, we consider a nonsmooth
Kalman formulation and compare the subgradient method, Chambolle-Pock, and 
interior point methods. In Section~\ref{sec:examples}, we show how nonsmooth formulations
can be used to address the motivating examples in the introduction. 
Finally, in Section~\ref{sec:generalPLQ}, we show how to construct
general piecewise linear quadratic Kalman smoothers (with constraints)
using the open-source package IPsolve.

\subsection{Algorithms and convergence rates}
\label{sec:rates}

In this section, we consider a particular signal tracking problem, where the 
underlying smooth signal is a sine wave, and a portion of the measurements 
are outliers. 

The synthetic ground truth function is given by $x(t) = \sin(-t)$. 
We reconstruct it from direct  noisy samples taken 
at instants multiple of $\Delta t$. 
We track this smooth signal 
by modeling it as an integrated Brownian motion which is equivalent to 
using cubic smoothing splines \cite{Wahba1990}. 
The state space model (sampled at instants where data are collected) is given by~\cite{Jaz,Oks,Bell2008}
$$
\begin{bmatrix}
\dot x_{t+1} \\ x_{t+1}
\end{bmatrix} 
= 
\begin{bmatrix}
1 & 0 
\\
\Delta t & 1
\end{bmatrix}
\begin{bmatrix}
\dot x_{t} \\ x_{t}
\end{bmatrix} 
+ v_t
$$
where the model covariance matrix of $v_t$ is 
$$
Q_t = \begin{bmatrix}\Delta t & \Delta t^2/2 \\ \Delta t^2 / 2 & \Delta t^3/ 3 \end{bmatrix}\;.
$$
The goal is to reconstruct the signal function from direct noisy measurements $y_t$, 
given by 
$$
y_t = C_t x_t + e_t\;, \quad C_t  = \begin{bmatrix} 0 & 1 \end{bmatrix}\;.
$$
We solve the following constrained modification of~\eqref{eq:KalmanL1}:
\begin{equation}
\label{eq:KalmanL1Con}
\min_{x\in \mathcal{C}} \| R^{-1/2}(y - Cx)\|_1 + \frac{1}{2} \|Q^{-1/2}(z - Ax )\|^2,
\end{equation}
where $z$ is constructed as in~\eqref{defs}.
For the sine wave, $\mathcal{C}$ is a simple bounding box, forcing each component to be in $[-1,1]$. 
Our goal is to compare three algorithms discussed in Section~\ref{sec:algo}: 
\begin{enumerate}
\item Projected subgradient method. We use the step size $\alpha_\iterc := \frac{1}{\iterc}$, 
and apply projected subgradient:
\[
x^{\iterc+1} := \proj_{\mathcal{C}}\left(x^\iterc - \frac{1}{\iterc}v^\iterc\right),
\]
where $v^\iterc \in \partial f(x^\iterc)$ is any element in the subgradient. 
\item  Chambolle-Pock (two variants described below).

\item Interior point formulation for~\eqref{genELQP}.
\end{enumerate}

Multiple splitting methods can be applied, including ADMM
(customized to deal with two nonsmooth terms), or the 
three-term splitting algorithm of~\cite{davis2015three}. 
We focus instead on a simple comparison 
of two variants of Chambolle-Pock with extremely different behaviors.  

To apply Chambolle-Pock, we first write the optimization problem~\eqref{eq:KalmanL1Con} 
using the template 
\[
\min_x f(Kx - r) + g(x).
\]
The Chambolle-Pock iterations (see Algorithm~\ref{alg:CP}) are given by
\[
\begin{aligned}
\omega^{\iterc+1} &:= r + \prox_{\sigma f^*} (\omega^\iterc + \sigma K (2 x^{\iterc} - x^{\iterc-1})-r)\\
x^{\iterc+1} &:= \proj_{\tau g} (x^\iterc - \tau K^T\omega^{\iterc+1}),
\end{aligned}
\]
where $\tau$ and $\sigma$ are stepsizes that must satisfy $\tau\sigma L <1$, and $L$ 
is the squared operator norm of $K$. Choices for $K$ give rise to different CP algorithms,  
and we two variants CP-V1 and CP-V2 below.

\noindent
{\bf CP-V1.} One way to make the assignment is as follows: 
\begin{eqnarray*}
f(\omega_1, \omega_2)  &=& \|\omega_1\|_1 + \frac{1}{2}\|\omega\|_2^2, \quad g(x) = \indicator{x}{\mathcal{C}}\\
f^*(\eta_1, \eta_2) &=& \indicator{\eta_1}{\mathbb{B}_{\infty}} + \frac{1}{2}\|\eta_2\|^2.\\
K &=& \begin{bmatrix} R^{-1/2} C\\ Q^{-1/2}A \end{bmatrix}, \quad r = \begin{bmatrix} R^{-1/2}y \\ Q^{-1/2} z\end{bmatrix}.
\end{eqnarray*}
The conjugate of $\|\cdot\|_1$ is computed in~\eqref{eq:ell1conj}, and it is easy to see that the 
function $\frac{1}{2}\|\cdot\|^2$ is its own conjugate using definition~\ref{eq:conjIntro}.

To understand the $\omega$-step, observe that
\[
\prox_{\sigma (f_1^*(x_1) + f_2^*(x_2)) }\!\left(\begin{bmatrix} y_1 \\ y_2\end{bmatrix}\right) 
\!=\! \begin{bmatrix} \prox_{\sigma f_1^*}(y_1)\\ \prox_{\sigma f_2^*}(y_2) \end{bmatrix}
\!=\! \begin{bmatrix} \proj_{ \mathbb{B}_{\infty}}(y_1)\\ \frac{1}{1+\sigma} \; y_2 \end{bmatrix}.
\]
The proximity operator for the indicator function is derived in~\eqref{eq:proxproj}, and the proximity
operator for $\frac{1}{2}\|\cdot\|^2$ is left as an exercise for the reader.  
The $x$-step requires a projection onto the set $\mathcal{C}$, which is the unit box for the sine example.  

\begin{figure*}
  \begin{center}
   \begin{tabular}{cccc}
\hspace{-.1in}
 { \includegraphics[scale=0.46]{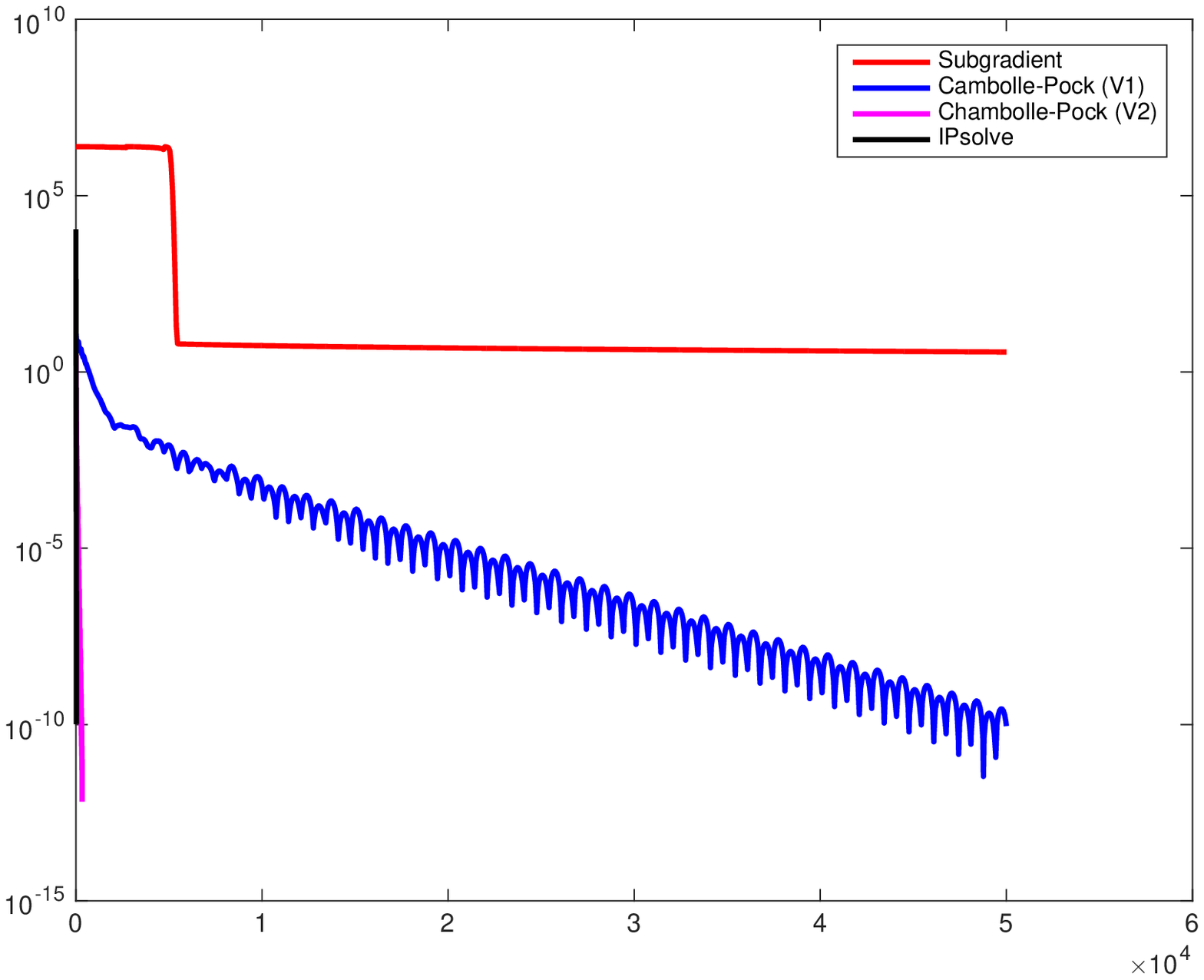}} 
\hspace{-.1in}
 { \includegraphics[scale=0.46]{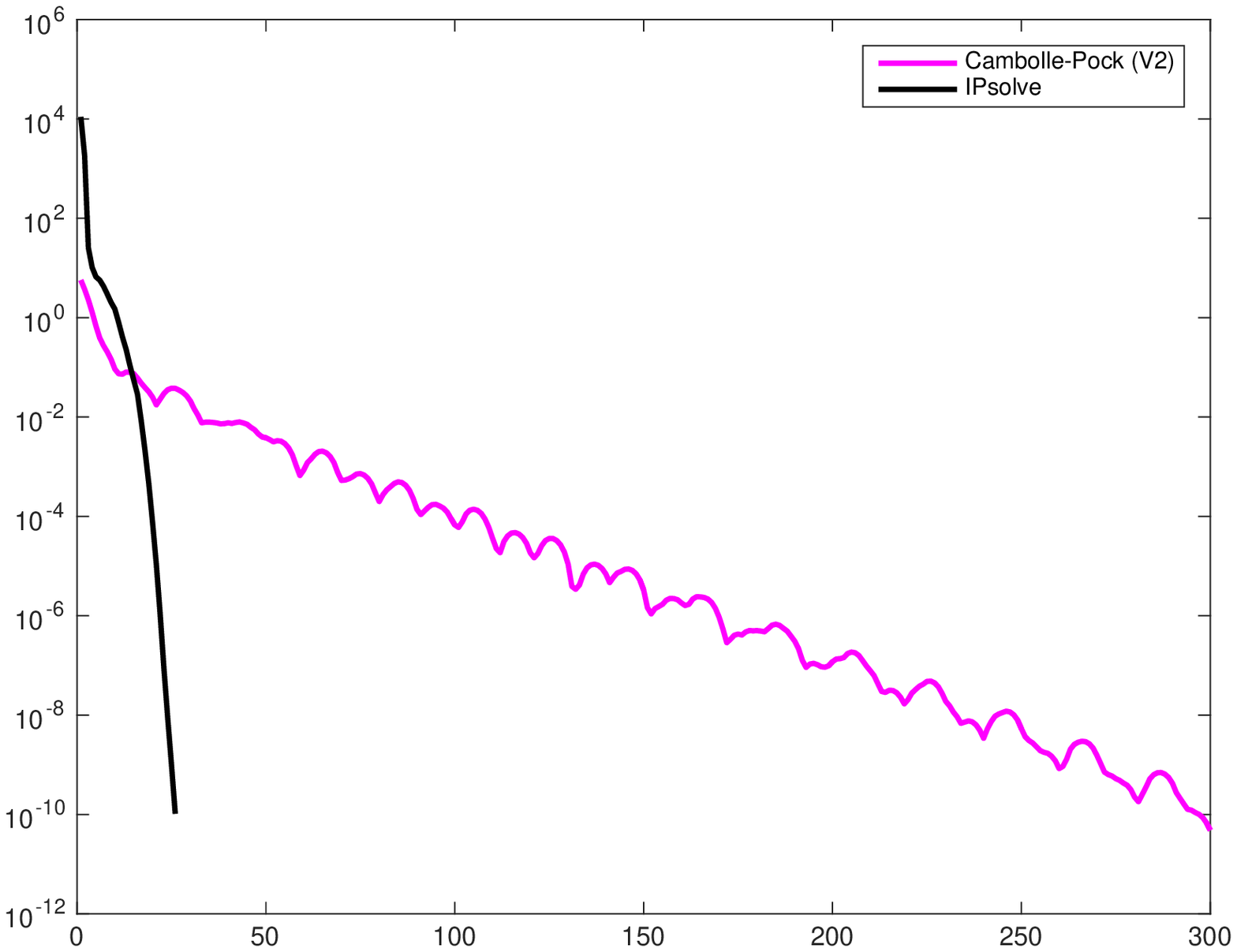}} 
    \end{tabular}
 \caption{{\bf{Convergence rate comparisons}}.  The $y$-axis 
 shows $f(x^t) - f(x^*)$, while $x$-axis shows the iteration count.
  {\it{Left:}} Convergence rates 
 for subgradient, CP-V1, CP-V2, and Interior Point methods, 
 after 50,000 iterations.   
  {\it{Right:}} Comparison for CP-V2 and IPsolve, after 300 iterations. Note that the methods have different complexities: subgradient and CP-V1 use only matrix vector products; CP-V2 requires a single factorization and then back-substitution at each iteration, and IPsolve solves linear systems at each iteration. } 
    \label{fig:rates}
     \end{center}
\end{figure*}

\noindent
{\bf CP-V2}. 
Here we treat $\frac{1}{2}\|Q^{-1/2}(Ax - z) \|^2$ as a unit, and assign in to $g$. As a result, the 
behavior of $A$ plays no role in the convergence rate of the algorithm.  
\begin{eqnarray*}
f(\omega_1, \omega_2) &=& \|\omega\|_1 + \indicator{\omega_2}{\mathcal{C}}, \quad g(x) = \frac{1}{2}\|Q^{-1/2}(Ax - z) \|^2\\ 
f^*(\eta_1,\eta_2) &=&  \indicator{\eta_1}{\mathbb{B}_{\infty}} + \|\eta_2\|_1.\\ 
K &=& \begin{bmatrix} R^{-1/2} C \\ I\end{bmatrix}, \quad r = \begin{bmatrix} R^{-1/2}y \\ 0\end{bmatrix}.
\end{eqnarray*}

The proximity operator for $g$ is obtained by solving a linear system:
\[
\prox_{\tau g}(y) = (\tau A^T Q^{-1}A + I)^{-1}(y + \tau A^TQ^{-1}z). 
\]
The linear system  $\tau A^T Q^{-1}A + I$ is block tridiagonal positive definite, and its eigenvalues are 
bounded away from $0$. Since it does not change between iterations, we compute its Cholesky factorization
once and use it to implement the inversion at each iteration. This requires a single factorization using $O(n^3N)$
arithmetic operations, followed by multiple $O(n^2N)$ iterations (same cost as matrix-vector
products with a block tridiagonal system). 

The $\omega$-step for CP-V2 is also different from the $\omega$-step in CP-V1, but still very simple and efficient:  
\begin{eqnarray*}
\prox_{\sigma (f_1^*(x_1) + f_2^*(x_2)) }\left(\begin{bmatrix} y_1 \\ y_2\end{bmatrix}\right) 
&=& \begin{bmatrix} \prox_{\sigma f_1^*}(y_1)\\ \prox_{\sigma f_2^*}(y_2) \end{bmatrix} \\
&=& \begin{bmatrix} \proj_{ \mathbb{B}_{\infty}}(y_1)\\ \prox_{\sigma \|\cdot\|_1} (y_2) \end{bmatrix}.
\end{eqnarray*}
The proximity operator for $\sigma \|\cdot\|_1$ is derived in~\eqref{eq:proxell1}.

The results are shown in Fig.~\ref{fig:rates}. The subgradient method is disastrously slow, and difficult 
to use. Given a simple step size schedule, e.g. $\alpha_\iterc = \frac{1}{\iterc}$, it may waste tens of thousands
of iterations before the objective starts to decrease. In the left panel of Fig.~\ref{fig:rates}, 
it took over 10,000 iterations before any noticeable impact. Moreover, as the step sizes become small, 
it can stagnate, and while in theory it should continue to slowly improve the objective, in practice it 
stalls on the example problem.\\  
CP-V1 is able to make some progress, but the results are not impressive. 
Even though the algorithm requires only matrix-vector products, it is adversely impacted by the conditioning 
of the problem. In particular, the ODE term for the Kalman smoothing problem (i.e. the $A$) can be
poorly conditioned, and in the CP-V1 scheme, it sits inside $K$. As a result, we see very slow 
convergence. Interestingly, the rate itself looks linear, but the constants are terrible, so it requires 50,000 iterations to fully solve the problem.\\ 
In contrast, CP-V2 performs extremely well. The algorithm treats the quadratic ODE term 
as a unit, and the ill-conditioning of $A$ does not impact the convergence rate. The price we pay is
having to solve a linear system at each iteration. However, since the system does not change, we factorize 
it once, at a cost of $O(n^3N)$, and then use back-substitution to implement $\prox_g$ at each iteration. 
The resulting empirical convergence rate is also linear, but with a significant 
improvement in the constant: CP-V2  needs only 300 iterations to reach $10^{-10}$ accuracy 
(gap to the minimum objective value), see the right plot of Fig.~\ref{fig:rates}.\\
Finally, IPsolve has a {\it super-linear} rate, and finishes in 27 iterations. 
It is not possible to pre-factorize any linear systems, so the complexity is $O(n^3N)$
for each iteration. For moderate problem sizes 
(specifically, smaller $n$), this approach is fast and generalizes to any PLQ losses 
$V$ and $J$ and any constraints.  
For large problem sizes, CP-V2 will win; however, it is very specific to the current problem. 
In particular, if we change $J$ in~\eqref{fullGen} 
from the quadratic to the 1-norm or Huber, we would need to develop a different splitting approach.  
The more general CP-V1 approach is far less effective. 

The following sections focus on modeling and the resulting behavior of the estimates. 
Section~\ref{sec:examples} presents the results for the motivating examples 
in the introduction.  

\begin{figure*}
  \begin{center}
   \begin{tabular}{cccc}
\hspace{.1in}
 { \includegraphics[scale=0.46]{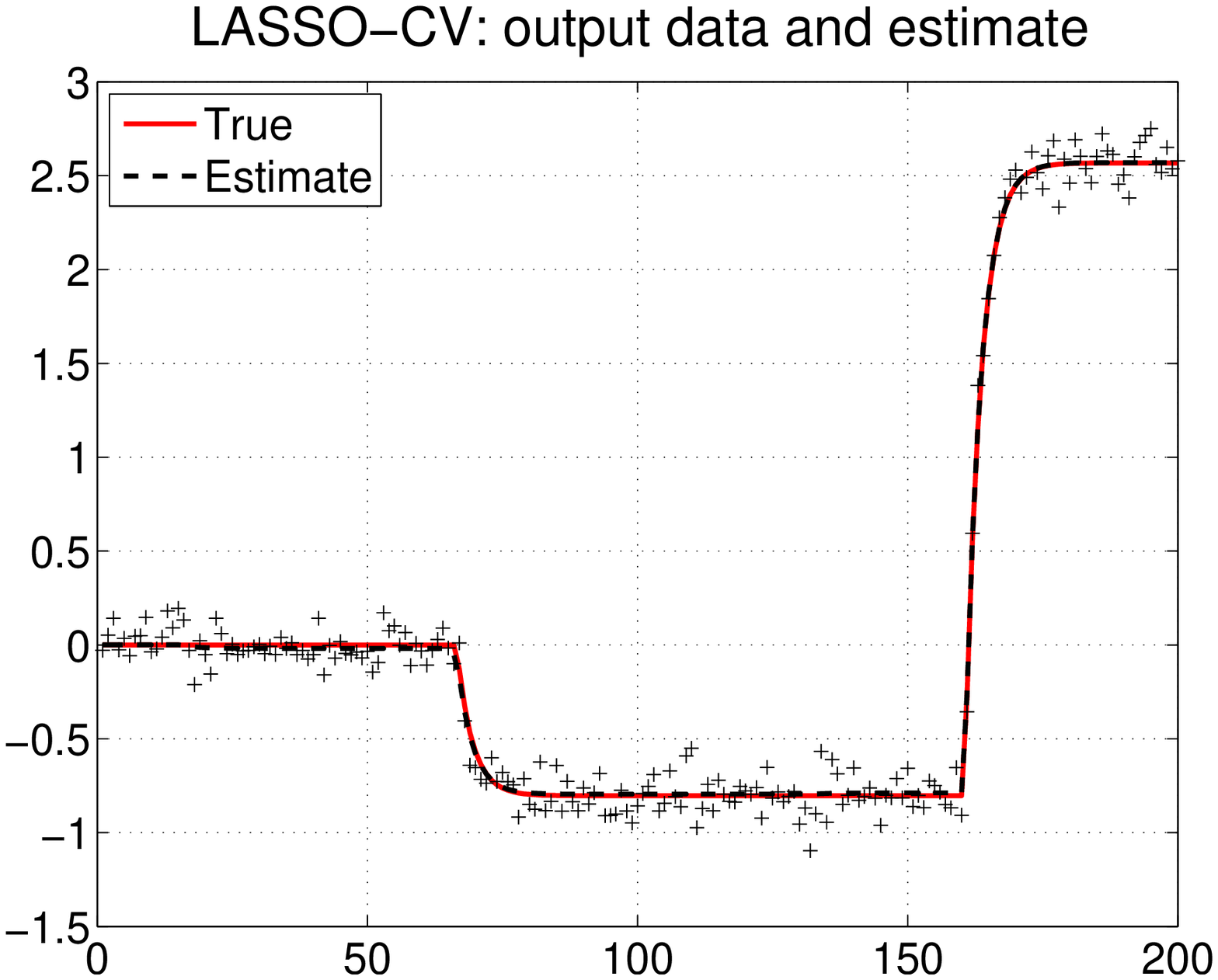}} 
\hspace{.1in}
 { \includegraphics[scale=0.46]{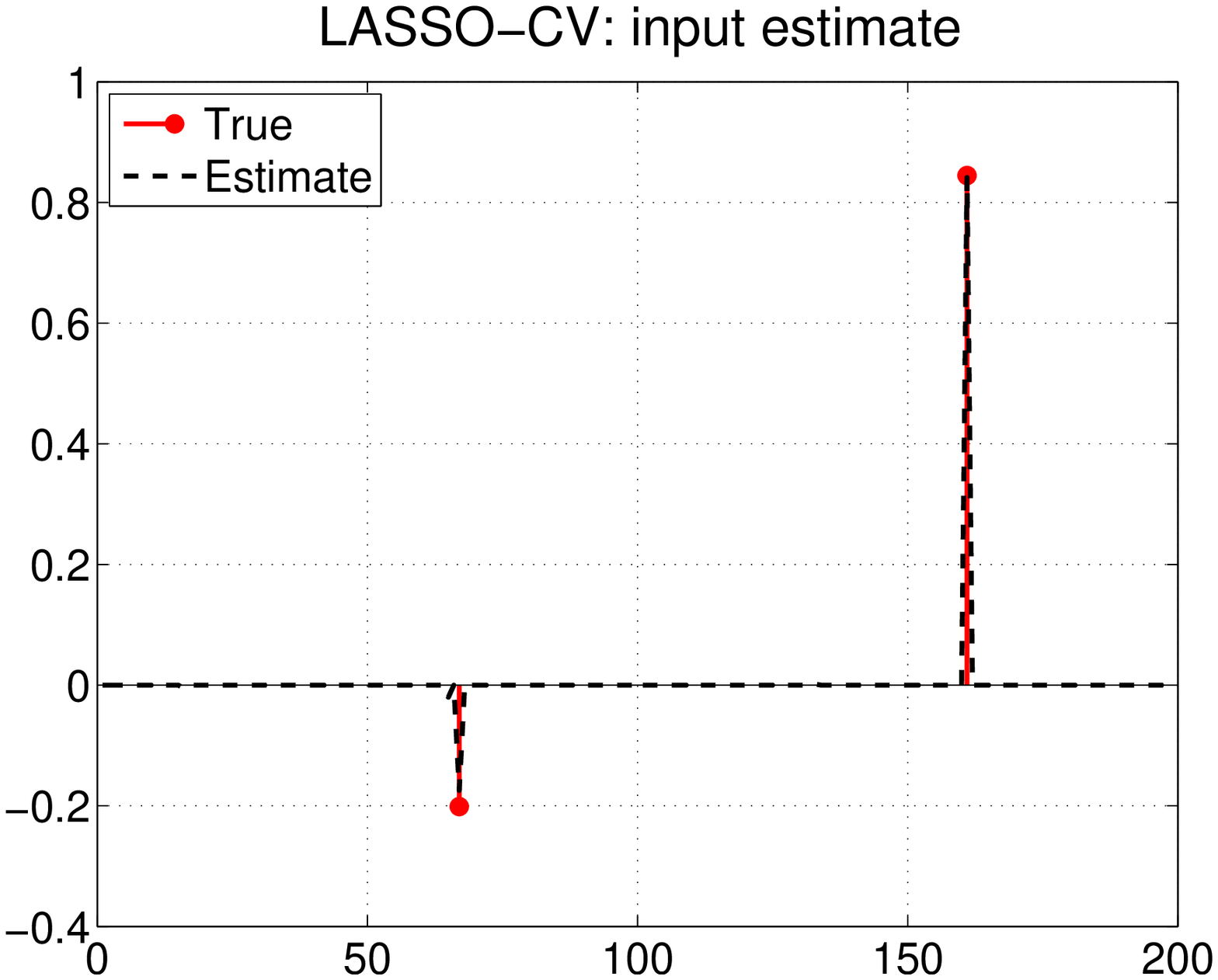}} 
    \end{tabular}
 \caption{{\bf{DC motor and impulsive disturbances}}. {\it{Left:}} noiseless output (solid line), measurements  ($+$) and output reconstruction by the LASSO smoother (dashed line). 
 {\it{Right:}} impulsive disturbance and reconstruction by the LASSO smoother (dashed line). } 
    \label{Fig34Imp}
     \end{center}
\end{figure*}

\subsection{DC motor: robust solutions using $\ell_1$ losses and penalties}
\label{sec:examples}

We now solve the problems described in subsection \ref{MotEx}
using two different smoothing formulations based on the $\ell_1$ norm.

{\bf{Impulsive inputs:}} Let $E_1=\left(\begin{array}{cc}1 & 0\end{array}\right), \ E_2=\left(\begin{array}{cc}0 & 1\end{array}\right).$
To reconstruct the disturbance torque $d_t$
acting on the motor shaft, we use the LASSO-type estimator 
proposed in \cite{Ohlsson2011}:
\begin{equation}
\label{logMAPlasso}
\begin{aligned}
& \min_{x_1,\ldots,x_N}  
\sum_{t=1}^N   \left (y_t - E_2 x_t\right)^2 
+ \gamma \sum_{t=0}^{N-1}  | d_t |\\
&\text{subject to the dynamics}  \  (\ref{DCmotor}) 
\end{aligned}
\end{equation}
Since $u_t=0$, this corresponds to the optimization problem 

$$
\begin{aligned}
 &\min_{x_1,\ldots,x_N}  &&
\sum_{t=1}^N \!  (y_t \!-\! E_2 x_t )^2  \\ 
& &&\!+\! \frac{\gamma}{2} \left[ \!\sum_{t=0}^{N-1}\!  \frac{| E_1 ( x_{t+1} \!-\! A_t x_t  )|}{11.81}   
\!+\! \frac{ | E_2 \left( x_{t+1} \!-\! A_t x_t \right)|}{0.625}\!    \right] \\
&\text{subject to}  &&  \frac{E_1 \left( x_{t+1} - A_t x_t  \right)}{11.81}  = \frac{E_2 \left( x_{t+1} - A_t x_t  \right)}{0.625} 
\end{aligned}
$$

The regularization parameter $\gamma$ is tuned using 5-fold cross validation on a grid consisting of 20 values, logarithmically  spaced between 0.1 and 10. The resulting smoother is dubbed LASSO-CV.

\begin{figure}
  \begin{center}
   \begin{tabular}{cccc}
\hspace{.1in}
 { \includegraphics[scale=0.45]{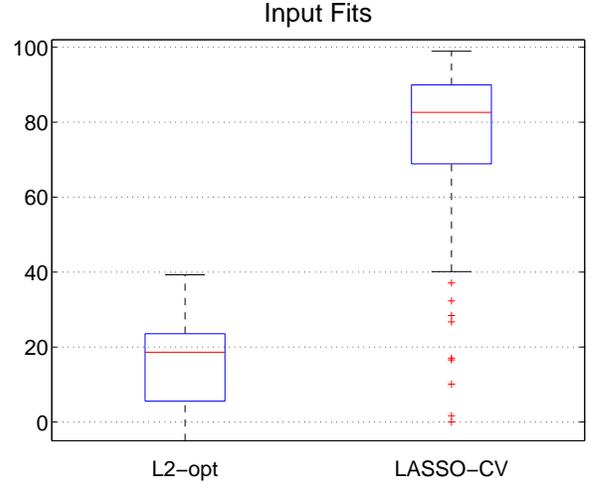}} 
    \end{tabular}
 \caption{{\bf{DC motor and impulsive disturbances}}. Boxplot of the fits returned by 
 optimal linear smoother (left) and by the LASSO smoother (right).} 
    \label{Fig5Imp}
     \end{center}
\end{figure}

The right panel of Fig. \ref{Fig34Imp} shows the estimate of 
$d_t$ obtained by LASSO-CV 
starting from the noisy outputs in the left panel. 
Note that we recover the impulsive disturbance, and that the
LASSO smoother outperforms the optimal linear smoother L$_2$-opt, 
shown in Fig. \ref{Fig12Imp}. To further exam the improved performance of the LASSO smoother
in this setting, 
we performed a Monte Carlo study of 200 runs,
comparing the fit measure
$$
100\left(1- \frac{\|\hat{d} - d \|}{\| d \|}\right),
$$
where $d=[d_1 \ldots d_{200}]$ is the true signal and $\hat{d}$ is the estimate
returned by  L$_2$-opt or by LASSO-CV. 
Fig. \ref{Fig5Imp} shows Matlab
boxplots of the 200 fits obtained by these estimators.
The rectangle contains the inter-quartile range ($25-75\%$ percentiles)  
of the fits, with median shown by the red line. 
The ``whiskers" outside the rectangle display the upper and
lower bounds of all the numbers,
not counting what are deemed outliers,
plotted separately as ``+".  
The effectiveness of the LASSO smoother is clearly supported by this study.

\begin{figure*}
  \begin{center}
   \begin{tabular}{cccc}
\hspace{.1in}
 { \includegraphics[scale=0.46]{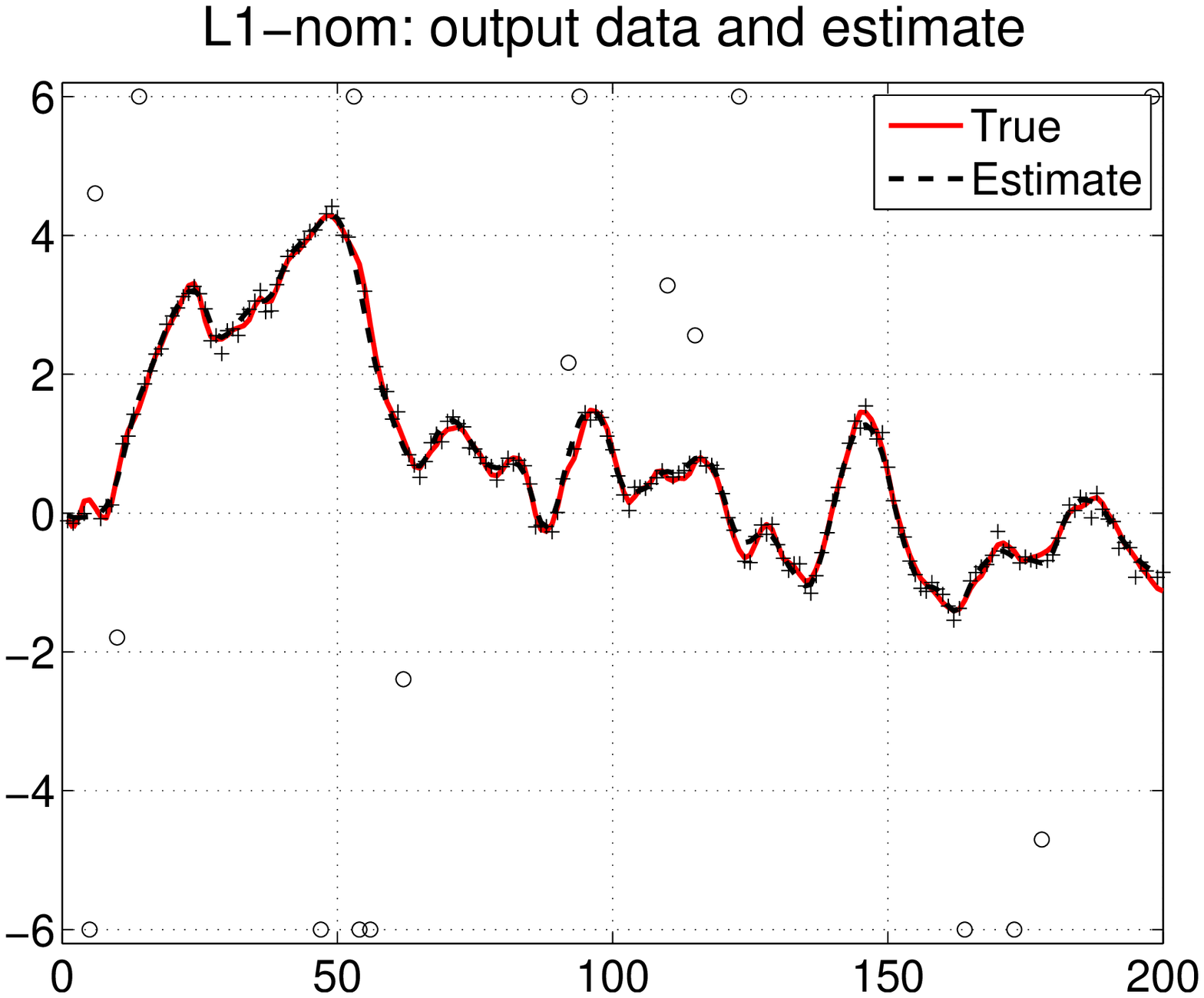}} 
\hspace{.1in}
 { \includegraphics[scale=0.46]{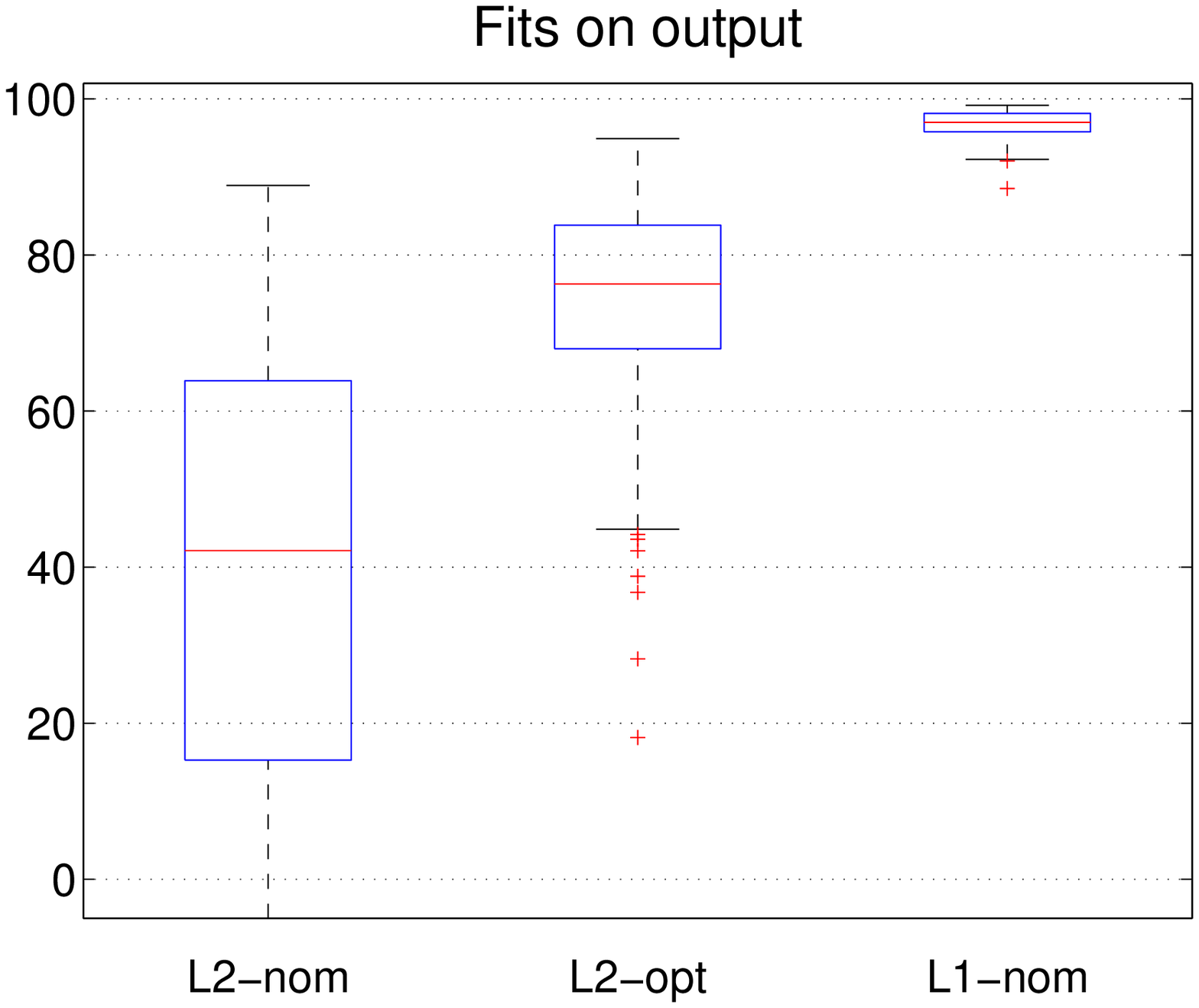}} 
    \end{tabular}
  \caption{{\bf{DC motor and outliers in the output measurements}}. {\it{Left:}} noiseless output (solid line), measurements  ($+$), outliers
   ($\circ$) and output reconstruction by the robust smoother equipped with the $\ell_1$ loss (dashed). {\it{Right:}} boxplot of the output fits returned by the nominal and
 optimal linear smoothers and by the robust smoother (both $L_2$- and $L_1$-nom use the nominal standard deviation $\sigma=0.1$ as weight for the residuals).} 
    \label{FigOut34}
     \end{center}
\end{figure*}

\begin{figure}
  \begin{center}
   \begin{tabular}{cccc}
\hspace{.1in}
 { \includegraphics[scale=0.46]{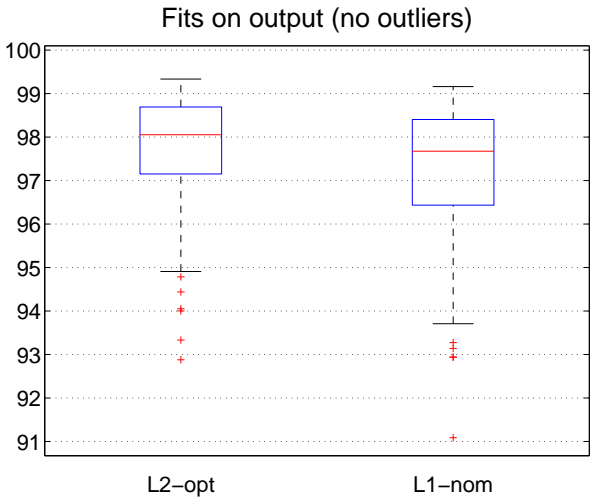}} 
    \end{tabular}
  \caption{{\bf{DC motor and output reconstruction without outliers corrupting the measurements}}. Boxplot of the output fits returned by the
 optimal linear smoother and by the robust smoother.} 
    \label{FigOut5}
     \end{center}
\end{figure}

{\bf{Presence of outliers:}} To reconstruct the angle velocity, we use the following
smoother based on the $\ell_1$ loss:
\begin{equation}
\label{logMAPlasso2}
\begin{aligned}
& \min_{x_1,\ldots,x_N}  
\sum_{t=1}^N  \frac{ | y_t - E_2 x_t | }{\sigma} + \frac{1}{0.1^2}  \sum_{t=0}^{N-1}  d_t^2\\
&\text{subject to the dynamics}  \  (\ref{DCmotor}) 
\end{aligned}
\end{equation}
Recall that $d_t \sim \G(0,0.1^2)$, so  now there is no impulsive input.
The $\ell_1$ loss used in  (\ref{logMAPlasso2})  
is shown in Fig. \ref{fig:1norm}. It can also be viewed as a limiting
case of Huber (Fig. ~\ref{fig:Huber}) and Vapnik (Fig.~\ref{fig:Vapnik}) losses, respectively,
when their breakpoints $\kappa$ and $\epsilon$ are set to zero.\\  
Over the state space domain, problem (\ref{logMAPlasso2}) is equivalent to 
$$
\begin{aligned}
& \min_{x_1,\ldots,x_N}  \quad
\sum_{t=1}^N  \frac{ | y_t - E_2 x_t | }{\sigma} \\
& + \frac{1}{0.1^2} 
\left[ \sum_{t=0}^{N-1}  \frac{ \left (E_1 ( x_{t+1} - A_t x_t ) \right)^2}{11.81}  
+  \frac{ \left (E_2 ( x_{t+1} - A_t x_t ) \right)^2}{0.625}   \right] \\
&\text{subject to}  \  
\frac{E_1 \left( x_{t+1} - A_t x_t  \right)}{11.81}  = \frac{E_2 \left( x_{t+1} - A_t x_t  \right)}{0.625}.
\end{aligned}
$$
Note that the $\ell_1$ loss uses the nominal standard deviation $\sigma=0.1$ as weight for the residuals, so that
we call this estimator L$_1$-nom.\\
The left panel of Fig. \ref{FigOut34} displays the estimate of the angle 
returned by L$_1$-nom. The profile is very close to truth, revealing the 
robustness of the smoother to the outliers.
Here, we have also performed a Monte Carlo study of 200 runs, 
using the fit measure
$$
100\left(1- \frac{\|\hat{y} - y \|}{\| y \|}\right), 
$$
where $y=[y_1 \ldots y_{200}]$ is the true value while $\hat{y}$ are the estimates
returned by L$_2$-nom, L$_2$-opt or L$_1$-nom. 
The boxplots in the right panel of Fig. \ref{FigOut34} compare 
the fits of the three estimators, and illustrate the
robustness of L$_1$-nom.\\
Finally, we repeated the same Monte Carlo study setting
$\alpha=0$, generating no outliers in the output measurements.
Under these assumptions, L$_2$-nom and L$_2$-opt 
coincide and represent the best estimator among all the possible smoothers. 
Fig. \ref{FigOut5}  shows Matlab
boxplots of the 200 fits obtained by L$_2$-nom and L$_1$-nom.
Remarkably, the robust smoother has nearly identical performance to 
the optimal smoother, so there is little loss of performance under nominal conditions. 

\subsection{Modeling with PLQ using IPsolve}
\label{sec:generalPLQ}

\begin{figure*}
  \begin{center}
   \begin{tabular}{cccc}
\hspace{.1in}
 { \includegraphics[scale=0.43]{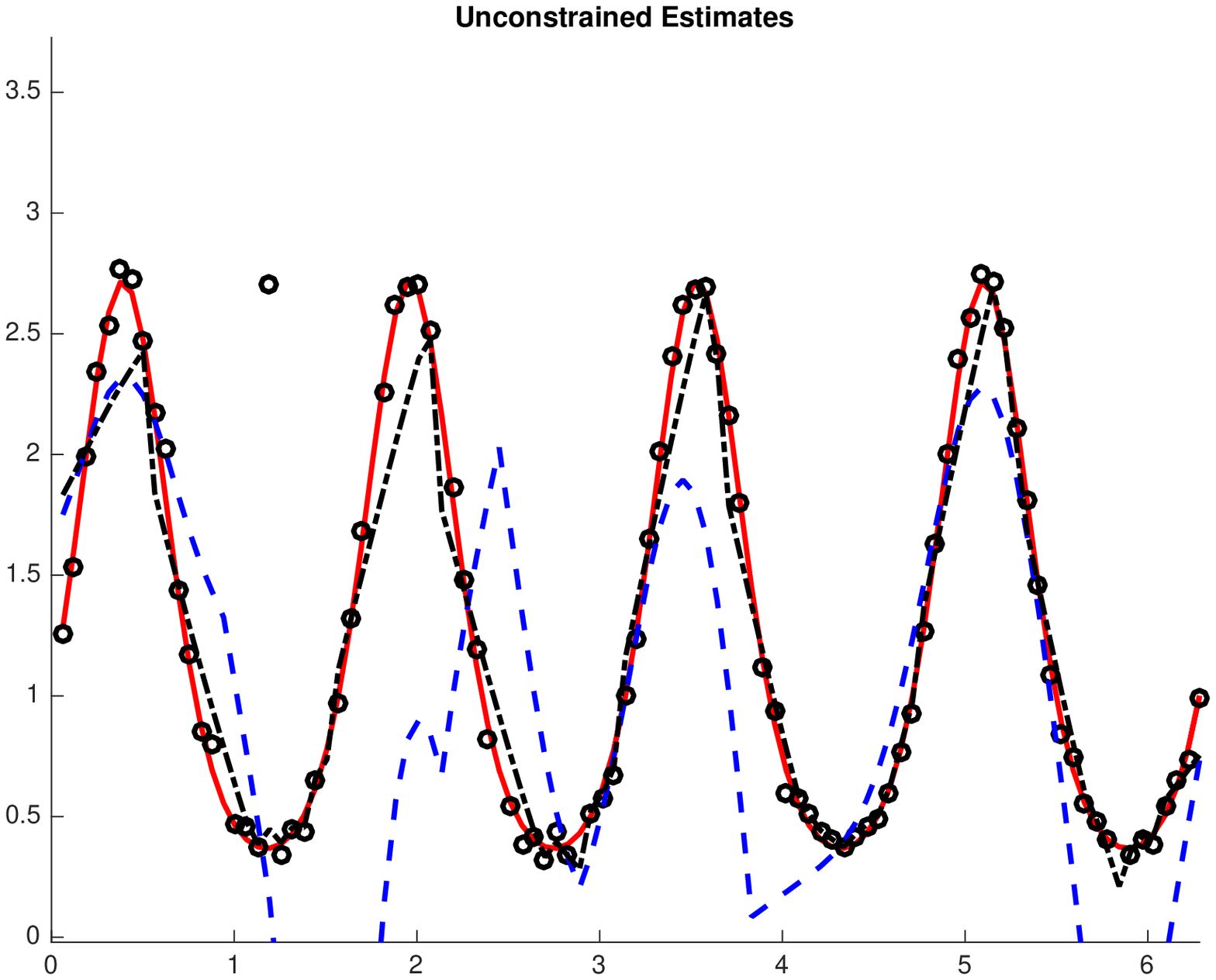}} 
\hspace{.1in}
 { \includegraphics[scale=0.43]{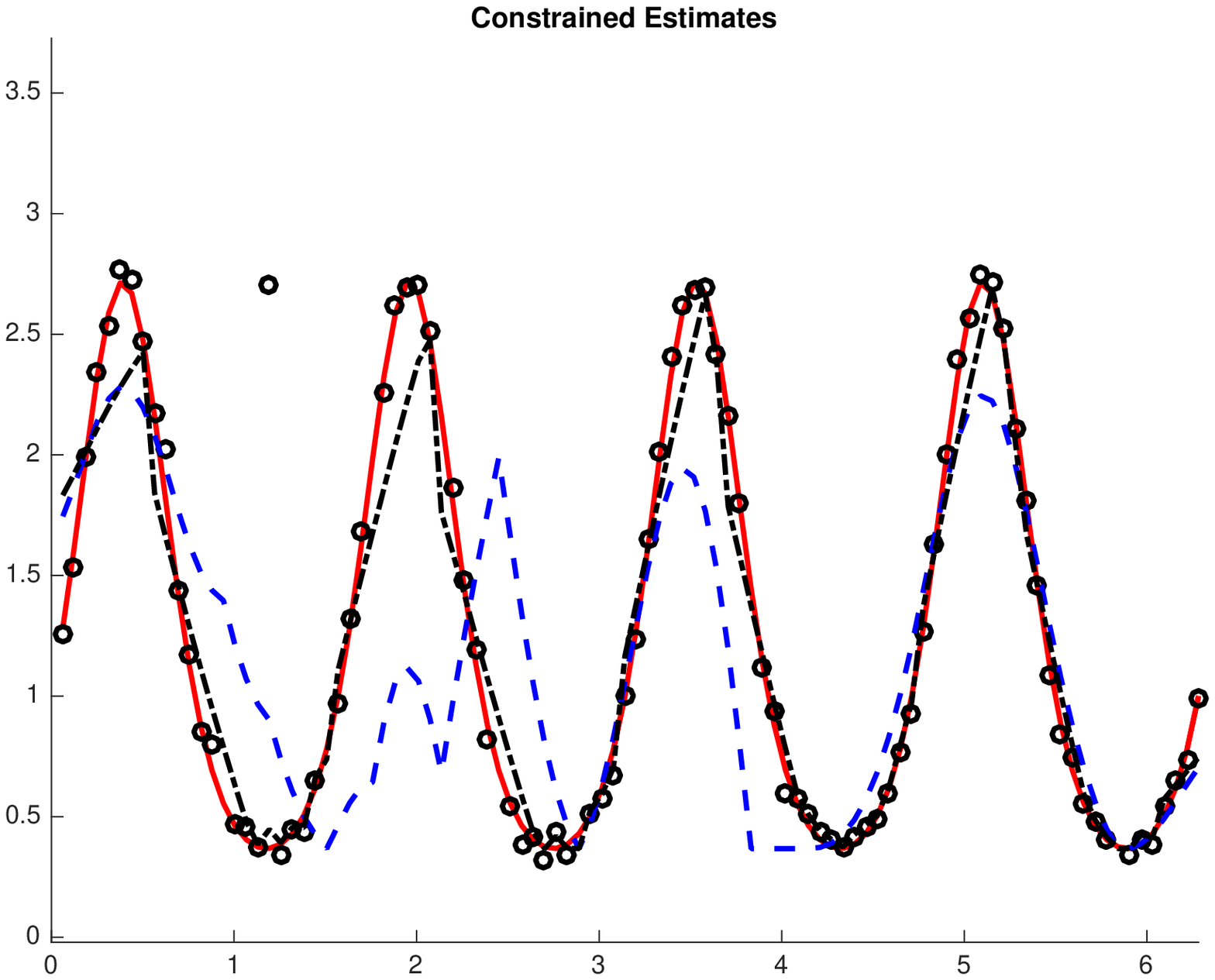}} 
    \end{tabular}
 \caption{{\bf{Results of four smoothers}}.  
  {\it{Left:}} Ground truth (solid red) and unconstrained results for L2 (dashed blue) and Huber (densely dashed black).    
  {\it{Right:}} Ground truth (solid red) and constrained results for cL2 (dashed blue) and cHuber (densely dashed black). 
  Constraints can be very  helpful in dealing with contamination. Best results are obtained when we use both robust penalties
  and constraints on the domain.} 
    \label{fig:PLQmodels}
     \end{center}
\end{figure*}

In this section, we include several modeling examples that combine robust
penalties with constraints. Each example is implemented using IPsolve.
The solver and 
examples are available online at 
\scriptsize \texttt{https://github.com/saravkin/IPsolve}, \normalsize
see in particular \scriptsize \texttt{blob/master/515Examples/KalmanDemo.m} 
\normalsize inside the folder \verb{IPsolve{.\\
In all examples, 
the ground truth function of interest is given by $x(t) = \exp(\sin(4t))$, 
and we reconstruct it from direct and noisy samples taken 
at instants multiple of $\Delta t$. 
The function $x(t)$ is smooth and periodic,
but the exponential accelerates the transitions around the maximum and minimum values.
The process and measurement models are the same as in Section~\ref{sec:rates}.
\noindent Four smoothers~\eqref{fullGen} are compared in this example using IPsolve. 
The L$_2$ smoother uses the quadratic penalty for both $V$ and $J$, and no constraints. 
The cL$_2$ smoother uses least squares penalties with constraints
including the information that $\exp(-1) \leq x(t) \leq \exp(1) \ \forall t$. 
The Huber smoother uses Huber penalties ($\kappa =1$) for both $V$ and $J$, without constraints,
while cHuber uses Huber penalties ($\kappa = 1$) together with constraints.
The results are shown in Fig.~\ref{fig:PLQmodels}.  
90\% of the measurement errors are generated from a Gaussian with nominal standard deviation $0.05$, 
while 10\% of the data are large outliers generated using a Gaussian with standard deviation $10$. 
The smoother is given the nominal standard deviation. 

The least squares smoother L$_2$ without constraints does a very
poor job. The Huber smoother
obtains a much better fit. Interestingly, cL$_2$ is much better than L$_2$, 
indicating that domain constraints 
can help a lot, even when using quadratic penalties. Combining constraints and robustness in cHuber 
gives the best fit since the inclusion of constraints eliminates the constraint
violations of Huber at 3 and 6 seconds in the left plot of Fig.~\ref{fig:PLQmodels}.\\

The calls to IPsolve are given below: 
\scriptsize
\begin{enumerate}
\item L$_2$: 
\begin{verbatim}
params.K = Gmat; params.k = w;
L2 = run_example( Hmat, meas, 'l2', 'l2', ...
[], params );
\end{verbatim}

\item Huber: 
\begin{verbatim}
params.K = Gmat; params.k = w;
Huber = run_example( Hmat, meas, 'huber', ...
'huber', [], params );
\end{verbatim}
The only difference required to run the HH smoother is to replace the names of the PLQ penalties 
in the calling sequence. 

\item cL$_2$: 
\begin{verbatim}
params.K = Gmat; params.k = w;
params.constraints = 1; conA = [0 1; 0 -1]; 
cona = [exp(1); -exp(-1)];
params.A = kron(speye(N), conA)'; 
params.a = kron(ones(N,1), cona);
cL2 = run_example( Hmat, meas, 'l2', 'l2',...
 [], params );
\end{verbatim}
For constraints, we need to create the constraint matrix and also pass it in 
using the \verb{params{ structure. 

\item cHuber:
\begin{verbatim}
params.K = Gmat; params.k = w;
params.constraints = 1; conA = [0 1; 0 -1]; 
cona = [exp(1); -exp(-1)];
params.A = kron(speye(N), conA)'; 
params.a = kron(ones(N,1), cona);
cHuber = run_example( Hmat, meas, 'huber',...
 'huber',[], params );
\end{verbatim}
The constrained Huber call sequence requires only name change for the PLQ penalties. 

\end{enumerate}
\normalsize

Above, one can see that the names of PLQ measurements are arguments to the file \verb{run_example{, 
which builds the combined PLQ model object that it passes to the interior point method. 
The measurement matrix and observations vector are also 
passed directly to the solver. The process terms are passed through the auxiliary 
\verb{params{ structure. Full details for constructing the matrices are provided in the online demo
\verb{KalmanDemo{ already cited above.

\section{Concluding remarks}

Various aspects of the state estimation problem in the linear system (1) have
been treated over many years in a very extensive literature. One reason for
the richness of the literature is the need to handle a variety of realistic
situations to characterize the signals $v$ and $e$ in (1). This has led to deviations
from the classical situation with Gaussian signals where the estimation problem
is a linear-quadratic optimization problem.
This survey attempts to give a comprehensive and systematic
treatment of the main issues in this large literature. The key has been to start
with a general formulation (15) that contains the various situations as special
cases of the functions $V$ and $J$. An important feature is that (15) still is a
convex optimization problem under mild and natural assumptions. This opens
the huge area of convex optimization as a fruitful arena for state
estimation. In a way, this alienates the topic from the original playground
of Gaussian estimation techniques and linear algebraic solutions. The survey
can therefore also be read as a tutorial on convex optimization techniques being
applied to state estimation.


{\bf Appendix}

\section*{A1. Optimization viewpoint on Kalman smoothing under correlated noise and singular covariances}

In some applications, the noises $\{e_t,v_t\}_{t=1}^N$ are correlated.
Assume that $e_t$ and $v_t$ are still jointly Gaussian, but 
with a cross-covariance denoted by $S_t$.
For $t=1,\ldots,N$, this implies that the last assumption in (\ref{eq:wgn}) can be replaced by
$$
\begin{aligned}
E(e_t v_s^\top)  = 
\left\{ \begin{array}{cl}
    S_t & \mbox{   if} ~ t=s \\
    0 & \mbox{   otherwise,}
\end{array} \right.
\end{aligned}
$$
while $v_0$ is assumed independent of $\{e_t,v_t\}_{t=1}^N$.\\
We now reformulate the objective (\ref{eq:opt}) under this more general
model. Define the process $\tilde{v}_0=v_0$ and 
$$
\tilde{v}_t = v_t - E(v_t | e_t) = v_t - S_t R_t^{-1} e_t, \quad  t \geq 1
$$
which, by basic properties of Gaussian estimation, is
independent of $e_t$ and consists of white noise with covariance
$$
\tilde{Q}_t = Q_t - S_t R_t^{-1} S_t^\top, \quad  t \geq 1.   
$$
Since $v_t$ is correlated only with $e_t$, we have that all
the $\{\tilde{v}_t\}$ and $\{e_t\}$ form a set of mutually independent Gaussian noises.  
Also, since $e_t=y_t-C_t x_t$, model (\ref{eq:Lin}) can be reformulated as
\begin{subequations}  \label{eq:LinPseudo}
\begin{align}
 x_{t+1} &=\tilde{A}_t x_t+B_t u_t + S_t R_t^{-1} y_t + \tilde{v}_t\\
y_t &=C_t x_t + e_t
\end{align}
\end{subequations}
where we define $\tilde{A}_0 x_0 + S_0 R_0^{-1} y_0 = A_0 x_0$ while 
$$
\tilde{A}_t = A_t - S_t R_t^{-1}  C_t, \quad  t \geq 1. 
$$
Note that (\ref{eq:LinPseudo}) has the same form as the original system (\ref{eq:Lin})
except for the presence of an additional input 
given by the output injection $S_t R_t^{-1} y_t$.\\
Assuming also the initial condition $x_0$ independent of the noises,
the joint density of $\{\tilde{v}_t\},\{e_t\}$ and $x_0$ turns out
\begin{equation*}
\Bp\left(x_0,\{e_t\},\{\tilde{v}_t\} \right) = 
\Bp\left(x_0 \right) \prod_{t=1}^N \Bp_{e_t}  \left( e_t  \right) \prod_{t=0}^{N-1} \Bp_{\tilde{v}_t} \left(\tilde{v}_t  \right), 
\end{equation*}
where we use $\Bp_{e_t}$ and $\Bp_{\tilde{v}_t}$ to denote the densities corresponding to $e_t$ and $\tilde{v}_t$.
Since $\{x_t\}_{t=0}^N$ and $\{y_t\}_{t=1}^N$ are a linear transformation of $\{v_t\}_{t=0}^N$, $\{e_t\}_{t=1}^N$ and $x_0$,
the joint posterior of states and outputs is proportional to
\begin{equation*}
\Bp\left(x_0 \right) \prod_{t=1}^N \Bp_{e_t}  \left( y_t - C_t x_t  \right) \prod_{t=0}^{N-1} \Bp_{\tilde{v}_t} \left( x_{t+1} -\tilde{A}_t x_t - S_t R_t^{-1} y_t - B_t u_t   \right).
\end{equation*} 
Consequently, maximizing the posterior of the states given the output measurements 
is equivalent to solving
\begin{equation}
  \label{eq:optPseudo}
\begin{aligned}
\min_{x_0,\ldots,x_N} &\| \Pi^{-1/2} (x_0 - \mu ) \|^2 +\sum_{t=1}^N\|R_t^{-1/2}( y_t-C_t x_t)\|^2 & \\
& + \sum_{t=0}^{N-1}\|  \tilde{Q}_t^{-1/2}(x_{t+1} -\tilde{A}_t x_t - S_t R_t^{-1} y_t - B_t u_t)\|^2. 
\end{aligned}
\end{equation}
Next consider the case where some of the covariance matrices are singular.
If some of the matrices $Q_t$ or $R_t$ are not invertible, problems
(\ref{eq:optPseudo}) and (\ref{eq:opt}) are not well-defined. 
In this case, one can proceed as follows. First, $\tilde{v}_t,\tilde{Q}_t$ and $\tilde{A}_t $
can be defined in the same way where $R_t^{-1}$ is replaced by its pseudoinverse 
$R_t^{\dag}$. The objective can then be reformulated by
replacing $\tilde{Q}_t^{-1}$ and $R_t^{-1}$ 
by $\tilde{Q}_t^{\dag}$ and $R_t^{\dag}$, respectively.  Linear constraints can be added to
prevent the state evolution in the null space of $\tilde{Q}_t$ and $R_t$. 
By letting $I_Q$ and $I_R$ be the sets with the time instants associated with singular 
$\tilde{Q}_t$ and $R_t$, problem (\ref{eq:optPseudo}) can be rewritten as
\begin{equation}
  \label{eq:optPseudo2}
\begin{aligned}
\min_{x_0,\ldots,x_N} &\| \Pi^{-1/2} (x_0 - \mu ) \|^2 +\sum_{t=1}^N\|(R_t^{\dag})^{1/2}( y_t-C_t x_t)\|^2 & \\
& + \sum_{t=0}^{N-1}\| (\tilde{Q}_t^{\dag})^{1/2}(x_{t+1} -\tilde{A}_t x_t - S_t R_t^{-1} y_t - B_t u_t)\|^2\\
&\text{subject to} \;R_t^{\perp} \left(y_t - C_t x_t \right) = 0 \  \text{for} \ t \in I_R\ \mbox{ and}\\ 
& \tilde{Q}_t^{\perp} \left(x_{t+1} - A_t x_t - S_t R_t^{\dag} y_t - B_t u_t \right) = 0  \  \text{for} \  t \in I_Q ,
\end{aligned}
\end{equation}

where $R_t^{\perp} = I-R_t R_t^{\dag}$ and $\tilde{Q}_t^{\perp} = I - \tilde{Q}_t \tilde{Q}_t^{\dag}$ provide the projections onto the null-space of $R_t$ and $\tilde{Q}_t$,
respectively.

\section*{A2. Convex analysis and optimization}\label{sec:opt}

Some of the background in convex analysis and optimization used
in the previous sections is briefly reviewed in this section. 
In particular, the fundamentals used in the development 
and analysis of algorithms for~\eqref{fullGen} is reviewed.\\ 
Many members of the broader class of penalties~\eqref{fullGen}
do not yield least squares objectives  
since they include nonsmooth penalties and constraints; however, they are convex. 
Convexity is a fundamental notion in optimization theory and practice and gives access to globally
optimal solutions as well as extremely efficient and reliable numerical solution techniques that scale to
high dimensions. The relationship between convex sets and functions was presented in Section~\ref{sec:convexBasics}.

\subsection*{Fundamental objects in convex analysis}

We begin by developing a duality theory for the 
general objective~\eqref{fullGen}. This is key for both algorithm design and sensitivity analysis.  
Duality is a consequence of the separation theory for convex sets. 

{\bf Separation:} We say that a hyperplane (i.e. an affine set
of co-dimension 1) separates two sets if they lie on opposite sides of the hyperplane.
To make this 
idea precise, we introduce the notion of {\it relative interior}.
The affine hull of a set $\calE$, denoted $\aff{\calE}$, is the intersection of all affine sets that contain $\calE$.

Given 
$\calE\subset\R^n$ the relative interior of $\calE$ is 
$$
\ri{\calE}:=
\set{x\in \calE}{\exists\eps>0\mbox{ s.t. }
(x+\eps\bB)\cap\aff{\calE}\subset \calE}.
$$
For example, $\mathrm{ri}\set{(2,x)}{-1\le x\le 1}=\set{(2,x)}{-1<x<1}$.

Let $\cl \calE$ denote the closure of set $\calE$, and $\intr \calE$ denote the interior. 
Then the boundary of $\calE$ is given by $\bdry{\calE} := \cl{\calE} \setminus \intr{\calE}$,
and the relative boundary $\rbdry{\calC}$ is given by $\cl{\calC}\setminus\ri{\calC}$. 

\begin{theorem}[Separation]\label{thm:separation}
Let $\calC\subset\R^n$ be nonempty and convex, and suppose $\bar y\notin\ri{\calC}$.
Then there exist $z\ne 0$ such that
$$
\ip{z}{\bar y}>\ip{z}{y}\quad\forall\ y\in \ri{\calC}.
$$
\end{theorem}

{\bf Support Function:} Apply Theorem \ref{thm:separation} to a point
$\bar x\in\rbdry{\calC}$ to obtain a nonzero vector $z$ for which
\begin{equation}\label{eq:separation2}
\ip{z}{\bx}=\support{z}{\calC}:=\sup\set{\ip{z}{x}}{x\in \calC} > \inf\set{\ip{z}{x}}{x\in \calC}.
\end{equation}
The function $\sig_\calC$ is called
the {\it support function} for $\calC$, and the nonzero vector $z$ is said to be
a support vector to $\calC$ at $\bx$.
When $\calC$ is polyhedral, $\sig_\calC$ is an example
of a  PLQ function, with~\eqref{eq:separation2} a special case of~\eqref{PLQpenalty} 
with $M =0$.

{\bf Example: dual norms}. Given a norm $\norm{\cdot}$ on $\Rn$ with unit ball $\mathbb{B}$,
the dual norm is given by
\[
\dnorm{z}:=\sup_{\norm{x}\le 1}\ip{z}{x} = \support{z}{\mathbb{B}}.
\]
For example, the 2-norm is self dual, while the dual norm for $\|\cdot\|_1$ is $\|\cdot\|_\infty$.

This definition implies that 
$\norm{x}=\sig_{\bB^\circ}(x)$, where 
\[ \bB^\circ:=\set{z}{\ip{z}{x}\le 1\ \forall\, x\in\bB}.\]
The set $\bB^\circ$ is the closed unit ball for the dual norm $\dnorm{\cdot}$. 
This kind of relationship between the unit ball of a norm and that of its
dual generalizes to {\it polars} of sets and cones. %

{\bf Polars of sets and cones:} For any set $\calC$
in $\Rn$, the set
\[ \calC^\circ:=\set{z}{\ip{z}{x}\le 1\ \forall\, x\in \calC}\]
is called the \emph{polar} of $\calC$, and we have $(\calC^\circ)^\circ = \cl{\conv{\calC\cup\{0\}}}$.  
Hence, if $\calC$ is a closed convex set containing the origin, then $(\calC^\circ)^\circ =\calC$.
If $\calK\subset\Rn$ is a convex cone 
($\calK$ is a convex and $\lam \calK\subset \calK$ for all $\lam >0$), 
then, by rescaling,  
\[
\calK^\circ=\set{z}{\ip{z}{x}\le 0\ \forall\, x\in \calK}\ \mbox{ and }\ (\calK^\circ)^\circ=\cl{\calK}.
\]
In particular, this implies that 
$\sig_\calK=\delta_{\calK^\circ}$. 

{\bf Subdifferential:} For nonsmooth convex functions, the notion of derivative can be captured 
by examining support vectors to their epigraph. 
Define the domain of the function $f$ to be the set
$\dom{f}:=\set{x}{f(x)<\infty}$.  
Using the fact that 
$$
\ri{\epi{f}}=\set{(x,\mu)}{x\in\ri{\dom{f}}\mbox{ and }f(x)<\mu},
$$
Theorem \ref{thm:separation} tells us  that, for every $\bx\in\ri{\dom{f}}$,
there is a support vector to $\epi{f}$ at $(\bx,f(\bx))$ of the form $(z,-1)$, which separates 
the points in the epigraph from the points in a half space below the epigraph: 
\[ 
\ip{(z,-1)}{(\bx,f(\bx))}\ge \ip{(z,-1)}{(x,f(x))}\quad\forall\ x\in\dom{f},
\] 
or equivalently,
\begin{equation}\label{eq:sd1}
f(\bx)+\ip{z}{x-\bx}\le f(x)\quad\forall\ x\in\dom{f}.
\end{equation}
This is called the {\it subgradient inequality}.
The vectors $z$ satisfying \eqref{eq:sd1} are said to be subgradients of $f$ at $\bx$, and
the set of all such subgradients is called the {\it subdifferential} of $f$ at $\bx$, denoted $\partial f(\bx)$.
This derivation shows that $\partial f(\bx)\ne \emptyset$ for all $\bx\in\ri{\dom{f}}$ when $f$ is {\it proper}, 
i.e. $\dom f$ is nonempty, with $f(x)> -\infty$.
In addition, it can be shown that $\partial f(\bx)$ is a singleton if and only if $f$ is differentiable at $\bx$
with the gradient equal to the unique subgradient.\\ 
For example, the absolute value function on $\R$ is not differentiable
at zero so there is no tangent line to its graph at zero; however, every line passing through the
origin having slope between $-1$ and $1$ defines a support vector to the epigraph
at the origin. In this case, we can replace the notion of derivative by the set of slopes of hyperplanes at the origin.  
Each of these slopes is a subgradient, and the set of all these is the {\it subdifferential} of $|\cdot|$ at the origin. 

{\bf Necessary and Sufficient Conditions for Optimality:} An immediate consequence of the subgradient
inequality is that 
\[
0\in\partial f(\bx) \quad \text{ if and only if} \quad \bx\in\argmin f.
\] 
That is, a first-order necessary
and sufficient condition for optimality in convex optimization is that the zero vector is an element
of the subdifferential. Returning to the absolute value function on $\R$,
note that the zero slope hyperplane supports the epigraph at zero and zero is the global minimizer of $\abs{\cdot}$. 

\begin{theorem}[Convex Optimality]\label{thm:cvx optimality}
Let $\map{f}{\Rn}{\R\cup\{+\infty\}}$ be a closed proper convex function. 
Then the following conditions are equivalent:
\begin{enumerate}
\item[(i)] $\bx$ is a global solution to the problem $\min_x f$.
\item[(ii)] $\bx$ is a local solution to the problem $\min_x f$.
\item[(iii)] $0\in\partial f(\bx)$.
\end{enumerate}
\end{theorem}

{\bf Convex conjugate:} 
Again consider the support functions defined in~\eqref{eq:separation2}. 
By construction, $z\in\partial f(x)$ if and only if
\[
\ip{(z,\!-1)}{(x,f(x))}\!=\!\support{(z,1)}{\epi{f}}\!=\!\sup_y(\ip{z}{y}\!-f(y)) \!=\! f^*(z),
\] 
or equivalently, $f(x)+f^*(z)=\ip{z}{x}$. 
When $f$ is a proper convex function, the conjugate function $f^*$ (defined in~\eqref{eq:conjIntro}), is a closed,
proper, convex function, since it is the pointwise supremum of the affine functions $z\to\ip{z}{y}-f(y)$
over the index set $\dom{f}$. Consequently we have  
\[
\partial f(x)=\set{z}{f(x)+f^*(z)\le\ip{z}{x}}.
\] 
Due to the symmetry of this expression for
the subdifferential, it can be shown that $(f^*)^*=f$ and $\partial f^* =(\partial f)^{-1}$
(i.e. $z\in\partial f(x)\iff x\in\partial f^*(z)$)
whenever $f$ is a closed proper convex function.
These relationships guide us to focus on the class of functions
\[
\Gam_n:=\set{\map{f}{\R^n}{\R\cup\{\infty\}}}{\mbox{$f$ is closed proper and convex}}.
\]

For example, if $\calC\subset\R^n$ is a nonempty closed convex set, then 
$\del_{\calC}\in \Gam_n$, where $\del_{\calC}$ is defined in \eqref{eq:indicator}.
It is easily seen that $\del_\calC^*=\sig_\calC$ and, for $x\in \calC$,  
\[
\partial \indicator{x}{\calC} = \set{z}{\ip{z}{y-x}\le 0\ \forall\, y\in \calC}=:\ncone{x}{\calC},
\]
where
\(
\ncone{x}{\calC}
\)
is called the normal cone to $\calC$ at $x$. 

{\bf Calculus for PLQ}: 
Just as in the smooth case, subdifferentials and conjugates become useful in practice by developing a calculus 
for their ready computation.
Here we focus on calculus rules for PLQ functions $\rho$ defined in \eqref{PLQpenalty} 
which are well established in \cite{RTRW}. In particular, if we set $q(v):=\half v^TMv+\indicator{v}{\calV}$,
then, by \cite[Corollary 11.33]{RTRW}, 
either $\rho \equiv \infty$ or
\begin{equation}\label{eq:PLQ conj}
\rho^*(y)=
\inf_{B^Tv=y}\left[q(v)-\ip{b}{v}\right]\ \mbox{ and }\  
\partial \rho (z)=
B^T\partial q^*(Bz+b),
\end{equation}
which can be reformulated as
$$\partial \rho (z)=\set{B^Tv}{v\in \calV\ \mbox{ and }\ Bz-Qv+b\in\ncone{v}{\calV}}.$$
In addition, we have from \cite[Theorem 3]{JMLR:v14:aravkin13a} that
\begin{equation}\label{eq:dom PLQ}
\dom{\rho^*}=B^T\calV \ \mbox{ and } \ 
\dom{\rho}=B^{-1}\left(\polar{[\calV^\infty\cap\Nul{M}]}-b\right),
\end{equation}
where
$\calV^\infty$ is the {\it horizon cone} of $\calV$. 
As the name suggests, $\calV^\infty$ is a closed cone, and, when  
$\calV$ is nonempty convex, it is a nonempty closed convex cone satisfying
$\calV^\infty=\set{w}{\calV+w\subset \calV}$. In particular, $\calV$ is bounded if and only if $\calV^\infty = \{0\}$. 

The reader can verify by inspection of figs.~\ref{fig:quadratic}-\ref{fig:enet} that the domain of each scalar PLQ is $\mathbb{R}$. This is also immediate from~\eqref{eq:dom PLQ}. 
Four of the six penalties have bounded sets $\calV$, so that $\calV^\infty = \{0\}$, 
the polar is the range of $B$, and so the result follows immediately. 
The quadratic penalty has $\calV^\infty = \mathbb{R}$, but $\Nul{M} = \{0\}$. We leave the elastic net as an exercise.\\  
More importantly,~\eqref{eq:PLQ conj} gives explicit expressions for derivatives and subgradients of PLQ functions in terms of $v$. Consider the Huber function, fig.~\ref{fig:Huber}. From~\eqref{eq:PLQ conj}, we have
\[
\partial \rho(z)=\set{v}{v\in \kappa[-1,1]\ \mbox{ and }\ z-v\in\ncone{v}{\kappa[-1,1]}}.
\]
From this description, we immediately have $\partial\rho(z) = \nabla \rho(z) = z$ for 
$|z| <\kappa$, and $\kappa\, \mbox{sgn}(z)$ for $|z|>\kappa$. 

\subsection*{Convex duality}

There are many approaches for convex duality theory \cite{RTRW}. 
For our purposes, we choose one based on the convex-composite Lagrangian \cite{Burke1985}.

{\bf Primal objective:} 
Let $f\in\Gam_{m},\ g\in \Gam_{n}$, and $K\in\R^{m\times n}$
and consider the primal convex optimization
problem
\begin{equation}
\label{eq:plqPrimal}
\frakP\qquad\min_x \sfp(x):= f(Kx)+g(x),
\end{equation}
where we call $\sfp(x)$ the {\it primal} objective.

The structure of the problem~\eqref{eq:plqPrimal} is the same as that used to develop the celebrated
Fenchel-Rockafellar Duality Theorem \cite[Section 31]{RTR} (Theorem \ref{F-R duality theorem} below). 
It is sufficiently general to allow an easy translation
to several formulations of the problem~\eqref{fullGen}  
depending on how one wishes to construct an algorithmic framework. This variability in
formulation is briefly alluded to
in Section~\ref{sec:splitting}. In this section, we focus on general
duality results for~\eqref{eq:plqPrimal} leaving the discussion of specific reformulation 
of~\eqref{fullGen}
to the discussion of algorithms.

We now construct the {\it dual} to the convex optimization problem $\frakP$. 
In general,
the dual is a concave optimization problem, but, as we show, it is often beneficial
to represent it as a
convex optimization problem.

{\bf Lagrangian:} First, define the 
{\it Lagrangian} $\map{\mathcal{L}}{\R^n\times\R^{m}\times\R^{n}}{\R\cup\{-\infty\}}$ 
for $\frakP$ by setting 
\[
\mathcal{L}(x,w, v):=\ip{w}{Kx}-f^*(w)+\ip{v}{x}-g^*(v).
\]
The definition of the conjugate immediately tells us that the primal objective is given by
maximizing the Lagrangian over the dual variables:
\[f(Kx)+g(x)=\sup_{w,v}\mathcal{L}(x,w,v).\]
{\bf Dual objective}: Conversely, the dual objective is is obtained by minimizing the Lagrangian
over the primal variables:
\[
\sfd(w,v):=\inf_x\, \mathcal{L}(x,w,v)=
\begin{cases}-f^*(w)
-g^*(v),& K^Tw+v=0,\\ 
-\infty,&K^Tw+v\ne 0.\end{cases}
\]
The corresponding dual optimization problem is
\[
\max_{w,v}\, \sfd(w,v) \quad= \quad
\max_{K^Tw+v=0}-f^*(w)
-g^*(v).
\]
One can eliminate $v$ from the dual problem and reverse sign to obtain 
a simplified version of the dual problem:
\begin{equation}
\label{eq:plqDual}
\frakD\qquad\min_{w}\,  
\tilde \sfd(w):=f^*(w)
+g^*(-K^Tw).
\end{equation}
Three examples of primal-dual problems pairs are given in Table 1.

{\bf Weak and strong duality:} By definition, $\max \sfd(w,v)\le\min \sfp(x)$,
or equivalently, $0\le (\min \tilde \sfd(w)) +  (\min \sfp(x))$. This inequality is called {\it weak duality}. 
If equality holds,
we say the duality gap is zero. If solutions to both $\frakP$ and $\frakD$ exist with zero duality gap,
then we say {\it strong duality} holds. In general, a zero duality gap and strong duality require
additional hypotheses called {\it constraint qualifications}. Constraint qualifications
for the problem $\frakP$ are given as conditions (a) and (b) in the following theorem.

\begin{table*} 
\label{ex pd pairs}
\hspace{2cm}  \begin{tabular}{|c|c|c|l|l|}
\hline 
&$\ \stackrel{f}{g}$&$\ \stackrel{f^*}{g^*}$&$\frakP$&$\frakD$\\ \hline 
Basis&
$\indicator{\cdot-s}{\tau\bB_2}$&$\tau\tnorm{\cdot}+\ip{w}{s}$&$\min\ \onorm{x}$&$\min\ \tau\tnorm{w}+\ip{w}{s}$\\
Pursuit~\cite{BergFriedlander:2008}&
$\onorm{\cdot}$&$\indicator{\cdot}{\bB_\infty}$&$\mathrm{s.t. }\tnorm{Kx-s}\le \tau$&$\mbox{s.t. }\inorm{K^Tw}\le 1$\\ &&&&\\
\hline
LASSO&
$\half\tnorm{\cdot -s}^2$&$\ip{\cdot}{s}+\half\tnorm{\cdot}^2$&$\min\half\tnorm{Kx-s}^2$&
$\min \half\tnorm{w}^2+\kappa\inorm{K^Tw}+\ip{w}{s}$\\
&$\indicator{\cdot}{\kappa\bB_1}$&$\kappa\inorm{\cdot}$&$\mbox{s.t. }\onorm{x}\le \kappa$& \\ 
&&&&\\
\hline
Lagrangian&
$\half\tnorm{\cdot -s}^2$&$\ip{\cdot}{s}+\half\tnorm{\cdot}^2$&$\min \half\tnorm{Kx-s}^2+\lam \onorm{x} $
&$\min\half\tnorm{w+s}^2-\half\tnorm{s}^2$
\\
&$\lam\onorm{\cdot}$&$\indicator{\cdot}{\lam\bB_\infty}$&&$\mbox{s.t. }\inorm{K^Tw}\le \lam$\\ 
\hline 
\end{tabular}
\caption{We show three common variants of sparsity promoting formulations, and compute the dual is in each case using 
the relationships between~\eqref{eq:plqPrimal} and~\eqref{eq:plqDual}. Strong duality holds for all three examples. }
\end{table*}

\begin{theorem}[Fenchel-Rockafellar Duality Theorem]\cite[Corollary 31.2.1]{RTR}
\label{F-R duality theorem}
Let $f\in\Gam_m,\ g\in \Gam_n$, and $K\in\R^{m\times n}$. If either
\begin{enumerate}
\item[(a)] there exists $x\in \ri{\dom{g}}$ with $Kx\in\ri{\dom{f}}$, or
\item[(b)] there exists $w\in\ri{\dom{f^*}}$ with $-K^Tw\in\ri{\dom{g^*}}$,
\end{enumerate}
hold, 
then $\min \sfp\, +\,\min\tilde\sfd=0$ with finite optimal values. Under condition (a), 
$\argmin \tilde\sfd$ is nonempty, while under (b), $\argmin \sfp$ is nonempty.
In particular, if both (a) and (b) hold, then strong duality between $\frakP$ and $\frakD$ holds in the sense
that $\min \sfp\, +\,\min\tilde\sfd=0$ with finite optimal values that are attained in both $\frakP$ and $\frakD$.
In this case, optimal solutions are characterized by
$$
\left\{
\begin{aligned}
\bx\mbox{  solves  }\frakP\\ \bw\mbox{  solves  }\frakD\\ \min \sfp\, +\,\min\tilde\sfd=0
\end{aligned}
\right\}
\iff
\left\{
\begin{aligned}
\bw\in \partial f(K\bx)\\ -K^T\bw\in\partial g(\bx)
\end{aligned}
\right\}
$$
$$
\iff
\left\{
\begin{aligned}
\bx\in\partial g^*(-K^T\bw)\\ K\bx\in \partial f^*(\bw)
\end{aligned}
\right\}.
$$
\end{theorem}

\bibliographystyle{abbrv}
\bibliography{kalmanSurvey}

\end{document}